\documentclass[11pt, reqno]{amsart}

\usepackage[top=3.75cm, bottom=3cm, left=3cm, right=3cm]{geometry}
\usepackage{amsmath}
\usepackage{amsthm}
\usepackage{amsfonts}
\usepackage{amssymb}
\usepackage{mathtools}
\usepackage{bm}
\usepackage[colorlinks=true, citecolor=darkblue, urlcolor=darkblue, linkcolor=darkblue, linktocpage = true]{hyperref}
\usepackage{color}
\usepackage[all,cmtip]{xy}
\usepackage{enumitem}
\usepackage{tikz}  
\usepackage{graphicx}
\usepackage{scalerel}
\usepackage{comment}
\usepackage[mathscr]{eucal}
\usepackage{stmaryrd}
\usepackage{tikz-cd}
\usepackage{appendix}
\usepackage{adjustbox}

\definecolor{darkblue}{rgb}{0,0,.75}

\numberwithin{equation}{section}

\theoremstyle{plain}
\newtheorem{theorem}{Theorem}[section]
\newtheorem{lemma}[theorem]{Lemma}
\newtheorem{proposition}[theorem]{Proposition}
\newtheorem{corollary}[theorem]{Corollary}

\newtheorem{question}[theorem]{Question}

\theoremstyle{definition}
\DeclareMathOperator{\Def}{Def}

\newtheorem{definition}[theorem]{Definition}
\newtheorem{example}[theorem]{Example}
\newtheorem{remark}[theorem]{Remark}

\newtheorem{construction}[theorem]{Construction}

\makeatletter
\newtheoremstyle{italicsname}
 {3pt}
 {3pt}
 {\itshape}
 {}
{\bf}
 {.}
 {.5em}
 {\thmname{#1}\thmnumber{\@ifnotempty{#1}{ }#2}%
 \thmnote{ {\the\thm@notefont(#3)}}}
\makeatother
\theoremstyle{italicsname}

\setlist[itemize]{topsep=5pt,itemsep=3pt}
\setlist[enumerate]{topsep=5pt,itemsep=3pt}

\newcommand{\Perf}{\mathrm{Perf}}

\newcommand{\Dperf}{\mathrm{D}_{\mathrm{perf}}}
\newcommand{\Dqc}{\mathrm{D}_{\mathrm{qc}}}

\DeclareMathOperator{\QCoh}{QCoh}

\DeclareMathOperator{\CatDef}{CatDef}
\DeclareMathOperator{\CRing}{CRing}
\DeclareMathOperator{\dGeomSt}{GeomdSt}
\DeclareMathOperator{\fppf}{fppf}

\DeclareMathOperator{\perf}{perf}

\DeclareMathOperator{\colim}{colim}
\DeclareMathOperator{\gerby}{gerby}

\newcommand{\Mod}{\mathrm{Mod}}

\DeclareMathOperator{\sqz}{sqz}
\DeclareMathOperator{\PrSt}{Pr^{St}}

\newcommand{\Fun}{\mathrm{Fun}}
\newcommand{\Ind}{\mathrm{Ind}}

\newcommand{\coarse}{\mathbf{CS}}

\DeclareMathOperator{\Grp}{Grp}
\newcommand{\dSt}{\mathrm{dSt}}

\newcommand{\dArtSt}{\mathrm{ArtdSt}}
\DeclareMathOperator{\cocart}{cocart}
\DeclareMathOperator{\DAlg}{DAlg}
\DeclareMathOperator{\QAlg}{QAlg}
\DeclareMathOperator{\DAlgMod}{DAlgMod}
\DeclareMathOperator{\CAlgMod}{CAlgMod}
\DeclareMathOperator{\StMod}{StMod}
\DeclareMathOperator{\Der}{Der}
\newcommand{\bAt}{\mathbf{At}}
\DeclareMathOperator{\cl}{cl}

\newcommand{\Gr}{\mathrm{Gr}}

\DeclareMathOperator{\At}{At}

\DeclareMathOperator{\Spec}{Spec}

\newcommand{\op}{\mathrm{op}}

\DeclareMathOperator{\RHom}{RHom}
\DeclareMathOperator{\Hom}{Hom}

\DeclareMathOperator{\Aut}{Aut}

\DeclareMathOperator{\HH}{HH}

\DeclareMathOperator{\Map}{Map}

\DeclareMathOperator{\GL}{GL}

\DeclareMathOperator{\Isom}{Isom}

\DeclareMathOperator{\AffPushout}{AffPushout}

\newcommand{\id}{\mathrm{id}}
\newcommand{\pr}{\mathrm{pr}}

\newcommand{\Cat}{\mathrm{Cat}}

\DeclareMathOperator{\characteristic}{char}

\newcommand{\an}{\mathrm{an}}

\newcommand{\et}{\mathrm{\acute{e}t}}

\newcommand{\triv}{\mathrm{triv}}

\newcommand{\ev}{\mathrm{ev}}
\newcommand{\coev}{\mathrm{coev}}
\DeclareMathOperator{\PSh}{Psh}

\DeclareMathOperator{\Bund}{Bund}
\DeclareMathOperator{\Gerbe}{Gerbe}
\DeclareMathOperator{\Out}{Out}
\DeclareMathOperator{\act}{act}
\DeclareMathOperator{\qcqs}{qcqs}

\DeclareMathOperator{\CAlg}{CAlg}
\DeclareMathOperator{\art}{art}
\newcommand{\CAlgan}{\CAlg^{\an}}
\newcommand{\CAlganArtLambdaK}{\CAlg^{\an,\art}_{\Lambda//k}}

\DeclareMathOperator{\TwForm}{TwForm}

\DeclareMathOperator{\fib}{fib}

\newcommand{\cO}{\mathcal{O}}
\newcommand{\cA}{\mathcal{A}}
\newcommand{\cC}{\mathcal{C}}

\newcommand{\cH}{\mathcal{H}}

\newcommand{\cJ}{\mathcal{J}}

\newcommand{\cP}{\mathcal{P}}

\newcommand{\cS}{\mathcal{S}}

\newcommand{\cX}{\mathcal{X}}
\newcommand{\cY}{\mathcal{Y}}

\newcommand{\rE}{\mathrm{E}}

\newcommand{\rH}{\mathrm{H}}

\newcommand{\rL}{\mathrm{L}}

\newcommand{\rR}{\mathrm{R}}

\newcommand{\rT}{\mathrm{T}}

\newcommand{\bL}{\mathbb{L}}

\newcommand{\bB}{\mathbf{B}}
\newcommand{\bC}{\mathbf{C}}

\newcommand{\bE}{\mathbf{E}}
\newcommand{\bG}{\mathbf{G}}
\newcommand{\bH}{\mathbf{H}}

\newcommand{\bZ}{\mathbf{Z}}

\newcommand{\sMap}{\mathcal{M}ap}
\newcommand{\sIsom}{\mathcal{I}som}
\newcommand{\sAut}{\mathcal{A}ut}
\newcommand{\ZModTwo}{{\underline{\bZ/2\bZ}}}
\DeclareMathOperator{\Mor}{Mor}
\newcommand{\Dqccn}{\mathrm{D}_{\mathrm{qc},\geq 0}}

\DeclareMathOperator{\nc}{nc}
\DeclareMathOperator{\FM}{FM}
\DeclareMathOperator{\HT}{HT}
\DeclareMathOperator{\HKR}{HKR}
\DeclareMathOperator{\scalloped}{scalloped}

\usepackage{xpatch}
\makeatletter   
\xpatchcmd{\@tocline}
{\hfil\hbox to\@pnumwidth{\@tocpagenum{#7}}\par}
{\ifnum#1<0\hfill\else\dotfill\fi\hbox to\@pnumwidth{\@tocpagenum{#7}}\par}
{}{}
\makeatother

\makeatletter
\renewcommand\part{%
   \if@noskipsec \leavevmode \fi
   \par
   \addvspace{4ex}%
   \@afterindentfalse
   \secdef\@part\@spart}

\def\@part[#1]#2{%
    \ifnum \c@secnumdepth >\m@ne
      \refstepcounter{part}%
      \addcontentsline{toc}{part}{Part \thepart.\hspace{1em}#1}%
    \else
      \addcontentsline{toc}{part}{#1}%
    \fi
    {\parindent \z@ \raggedright
     \interlinepenalty \@M
     \normalfont
     \ifnum \c@secnumdepth >\m@ne
     \centering 
     \Large\bfseries \partname\nobreakspace\thepart     
       \nobreak. 
     \fi
     \Large \bfseries { #2}%
     \par}%
    \nobreak
    \vskip 3ex
    \@afterheading}
\def\@spart#1{%
    {\parindent \z@ \raggedright
     \interlinepenalty \@M
     \normalfont
     \huge \bfseries #1\par}%
     \nobreak
     \vskip 3ex
     \@afterheading}
\makeatother

\renewcommand{\thepart}{\Roman{part}}

\makeatletter 
\def\l@subsection{\@tocline{2}{0pt}{3pc}{6pc}{}} 
\makeatother

\begin{document}

\title{Derived Deformation Theory of Gerbes and Twisted Derived Categories}

\author{Riku Kurama }
\address{Department of Mathematics, University of Michigan, Ann Arbor, MI 48109 \smallskip}
\email{rkurama@umich.edu}

\begin{abstract} 
As we infinitesimally deform principal $\bB\bG_m$-bundles of scalloped derived stacks, the associated twisted derived category also deforms. By formalizing this idea, we define what we call the twisted-to-noncommutative map and use it to answer a special case of Lieblich's question on derived invariance of liftability. 
Along the way, we also develop the derived deformation theory of principal $\bB G$-bundles and $G$-gerbes, as well as the theory of Atiyah classes for derived stacks and its twisted analogue. 
Other results of their own interests include: (1)  derived deformations of scalloped derived stacks form formal moduli problems, and (2) for a smooth proper  variety $X$ with $\rH^0(X,\bigwedge^2\rT_X)=0$, the compactly generated deformations of $\Dqc(X)$ over classical thickenings define a deformation functor.
\end{abstract}

\maketitle

\setcounter{tocdepth}{2}
\newcounter{appendixFix}
\tableofcontents

\section{Introduction}

\subsection{Motivations}\label{subsection-intro-Motivations}
In noncommutative algebraic geometry, we try to understand varieties via their derived categories, and as we try to understand the behaviors of categories in families, it is natural to develop their deformation theory.    
The work of Toda (\cite{Toda}) illuminates  general expectations in the geometric setting. As he explains in the paper, it was expected for some time that the Hochschild cohomology controls the deformations of categories, and he 
gave a construction realizing this philosophy;  
given a smooth projective complex variety $X$ and a class $u\in\HH^2(X/\bC)\simeq \rH^0(X,\bigwedge^2\rT_{X/\bC})\oplus\rH^1(X,\rT_{X/\bC})\oplus \rH^2(X,\cO_X)$, he constructed a deformation $\mathrm{D^b_{coh}}(X,u)$ of $\mathrm{D^b_{coh}}(X)$ over $\bC[\epsilon]/(\epsilon^2)$. Here, the identification comes from the celebrated HKR isomorphism, and he explains that the three components of this decomposition morally correspond to the
(purely) noncommutative deformations, commutative deformations (i.e.\ deformations of schemes) and the gerby deformations of the twists (i.e.\ Brauer classes) respectively. 

Parts of this picture have been formalized using the language of derived algebraic geometry, and
the most relevant work in the literature for our purposes is the paper of Moore (\cite{Moore2}).
In this paper, for a smooth proper $k$-variety $Y$, Moore studied what he called the ``commutative-to-noncommutative period map'' 
\[
\theta:
\Def_{Y}\to \CatDef_{\Dqc(Y)}.
\]
$\Def_{Y}$ is the formal moduli problem (i.e.\ the derived analogue of  a deformation functor) of $Y$, and $\CatDef_{\Dqc(Y)}$ is the functor parameterizing the infinitesimal derived deformations of $\Dqc(Y)$, as defined by Lurie in \cite{SAG}. The above map $\theta$ is simply the association from a (derived) lift of $Y$ to its derived category. Using the work of Lurie (\cite{SAG}), we can approximate $\CatDef_{\Dqc(Y)}$ by a genuine formal moduli problem, and this approximation yields
\[
\theta^{\wedge}:
\Def_{Y}\to \CatDef_{\Dqc(Y)}^{\wedge}.
\]
This morphism induces a map of tangent fibres, which are an analogue of tangent spaces in the derived setting, and it takes the form
\[
T_{\theta^{\wedge}
}:
\rR\Gamma(Y,\rT_{Y/k})[1]\to \HH^\bullet(Y/k)[2]
\] where $\HH^\bullet(Y/k)[2]$ is a shift of the total Hochschild cohomology. 
On the other hand, if we assume for simplicity that $\characteristic(k)>\dim Y$ or $\characteristic(k)=0$, 
we have a morphism of Caldararu (\cite{Caldararu-Mukai-pairing-II})
\[
I:\Delta^*\Delta_*\cO_Y\to \bigoplus_{i}\Omega^i_{Y/k}[i]
\] (where the functors are derived) which induces the HKR isomorphism \[\bigoplus_i\rR\Gamma(Y,\bigwedge^i\rT_{Y/k})[-i]\simeq \HH^\bullet(Y/k).\]
One can describe the map $T_{\theta^{\wedge}}$ via the universal Atiyah class on $Y$, and use it to identify the map with the summand of the above HKR isomorphism with $i=1$. In this way, Moore described a precise relation between commutative deformations of $Y$ and noncommutative deformations of $\Dqc(Y)$ from the derived deformation theoretic perspective.

In the present paper, we develop an analogous theory which accounts for both the commutative deformations and the gerby deformations. 
As a direct analogue to the above story, for a reasonable principal $\bB\bG_m$-bundle $f:X\to Y$ of 1-Artin derived stacks, we construct the 
\textbf{``twisted-to-noncommutative''} maps
\[
\Theta:
\Def_X\to \CatDef_{\Dqc^1(X/Y)},\quad \Theta^\wedge:
\Def_X\to \CatDef_{\Dqc^1(X/Y)}^\wedge,
\]
where $\Dqc^1(X/Y)$ is the (1-)twisted derived category of $f.$
When $Y$ is a smooth proper $k$-variety, any infinitesimal deformation $\cY$ of $Y$ defines an the infinitesimal deformation $\cY\times\bB\bG_m$ of $X:=Y\times\bB\bG_m$, defining a map $\Def_Y\to\Def_{Y\times\bB\bG_m}$. Then, we can recover the map $\theta$ as the composition
\[
\Def_Y\to\Def_{Y\times\bB\bG_m}\xrightarrow[]{\Theta}\CatDef_{\Dqc^1(Y\times\bB\bG_m/Y)}\simeq \CatDef_{\Dqc(Y)}.
\]
This way, $\Theta$ remembers commutative deformations as special cases of twisted deformations. 

As an application, we will answer a special case of the following question of Lieblich: 
\begin{question}[{\cite[Question 2.7]{AntieauBragg_THH}}]\label{question-lieblich}
    Let $X,Y$ be smooth proper varieties over a perfect field of characteristic $p>0$ with equivalent derived categories.
    If $X$ lifts to characteristic 0 (or infinitesimally to rings like $W_2(k)$), does $Y$ also lift?
\end{question}
We can answer this question when one cohomological vanishing holds:
\begin{theorem}[Theorem {\ref{theorem-equivalence-classical-deformation-functors}(3)}]\label{theorem-main-theorem-preview}
Let $\Lambda$ be a complete local noetherian ring with residue field $k$ of positive characteristic.
Let 
$X,Y$ be smooth proper $k$-schemes 
that are derived equivalent, and assume $\rH^0(Y,\bigwedge^2\rT_{Y/k})=0.$ 
Let $A$
be an artin local classical $\Lambda$-algebra with an $\Lambda$-algebra map $A\to k$. If $X$ lifts to $A$ (or more generally if
$X\times\bB\bG_m$ lifts to $A$), then $Y$ also lifts to $A$.
\end{theorem}
\begin{remark}\label{remark-special-known-case-main-theorem}
A special case of this theorem where we assume
$
\rH^1(Y,\cO_Y)=\rH^2(Y,\cO_Y)=\rH^0(Y,\bigwedge^2\rT_{Y/k})=0$ and drop the statement on $X\times\bB\bG_m$ 
can be obtained by
using the commutative-to-noncommutative map instead.  
Accordingly, this special case 
was known previously 
in \cite[Corollary 1.4(1)]{Moore2} as well as in an unpublished joint work of the author with Alexander Perry (see \cite[\S1.2]{Moore2}) which used the same technique as Moore.
 \end{remark}
We also remark that certain variations of (derived or classical) formal deformations of gerbes (or schemes with Brauer classes) were previously considered in \cite{BranTael} and \cite{Bragg-twisted-K3-lift}.

\subsection{Highlights of other key results}
In the following, we always assume Artin derived stacks are qcqs (with affine stabilizers by our convention). Also, fix be a noetherian complete local ring $\Lambda$ with residue field $k$.
\subsubsection{Scalloped derived stacks}
When developing the derived deformation theory of gerbes, we necessarily need to study non-Deligne Mumford derived stacks. In addition to $\bG_m$-gerbes from earlier, it is also desirable to allow for $\mu_p$-gerbes which play a pivotal role in \cite{bragg-lieblich-twistor}, in which case the stabilizers are not smooth either.
However, arbitrary Artin derived stacks can be poorly behaved (see for example \cite[Theorem 1.0.1]{Stefanich} and \cite[Remark 1.3]{Moore2}).
For this reason, we will focus on the class of \textbf{(nicely) scalloped} (qcqs 1-Artin) derived stacks as studied in \cite{Khan-Ravi}. 
Roughly these are the 1-Artin derived stacks with nice stabilizers (in the sense of Definition \ref{definition-nice-embeddable-groups}).
They enjoy many pleasant properties:
\begin{enumerate}
    \item They are perfect (\cite{Khan-Ravi}) and hence have dualizable derived categories.
    \item Their derived categories are tamely compactly generated (Proposition \ref{proposition-tamely-compactly-generated}).
    \item For a nice commutative group scheme $G$, when $f:X\to Y$ is a principal $\bB G$-bundle where $Y$ is scalloped, then $X$ is also scalloped (Proposition \ref{principal-BG-bundle-over-scalloped-is-scalloped}).
    \item A version of local structure theorem holds (see \cite[Theorem 1.13]{Alper-Hall-HL-Rydh} and Proposition \ref{proposition-khan-ravi-scalloped-lift-cover}).
\end{enumerate}
Along the way, we prove the following key results about scalloped derived stacks.
\begin{theorem}[Existence of Ferrand pushouts along nilpotent immersions, Proposition \ref{proposition-existence-geometric-pushout}, Lemma \ref{lemma-geometric-pushout-pushout-among-qcqs-Artin-stacks}]
Consider a diagram of scalloped (qcqs 1-Artin) derived stacks:
$Y_1 \leftarrow Y_0 \to Y_2.$
Assume that the arrows are nilpotent closed immersions. Then, it admits a 
pushout in the $\infty$-category of qcqs 1-Artin derived stacks, and the pushout is scalloped.
\end{theorem}
\begin{theorem}[Proposition \ref{proposition-schlessinger-criterion-scalloped-d-stack}]
    Consider a pushout diagram of derived affine schemes:
    \[
    \xymatrix{
X_{0}\ar[r]\ar[d]&X_1\ar[d]\\
X_2\ar[r]& X_3
    }
    \]
    where $X_{0}\to X_i$  is a nilpotent closed immersion for $i=1,2$ (i.e.\ it is a nilpotent closed immersion on the classical locus).
    Then, pullback induces an equivalence
    \[
   \Phi: \dArtSt^{\scalloped}_{/X_3}\simeq  
\dArtSt^{\scalloped}_{/X_1}
\times_{\dArtSt^{\scalloped}_{/X_{0}}} 
\dArtSt^{\scalloped}_{/X_2}
    \] where $\dArtSt^{\scalloped}_{/T}$ is the $\infty$-categories of scalloped derived stacks over $T$. 
\end{theorem}
\begin{theorem}[Proposition \ref{proposition-Def-X-FMP}]\label{theorem-Def-X-FMP-preview} 
Let $X$ be a scalloped (qcqs 1-Artin) derived stack over $k$. Then, the functor $\Def_X:\CAlganArtLambdaK\to\widehat{\cS}$ sending an artinian derived ring $A$ to the space of deformations of $X$ over $A$ is a formal moduli problem. 
\end{theorem}
\begin{remark}
Analogous theorems are known for spectral DM stacks by the work of Lurie:
\begin{enumerate}
    \item Existence of pushouts is known in this setting for a pair of closed immersions (\cite[Theorem 6.1]{DAG-IX}).
    \item The clutching of spectral DM stacks is proved for bases living in a pushout diagram of spectral DM stacks with closed immersions 
(\cite[Theorem 9.1]{DAG-IX}).
\item Derived deformations of spectral DM stacks defines a formal moduli problem (\cite[Theorem 16.3.3.1]{SAG}).
\end{enumerate} 
Additionally, for classical algebraic stacks, the existence of Ferrand pushouts are known for one affine morphism and one closed immersion (\cite[Theorem 1.8]{Alper-Hall-HL-Rydh}).
\end{remark}
\subsubsection{Derived deformations of gerbes}
By Theorem \ref{theorem-Def-X-FMP-preview}, we can  talk about formal moduli problems of scalloped derived stacks, and analogously for  a morphism $f:X\to Y$ of scalloped derived stacks over $k$,  the formal moduli problem $\Def_f$ of the morphism $f$ makes sense. 
We are interested in deforming  principal $\bB G$-bundles (Definition \ref{definition-principal-G-bundle}) and $G$-gerbes (Definition \ref{definition-G-gerbe-without-banding}). The following theorem makes the picture clear when $G$ is nice and commutative. 
\begin{theorem}[Proposition \ref{proposition-equiv-FMP-f-X-Gmgerbe-BGmbundle}]
Let $G$ be a nice commutative group scheme over $\Lambda$. 
Suppose that $f:X\to Y$ is a principal $\bB G$-bundle of scalloped derived stacks over $k$. 
     The morphisms of $\cS$-valued functors
    \[
    \Def^{\Bund(\bB G)}_f\to \Def_f^{\Gerbe(G)}\to \Def_f\to \Def_X
    \]are all equivalences.
    Here, $\Def^{\Bund(\bB G)}_f$ is the formal moduli problem parametrizing deformations of $f$ as principal $\bB G$-bundles, and $\Def_f^{\Gerbe(G)}$ is a functor parametrizing its deformations as $\bG$-gerbes (without bandings).
    In particular, $\Def_f^{\Gerbe(G)}$ is a formal moduli problem.
\end{theorem}
Due to this theorem, we can define the twisted-to-noncommutative map $\Theta:\Def_X\simeq \Def_{f}^{\Bund(\bB G)}\to \CatDef_{\Dqc^1(X/Y)}$ which sends a derived deformation of $f$ to its twisted derived category.
\subsubsection{Atiyah class for derived stacks and nonconnective geometry}
In order to study the twisted analogue of the commutative-to-noncommutative map, we will construct the
\textbf{twisted Atiyah class}. One way to do so it is to restrict the universal Atiyah class on a principal $\bB\bG_m$-bundle (or equivalently a $\bG_m$-gerbe with trivialized banding) to the twisted derived category in an appropriate sense. This necessitates the formalism of (untwisted) Atiyah classes on Artin derived stacks. 
\begin{theorem}[Definition \ref{definition-atiyah-transformation}]
Let $k$ be a field.
    Let $X$ be a derived stack with a cotangent complex over $k$. Let $E\in\Dqc(X)$. Then, for each $\alpha\in\Map_X(\bL_{X/k},E)$, we can associate the \textbf{universal $\alpha$-Atiyah transformation of $X$} $\At_\alpha(-):\id\to(-)\otimes E[1]$ which is a morphism in the $\infty$-category $\Fun_k^\rL(\Dqc(X),\Dqc(X))$. 
\end{theorem}
\begin{remark}
    In \cite{Kuhn}, Kuhn constructed Atiyah classes for classical algebraic stacks. For geometric derived stacks in characteristic 0, there exists a construction of 
    Atiyah classes via the moduli of perfect complexes (see \cite{STV} and  \cite{Hennion}).
\end{remark}
While the theorem above might look different from the conventional Atiyah classes, if we choose $E=\bL_{X/k},\alpha=\id_{\bL_{X/k}}$ and some object $F\in\Dqc(X)$, the universal $\alpha$-Atiyah transformation specializes to $\At_{\id_{\bL_{X/k}}}(F):F\to F\otimes\bL_{X/k}[1]$ which has the same form as the conventional Atiyah class of $F$. Also, the above construction satisfies the usual properties like compatibility with morphisms of derived stacks and tensor products of modules under suitable hypotheses.

Our construction of the Atiyah transformation is similar to Lurie (\cite[\S19.2.2]{SAG}) and Moore (\cite[\S3]{Moore2}), where they consider the cases of spectral DM stacks and derived schemes respectively, but in the non-Deligne Mumford setting, certain difficulty arises. Namely, our construction uses the (trivial) square-zero extension $X[E]$ of $X$ by $E$, which is naturally a nonconnective derived stack. While one can require $E$ to be connective to avoid this issue, that would rule out $E=\bL_{X/k}$ from considerations. For this reason, the theorem above used nonconnective derived stacks as developed in \cite{blowup-normal-bundle-nonconnective-DAG} (where the affines come from the derived commutative rings in the sense of Bhatt-Mathew \cite{Raksit}).
\subsubsection{Twisted Atiyah class and the twisted-to-noncommutative map}
Similarly to the untwisted Atiyah class, we also construct the twisted version:
\begin{theorem}[Definition \ref{definition-twisted-atiyah-transformation}]
Let $k$ be a field.
    Let $f:X\to Y$ be a principal $\bB\bG_m$-bundle over a scalloped derived stack $Y$ over $k$, and let $E\in\Dqc^0(X/Y)\simeq\Dqc(Y)$ be $0$-twisted. Then, for each $\alpha\in\Map_X(\bL_{X/k},E)$, we can associate the \textbf{twisted universal $\alpha$-Atiyah transformation of $X/Y$},
    $
    \widetilde{\At}_{\alpha}(-):\id\to(-)\otimes E[1]
    $ which is a morphism in  $\Fun_k^\rL(\Dqc^1(X/Y),\Dqc^1(X/Y))$.
\end{theorem}
By fixing $E$ and varying $\alpha$, this gives us the\textbf{ twisted Atiyah map of $E$}:
\[
\widetilde{\bAt}_E:\Map_{X//k}(X[E],X)\simeq\Map_X(\bL_{X/k},E)\to \Map_{\Fun_k^{\rL}(\Dqc^1(X),\Dqc^1(X))}(\id,(-)\otimes E[1]).
\]
We can upgrade this map to a version with mapping spectra, and setting $E=\cO_X[1]$ yields:
\[
\widetilde{\At}_{\cO_X[1]}^{\FM}:\RHom(\bL_{X/k},\cO_X)[1]\to\HH^\bullet(\Dqc^1(X/Y)/k)[2].
\]
This map connects to the twisted-to-noncommutative map:
\begin{theorem}[Proposition \ref{proposition-map-on-tangent-fibers-twisted-to-noncommutative}]\label{theorem-tangent-fiber-Theta-twisted-Atiyah}
    In the above situation, the map of tangent fibers induced by the twisted-to-noncommutative map $T_{\Theta^\wedge}:T_{\Def_X}\to T_{\CatDef^\wedge_{\Dqc^1(X/Y)}}$ is canonically identified with 
    \[
\widetilde{\At}_{\cO_X[1]}^{\FM}:
\RHom(\bL_{X/k},\cO_X)[1]\to\HH^\bullet(\Dqc^1(X/Y)/k)[2]
\]
\end{theorem}
\subsubsection{The twisted-to-noncommutative map in the split case}
Assume $Y$ is a smooth proper variety. 
When the principal $\bB G$-bundle $f:X\to Y$ is trivial, the map in Theorem \ref{theorem-tangent-fiber-Theta-twisted-Atiyah} becomes simpler.
Similarly to the discussion in \S\ref{subsection-intro-Motivations}, we have a map 
$I_{\leq 1}:\Delta^*\Delta_*\cO_Y\to \bigoplus_{i=0,1}\Omega^i_{Y/k}[i]
$ which appears in the HKR theory in any characteristic.
We  use it to define the map
\[ I^{\HKR}_{\leq1}[2]:
\RHom_Y(\bL_{Y/k},\cO_Y[1])
\oplus \RHom_Y(\cO_Y,\cO_Y[2])
\to 
\RHom_{Y\times_k Y}(\Delta_*\cO_Y,\Delta_*\cO_Y[2]).\]
\begin{theorem}[Theorem \ref{theorem-tangent-fibre-trivial-gerbe-hkr}]
Let $Y$ be a smooth proper variety over $k$, and consider the trivial $\bB\bG_m$-bundle
$f:X=Y\times\bB\bG_m\to Y$. The map  $T_{\Theta^\wedge}$ on the tangent fibres induced by the twisted-to-noncommutative map can be identified with $ I^{\HKR}_{\leq1}[2]$.
If $\rH^0(Y,\bigwedge^2\rT_{Y/k})=0$, 
$\rH^0(T_{\Theta^\wedge}),\rH^{-1}(T_{\Theta^\wedge})$ are isomorphisms, and $\rH^1(T_{\Theta^\wedge})$ is injective. 
\end{theorem}
\subsubsection{Applications to liftablity and classical deformations of categories}
By applying the results developed so far, we can deduce several interesting theorems, including Theorem \ref{theorem-main-theorem-preview} at the beginning of the paper:
\begin{theorem}[Theorem \ref{theorem-sqz-lift-twisted-Dqc}]\label{theorem-preview-sqz-lift-twisted-Dqc}
    Let $Y$ be a smooth proper $k$-scheme. Assume $\characteristic(k)\neq2$ or $\rH^0(Y,\bigwedge^2\rT_Y)=0.$ Let $A'\to A$ be a {small} extension of classical Artinian local $\Lambda$-algebras with residue field $k.$
    Suppose $X:=Y\times\bB\bG_m$ lifts to $\cX$ over $A$. 
    Then, the principal $\bB\bG_m$-bundle $X\to Y$ lifts uniquely to $\cX\to \cY$ over $A$ up to equivalence. Moreover, 
    $\cX$ lifts to $A'$ if and only if $\Dqc^1(\cX/\cY)$ lifts to $A'$. 
\end{theorem}
\begin{theorem}[Corollary \ref{corollary-unique-sqz-lift-FM-partner}]\label{theorem-preview-unique-sqz-lift-FM-partner}
 Let $X,Y$ be smooth proper varieties over $k$, and assume we have a Fourier-Mukai equivalence $\Phi:\Dqc(X)\simeq\Dqc(Y)$. Let $A\to k$ be a square-zero thickening as a classical $\Lambda$-algebra.
    Let $\cX$ be a deformation of $X\times\bB\bG_m$ over $A$ (which comes with a principal $\bB\bG_m$-bundle structure $\cX\to \cX_{\coarse}$ lifting the projection $X\times\bB\bG_m\to X$). Then, there is at most one deformation $\cY$ of $Y\times\bB\bG_m$ over $A$ up to equivalence (which comes with a principal $\bB\bG_m$-bundle structure $\cY\to \cY_{\coarse}$ lifting the projection $Y\times\bB\bG_m\to Y$) which admits an $A$-linear Fourier-Mukai equivalence $\Dqc^1(\cX/\cX_{\coarse})\simeq \Dqc^1(\cY/\cY_{\coarse})$ reducing to 
    \[\Phi:\Dqc^1(X\times\bB\bG_m)\simeq\Dqc(X)\simeq \Dqc(Y)\simeq\Dqc^1(Y\times\bB\bG_m).\]
    In particular, for any (untwisted) lift $\cX$ of $X$ over $A$, there is at most one (untwisted) lift $\cY$ of $Y$ up to equivalence which admits an $A$-linear Fourier-Mukai equivalence $\Dqc(\cX)\simeq\Dqc(\cY)$ reducing to $\Phi:\Dqc(X)\simeq\Dqc(Y)$.

\end{theorem}
Below, when $Y$ is smooth proper over $k$, we also consider the compactly generated variant $\CatDef_{\Dqc(Y)}^{\omega}$ parametrizing compactly generated derived deformations of $\Dqc(Y)$.
We also let $\CatDef_{\Dqc(Y)}^{\omega,\cl}$ denote the composition
    \[
\CRing_{\Lambda//k}^{\art}\hookrightarrow\CAlganArtLambdaK\xrightarrow{\CatDef_{\Dqc(Y)}^{\omega}}{\cS}\xrightarrow{\pi_0}{\mathrm{Set}}
    \] where $\CRing_{\Lambda//k}^{\art}\subset \CAlganArtLambdaK$ is the full subcategory of classical rings.
\begin{theorem}[Theorem {\ref{theorem-equivalence-classical-deformation-functors}}, Theorem \ref{theorem-compactly-generated-deformation-functor}]\label{theorem-preview-equivalence-classical-deformation-functors-compact-generation}
    \leavevmode
    \begin{enumerate}
        \item  Let $X$ be a smooth proper $k$-scheme with $\rH^0(X,\bigwedge^2\rT_{X/k})=0.$ Then, (the compactly generated variant of) the twisted-to-noncommutative map induces an isomorphism of the associated discrete functors on classical rings:
        \[
        \Def_{X\times\bB\bG_m}^{\cl}\simeq 
        \CatDef_{\Dqc(X)}^{\omega,\cl}\simeq \CatDef_{\Dqc(X)}^{\wedge,\cl}.
        \] 
        In particular, $\CatDef_{\Dqc(X)}^{\omega,\cl}$ is a deformation functor.
        \item When $X,Y$ are two such $k$-schemes, any (untwisted) derived equivalence between $X$ and $Y$ induces an isomorphism of classical deformation functors:
        \[
        \Def_{X\times\bB\bG_m}^{\cl}\simeq \Def_{Y\times\bB\bG_m}^{\cl}
        \] 
    \end{enumerate}
\end{theorem}
\begin{remark}
   Untwisted analogues of some of the above theorems were previously obtained by \cite{Moore2} (and the independent unpublished work by the author and Alexander Perry).
   \begin{enumerate}
       \item An untwisted analogue of Theorem \ref{theorem-preview-sqz-lift-twisted-Dqc} (for qcqs derived schemes) is proved in \cite[Theorem B(1)]{Moore2}.
       \item An untwisted analogue of Theorem \ref{theorem-preview-equivalence-classical-deformation-functors-compact-generation} with additional vanishing hypotheses was obtained in \cite[Corollary 1.4]{Moore2}, with the exception of  the claim about  $\CatDef^{\omega,\cl}_{\Dqc(X)}$.
   \end{enumerate}
\end{remark}
\begin{remark}
We emphasize that the assertion that $\CatDef_{\Dqc(X)}^{\omega,\cl}$ is a deformation functor appears quite nontrivial in general, since $\CatDef_{\Dqc(X)}^\omega$ is not known to be a formal moduli problem. See also Remark \ref{remark-deformation-functor-CatDef-question}.
\end{remark}

\subsection{Acknowledgment}
{During the preparation of this paper, the author was partially supported
by the Murata Overseas Scholarship. 
The author thanks his advisor 
Alexander Perry for continued support and helpful discussions.
The author thanks 
Max Lieblich,
Samuel Moore,
Saket Shah and
Gleb Terentiuk for helpful discussions.
The author also thanks 
Jakub Witaszek for first bringing up a version of Question \ref{question-lieblich} to the author and for the interests in this project.
\subsection{Tools disclosure}
 This paper was entirely written by the author, except that the author initially attempted to transcribe several hand-written diagrams 
to LaTeX codes in \S\ref{section-Atiyah-class-for-derived-stacks} 
using AI tools.
In addition, the author used 
ChatGPT up to 5.6 Sol/Gemini up to 3.1 for 
searching for results in the literature, 
checking correctness and revising. {In particular, the original proofs of Corollary \ref{corollary-BG-bundle-formal-moduli-problem} and Lemma \ref{lemma-affineness-derived-vs-classical} had issues and were subsequently fixed using suggestions from ChatGPT.}
The author is responsible for any mistakes.}
\subsection{Conventions}
\begin{enumerate}
    \item 
The paper uses the combination of homological (denoted by subscripts) and cohomological (denoted by superscripts) indexings. For instance, $\Dqccn(X)$ denotes the $\infty$-category of connective quasi-coherent modules.
\item
$\cS$ denotes the $\infty$-category of spaces.
\item 
For a morphism $f:X\to Y$ in an $\infty$-category with finite products, its Cech nerve is denoted by $C_{\bullet}(f)=C_{\bullet}(X/Y)$, and $C_{n}(f)=C_n(X/Y)$ is its $n$-th piece. 
\item
Standard functors like pushforward and pullback along a morphism $f:X\to Y$ of derived stacks are often derived by default, so that we write them as $f_*$ and $f^*$ (instead of $\rR f_*$ and $\rL f^*$).
\item
For a morphism $f:X\to Y$ of derived stacks relative to some base $S$, we write $[f]:f^*\bL_{Y/S}\to \bL_{X/S},[f]^A:\bL_{Y/S}\to f_*\bL_{X/S}$ for the induced morphisms of cotangent complexes (when they exist).
\item
We will often drop ``homotopy'' from our terminology for shortness (e.g.\ homotopy (co)limits and commutativity of diagrams up to homotopy).
\item
As mentioned in \S\ref{section-basic-derived-stacks-principal-bundles}, most Artin derived stacks considered in the paper are 1-Artin qcqs (with affine stabilizer by convention). 
\item For a nice commutative group scheme over a classical affine, the delooping $\bB G$ (in the $\infty$-topos of \'etale derived stacks) is computed in the fppf sense. 
\item
For $G=\bG_m,\mu_n$ (where $n$ is a positive integer), in our convention, the term ``$G$-gerbe'' refers to a $G$-gerbe without banding (\S\ref{subsection-banding}). We also use the terminology ``$G$-gerbes with trivialized bandings'' to distinguish banded gerbes (we show in \S\ref{section-basic-derived-stacks-principal-bundles} that they are the same as principal $\bB G$-bundles). 
\end{enumerate}
\part{Gerbes and twisted derived categories}
\section{Basic results on derived stacks and principal bundles}\label{section-basic-derived-stacks-principal-bundles}
We recall and prove some basics of derived stacks.
By default, following \cite{DAG-gaitsgory-1}, derived stacks in this paper are $\infty$-sheaves in the \'etale topology, and Artin (geometric) derived stacks have covers in the smooth topology as in \cite[Ch.2 2.1.6]{DAG-gaitsgory-1}, unless stated otherwise. 
In some places, fppf topology will be used, but as far as geometric derived stacks are concerned, this will not affect the theory (\cite[Remark 4.1.4]{DAG-gaitsgory-1}), and this is clarified in \S\ref{subsection-principal-BG-bundle-etale-fppf}.
Most algebraic derived stacks in the paper are assumed to be 1-Artin derived stacks (\cite[\S1.11]{Khan-Ravi}), except for rare instances like $\bB^2\bG_m$.
Let us quickly recall this convention. \textbf{$0$-Artin derived stacks (or $0$-geometric derived stacks)} are disjoint unions of derived affines. \textbf{$1$-Artin derived stacks (or $1$-geometric derived stacks)} have $0$-geometric diagonal morphisms (i.e. relatively representable by $0$-Artin derived stacks) and has a $0$-geometric smooth surjective cover by a $0$-Artin derived stack (note that smooth surjective morphisms make sense as $0$-geometric morphisms). We also call these smooth covers smooth atlases. 
In our convention, while separated schemes are $1$-Artin, not all schemes in the traditional sense are $1$-Artin.
The $\infty$-category of 1-Artin derived stacks is denoted by $\dArtSt,$ or $\dArtSt_{/S}$ if we work over the base derived affine $S$. 
\begin{remark}\label{remark-Toen-fppf-derived-geometric-stack}
    As mentioned above, by a result of Toen (\cite{toen-descent}, or \cite[Remark 4.1.4]{DAG-gaitsgory-1}),  (possibly higher) Artin stacks can equivalently formulated as (higher) geometric fppf derived stacks with fppf covers. More precisely, we let $\dArtSt^n$ denote the $\infty$-category of $n$-Artin derived stacks in the \'etale topology (with smooth covers), and let $\dGeomSt^n_{\fppf}$ be the $\infty$-category of $n$-geometric derived stacks in the fppf topology (with fppf covers). We have an equivalence
    \[
    \dArtSt^n\simeq \dGeomSt^n_{\fppf}
    \] where the functor from the left to the right is induced by fppf $\infty$-sheafification.
    Note also that any Artin derived stack is an fppf $\infty$-sheaf and that the inverse is induced by the inclusion functor from the $\infty$-category of fppf derived stacks to the $\infty$-category of \'etale derived stacks.

\end{remark}

A 1-Artin derived stack is \textbf{quasi-compact} if it admits a smooth atlas from a derived affine (\cite[4.1.9]{DAG-gaitsgory-1}, \cite[Definition 1.3.6.4]{HAG2}). A morphism $f:X\to Y$ of such derived stacks is \textbf{quasi-compact} if base changing $Y$ from derived affines $U\to Y$ yields quasi-compact base changes $X\times_YU$.
A 1-Artin derived stack is \textbf{qcqs} if it's quasi-compact and additionally the diagonal morphism is quasi-compact (in the $0$-representable sense). Note that the definition in \cite[4.1.9]{DAG-gaitsgory-1} looks more complicated, but in our case the additional condition that the diagonal itself be quasi-separated is automatic since the diagonal is $0$-geometric and quasi-compact (hence affine). The $\infty$-category of qcqs 1-Artin derived stacks is denoted by $\dArtSt^{\qcqs}$.
\begin{remark}\label{remark-affine-diagonal}
    As noted above, in our convention,  qcqs 1-Artin derived stacks have affine diagonals. In some other conventions, this is not automatic. 
\end{remark}

\begin{lemma}\label{lemma-qc-qs-check-on-classical-locus}
\leavevmode
Let $X$ be a 1-Artin derived stack.
\begin{enumerate}
    \item  $X$ is quasi-compact if and only if $X_{\cl}$ is.
    \item If $Y$ is another 1-Artin derived stack, a morphism $f:X\to Y$ is quasi-compact if and only if $f_{\cl}:X_{\cl}\to Y_{\cl}$ is.
    \item  $X$ is qcqs if and only if $X_{\cl}$ is.
\end{enumerate}
\end{lemma}
\begin{proof}
\leavevmode
\begin{enumerate}
    \item  If we have a cover $U\to X$ from a derived affine, its truncation is a cover of $X_{\cl}$ by a classical affine scheme, so ``only if'' implication holds. 

    Conversely, choose an atlas $U\to X$ by a 0-geometric derived stack. Write $U=\coprod_{i\in I} U_i$ where each $U_i$ is a derived affine. By quasi-compactness of $X_{\cl},$ the truncation $\coprod_{i\in I} U_{\cl}\to X_{\cl}$ can be refined to a finite cover for some finite subset $S\subset I$, and the corresponding original cover $\coprod_{i\in S}U_i\to X$ is an atlas.
    \item This is immediate from (1) once we note that the classical locus functor $(-)_{\cl}$ commutes with fiber products. 
    \item This follows from (1) and (2).
\end{enumerate}
\end{proof}
\begin{lemma}\label{lemma-infinitesimal-invariance-etale}
    Let $A\to A'$ be a morphism of derived rings inducing a nilpotent thickening $\pi_0(A)\to \pi_0(A')$. Then, the induced morphism of the $\infty$-categories of \'etale algebras $\CAlg_{A/}^{\an,\et}\to \CAlg_{A'/}^{\an,\et}$ is an equivalence.
\end{lemma}
\begin{proof}
    The result holds if the morphism $A\to A'$ is the map from $A$ to its discrete truncation
    (see \cite[Theorem 7.5.0.6]{HA}).
    The result is also well-known if both derived rings are classical.  
The general case 
follows by the combination of these cases if we use the truncation maps $A\to \pi_0(A)$ and $B\to\pi_0(B)$.
\end{proof}
\begin{definition}[{\cite[Definition 2.1]{Khan-Ravi}}]\label{definition-nice-embeddable-groups}
    Let $G$ be an fppf (classical) affine group scheme over a (classical) affine $S$.
    \begin{enumerate}
        \item $G$ is called \textbf{nice} if it's an extension of a finite \'etale group scheme, whose order is invertible on $S$, by a group of multiplicative type.
        \item $G$ is \textbf{embeddable } if it's a closed subgroup scheme of some $\GL_S(E)$ where $E$ is a locally free sheaf on $S$.
    \end{enumerate}
\end{definition}
\begin{remark}\label{remark-nice-group-closed-under-extension-subgroup-quotient}
 Over a field, subgroups, quotients and extensions of nice groups are nice. See the discussions after \cite[Definition 1.1]{Hall-Rydh}.
\end{remark}
The following lemma clarifies the behavior of ``niceness'' of group schemes under field extensions. 
\begin{lemma}\label{lemma-nice-group-field-extension}
    Let $G$ be a 
    group scheme over a field $k$. If $G_L$ is nice over $L$ for some algebraically closed field extension $L/k$, then $G$ is nice. In particular, niceness of group schemes over fields can be checked after any field extensions. 
\end{lemma}
\begin{proof}
Let $p=\characteristic(k).$
We first observe that $G\to \Spec k$ is fppf and affine, as it holds after the field extension $L/k$. 
By \cite[Lemma 9.5.1]{FGA-explained}, 
    the identity component $G^\circ\subset G$ is clopen and geometrically connected so that its formation commutes with field extensions. In particular, we can write $G_L^\circ$ which is simultaneously $(G^\circ)_L$ and $(G_L)^\circ$. 
    By assumption, we have a closed subgroup of multiplicative type $H\subset G_L$ such that $G_L/H$ is finite discrete, of order not divisible by $p$ if $p>0$. Explicitly, $H$ is a product of copies of $\bG_m$ and the $\mu_{n_i}$'s for some integers $n_i$. 
    $H$ is automatically open in $G_L,$ being the kernel of a map to a finite discrete group. 
    We have $G^\circ_L\subset H$, as the composition $G^\circ_L\hookrightarrow G_L\to G_L/H$ is trivial by discreteness of $G_L/H$. 
    Hence, $G^\circ_L$ is a group of multiplicative type. 
    $G^\circ$ is then a group of multiplicative type.
    In the short exact sequence
    \[
    1\to H/G^\circ_L\to G_L/G_L^\circ \to G_L/H\to 1,
    \]
    the left and right terms are both finite \'etale, of order not divisible by $p$ if $p>0$, since $H/G_L^\circ$ is a product of copies of the $\mu_{n_i}$'s, where $p$ doesn't divide $n_i$ if $p>0$.
    Hence, $G_L/G_L^\circ $ has the same property, and so does $G/G^\circ.$ We already showed that $G^\circ$ is of multiplicative type, so we are done.  
\end{proof}

\begin{lemma}\label{lemma-nice-groups-fiber-products}
    Let $G,H,K$ be nice group schemes over a field $k$. For some group homomorphisms $f:G\to K,g:H\to K$, form the fiber product $G\times_KH$. The resulting fiber product is nice.
\end{lemma}
\begin{proof}
    We first observe that the fiber product can be written as 
    \[
    G\times_KH\simeq 
    (G\times_kH)\times_{f\times g,K\times_kK,\Delta_K}K.
    \] Then, it's a closed subgroup of the nice group $G\times_kH$, since $K$ is separated.
    By Remark \ref{remark-nice-group-closed-under-extension-subgroup-quotient}, $G\times_KH$ is nice.
\end{proof}
\subsection{Principal $G$-bundles and $G$-gerbes}
Following \cite{NSS}, we will recall principal bundles and gerbes for an $\infty$-topos $\bH$. {The key case is $\bH=\dSt_{/S}$ or its slice, where $S$ is a derived affine, and $\dSt_{/S}$ is the $\infty$-topos of derived stacks for either the \'etale or the fppf topology. }
\begin{definition}
    A morphism $f:X\to Y$ in $\bH$ is an \textbf{effective epimorphism} if the natural morphism $C_\bullet(X/Y)\to Y$ from the Cech nerve of $f$ is the colimit diagram.
\end{definition}
\begin{definition}
    An \textbf{$\infty$-group in $\bH$ }is an $A_\infty$-algebra $G$ in $\bH$ such that $\pi_0(G)$ is a group object in $\tau_{\leq0}\bH$. We write $\Grp(\bH)$ to denote the $\infty$-category of $\infty$-group in $\bH$. 
\end{definition}
\begin{remark}
    Equivalently, $\infty$-groups can be defined as $\infty$-groupoid objects (Segal groupoid objects) whose $0$-th pieces are $*$ (\cite[Proposition 2.22]{NSS}).
\end{remark}
\begin{example}
In the case $\bH=\dSt_{/S}$, any group scheme $G$ over $S$ is an $\infty$-group.
\end{example}
\begin{example}\label{example-BG-group-when-G-commutative}
 Assume $G$ is commutative. In this case, the morphisms in the following Segal groupoid object are group homomorphisms.
 \[
    \xymatrix{
    ... \ar@<0.5ex>[r]\ar@<-0.5ex>[r]\ar@<1.5ex>[r]\ar@<-1.5ex>[r]& G\times G \ar@<1ex>[r]\ar[r]\ar@<-1ex>[r] &  G \ar@<0.5ex>[r]\ar@<-0.5ex>[r]  & S
    }
    \]
   For example, the three morphisms $G\times G\to G$ are given by the group homomorphisms $p_1,p_2,m:G\times G\to G$.
    Applying the delooping functor $\bB:\Grp(\dSt_{/S})\to(\dSt_{/S})^{*/}_{\geq1}$   (\cite[Theorem 2.19]{NSS}), we obtain a Segal groupoid object which exhibits $\bB G$ as an $\infty$-group:
     \[
    \xymatrix{
    ... \ar@<0.5ex>[r]\ar@<-0.5ex>[r]\ar@<1.5ex>[r]\ar@<-1.5ex>[r]&\bB G\times\bB G \ar@<1ex>[r]\ar[r]\ar@<-1ex>[r] & \bB G \ar@<0.5ex>[r]\ar@<-0.5ex>[r]  & S
    }
    \] 
\end{example}

\begin{remark}\label{remark-truncation-group}
Given an $\infty$-group $G$
  and an integer $n\geq0$, the truncation $\tau_{\leq n}$ (in the sense of $\infty$-topos) is still an $\infty$-group (\cite[Definition 4.54]{NSS}), and the morphism to the truncation (i.e.\ unit of adjunction $\id\to \iota_n\tau_{\leq n}$ for the adjoint pair $(\tau_{\leq n},\iota_n)$ where $\iota_n$ is the inclusion of $n$-truncated objects into $\bH$) $G\to \tau_{\leq n}G$ is a morphism of $\infty$-groups. 
\end{remark}
\begin{definition}\label{definition-principal-G-bundle}
Let $X$ be an object of $\bH$ and $G$ be an $\infty$-group in $\bH$. A \textbf{principal $G$-bundle over $X$} is a morphism $q:P\to X$ in $\bH$ and a $G$-action on $P$ (\cite[Definition 3.1]{NSS}) such that  $q$ is the quotient morphism exhibiting $X$ as the colimit $P\sslash G$ of the simplicial diagram describing the $G$-action.

A \textbf{morphism of {principal $G$-bundles over $X$}} is a morphism of $G$-actions fixing the quotient $X.$

The infinity category of principal $G$-bundles over $X$ is denoted by $\Bund_X(G)$, and it is in fact an $\infty$-groupoid (\cite[Remark 3.6]{NSS}).
\end{definition}

\begin{definition}\label{definition-G-gerbe-without-banding}
Let $G$ be an $\infty$-group in $\bH$, and let $S$ be the final object of $\bH$. 
    A \textbf{$G$-gerbe (without banding) (over $S$)} is an object $E$ such that for some effective epimorphism $U\to S,$ there exists an equivalence $E|_U\simeq \bB G|_U.$
    We define $\Gerbe_S(G)$ to be the $\infty$-groupoid of $G$-gerbes over $S$ (i.e.\ we restrict to the maximal $\infty$-groupoid).
\end{definition}
\begin{remark}
    In practice, the $\infty$-topos $\bH$ in the above definition is a slice of an ambient topos like $(\dSt_{/S})_{/X}$.
\end{remark}
\begin{proposition}[{\cite[Proposition 3.12]{NSS}}]
    Let $G$ be an $\infty$-group in an $\infty$-topos $\bH$. A principal $G$-bundle is locally trivial. More precisely, up to changing the base of the principal $G$-bundle along an effective epimorphism, this becomes equivalent to a trivial principal $G$-bundle. 
\end{proposition}
\begin{remark}
    In particular, when $G$ is commutative, principal $\bB G$-bundles (which make sense by Example \ref{example-BG-group-when-G-commutative}) are $G$-gerbes. 

\end{remark}
\begin{remark}
Our definition of $\Gerbe_S(G)$ appears to differ from the convention of \cite{NSS} (as we took the maximal $\infty$-groupoid) while the precise definition is not stated in their paper.
\end{remark}
\begin{proposition}[see for example {\cite[Corollary 4.51]{NSS}}]\label{proposition-gerbe-BAutBG}
We have an equivalence
{
 \[\Gerbe_S(G)\simeq\Map(S,\bB\sAut_{\bH}(\bB G)).\]
 }
\end{proposition}

\subsection{Principal $\bB G$-bundles in the geometric setting}
We will next specialize to the case of derived stacks. 

\begin{remark}\label{remark-segal-groupoid-assoc-group}
Let $S$ be a classical affine, and let $G$ be a nice group scheme over $S$. 
Then, the delooping $\bB G$
is the colimit of the following diagram
\[
    \xymatrix{
    ... \ar@<0.5ex>[r]\ar@<-0.5ex>[r]\ar@<1.5ex>[r]\ar@<-1.5ex>[r]& G\times G \ar@<1ex>[r]\ar[r]\ar@<-1ex>[r] &  G \ar@<0.5ex>[r]\ar@<-0.5ex>[r]  & S
    }
    \] 
    In this paper, { even in the context of \'etale derived stacks,  we construct such deloopings $\bB G$ } in the fppf topology (i.e.\ we compute the colimit in the $\infty$-topos of fppf derived stacks). 
    {For this reason, in this particular instance, the meaning of $\bB G$ diverges from the standard one internal to the $\infty$-topos. Similarly, when we work with principal $\bB G$-bundles and $G$-gerbes in the \'etale derived stack setting, strictly speaking, if we want to apply their abstract definitions from above, they should be considered as principal $\bB G$-bundles and $G$-gerbes in the fppf sense.}
    $G$ is fppf over $S$, so by \cite[Proposition 1.3.4.2]{HAG2} this is 1-geometric as an fppf derived stack with an fppf atlas. We can then view this as a 1-Artin  derived stack in our sense by Remark \ref{remark-Toen-fppf-derived-geometric-stack}. As noted in this remark, the associated \'etale derived stack is obtained by naturally viewing an fppf $\infty$-sheaf as an \'etale $\infty$-sheaf.

    Because of this fppf convention, when $X$ is an  $n$-Artin derived stack (with $n\geq1$), we have
    \[\Map_{\dSt}(X,\bB G)
    \simeq 
    \Map_{\dArtSt^n}(X,\bB G)
    \simeq 
    \Map_{\dGeomSt^n_{\fppf}}(X,\bB G)\simeq \Bund_{X,\fppf}(G),
    \] where $\Bund_{X,\fppf}(G)$ is the $\infty$-groupoid of principal $G$-bundles in the fppf topology. When $G$ is additionally commutative, a similar result holds for $\bB^2 G$ (see also \S\ref{subsection-principal-BG-bundle-etale-fppf}). 
\end{remark}
\begin{remark}
    {When $G$ is a nice commutative group scheme over an affine $S$, principal $\bB G$-bundles over Artin derived stacks can be equivalently considered in the fppf topos or the \'etale topos  if we use the convention as in Remark \ref{remark-segal-groupoid-assoc-group} for the \'etale case (see \S\ref{subsection-principal-BG-bundle-etale-fppf}). }
\end{remark}
\begin{lemma}\label{lemma-BG-qcqs-1-artin-smooth-surj}
Let $G$ be a nice commutative group over an affine base $S$. 
   Then, $\bB G$ (defined in the fppf sense as in Remark \ref{remark-segal-groupoid-assoc-group}) is qcqs 1-Artin and is smooth surjective over $S.$ Its diagonal morphism is fppf.
\end{lemma}
\begin{proof}
As in Remark \ref{remark-segal-groupoid-assoc-group}, $\bB G$ is a 1-geometric fppf derived stack, 
which we can view as a 1-Artin derived stack in this paper's sense. 
Since being flat and locally of finite presentation  is fppf local on the source, these properties make sense for geometric fppf derived stacks and can be checked on fppf atlases. 
The morphism $S\to \bB G$
(which exhibits $\bB G$ as the colimit of the diagram in Remark \ref{remark-segal-groupoid-assoc-group}) is an fppf atlas by the proof of $((3)\implies (1))$ in \cite[Proposition 1.3.4.2]{HAG2}.
Then, using this atlas shows that the morphism $\bB G\to S$ is flat and locally of finite presentation as a 1-geometric derived stack (because $\id_S$ is). The morphism $\bB G\to S$ is still flat and locally of finite presentation if we view it as a 1-Artin derived stack. 
Choose a smooth atlas $U\to \bB G$ by a $0$-geometric derived stack. We hope to show $U\to S$ is smooth surjective. By above, this morphism is flat and locally of finite presentation, so $U$ is classical, and we may check smoothness and surjectivity fiberwise on $S$. Thus, we reduce to the case where $S$ is the spectrum of a field. 
The result will follow once we check smoothness and surjectivity of some particular smooth atlas of $\bB G$ (over the field). This can be done by using a smooth surjection $\GL/G\to \bB G$ as in \cite[Example 2.2]{HKR-charp}. More precisely, the  quotient $\GL/G$ is a scheme 
by \cite[Remarque 9.3(b)]{SGA3-Expose-5-6} and is smooth over the base field. We can choose a Zariski cover of $\GL/G$ by a 0-geometric derived scheme  $U$, and the composition $U\to \GL/G\to \bB G$ is a smooth atlas. $U$ is smooth surjective over the base field, so this shows smoothness and surjectivity. 

As for the fact that $\bB G$ has an affine diagonal, 
since being affine is fppf local on the target, it suffices to check the base change of $\Delta_{\bB G}$ to some fppf atlas is affine. 
We can choose $S\times_S S\to \bB G\times_S\bB G$ as the atlas, and we get $(S\times_S S)\times_{\bB G\times_S \bB G}\bB G\simeq S\times_{\bB G}S\simeq G$, which is affine. Using the same atlas, we also see that $\Delta_{\bB G}$ is fppf. 

To see $\bB G$ is quasi-compact (as a 1-Artin derived stack), we choose a smooth atlas $\phi:V\to \bB G$. $V$ is 0-Artin, so $V=\coprod_{i\in I}V_i$ for some derived affines $V_i$. 
On the other hand, we have the fppf atlas $S\to \bB G$. Consider the base change of the smooth atlas $\phi_S:\coprod_{i\in I} V_i\times_{\bB G}S\to S$ where $V_i\times_{\bB G}S$ is 0-geometric. Since $\phi_S$ is smooth and $S$ is affine, we can find a finite subset $J\subset I$ for which the morphism $\coprod_{i\in J} V_i\times_{\bB G}S\to S$ is still surjective. Then, $\coprod_{i\in J}V_i\to \bB G$ is a smooth atlas from a derived affine, so this shows quasi-compactness. 
\end{proof}
\begin{proposition}\label{proposition-BG-bundle-1-Artin-smooth-surj}
Let $G$ be a commutative nice group scheme over an affine $S$.
Let $f:X\to Y$ be a principal $\bB G$-bundle. Then, $f$ is 1-Artin (i.e.\ for any derived affine $T\to Y$, the base change $X\times_YT$ is 1-Artin) and smooth surjective (in the 1-geometric sense \cite[4.1.8]{DAG-gaitsgory-1}).
\end{proposition}
\begin{proof}
We first show 1-geometricity. Take a morphism $T\to Y$ from a derived affine. 
{ Structure morphisms of principal $\bB G$-bundles are effective epimorphisms by an axiom of $\infty$-topoi (\cite[Definition 2.1(3)]{NSS}), so $f_T:X\times_YT\to T$ is an effective epimorphism.}
As in \cite[Proposition 2.5]{NSS}, it induces an epimorphism of discrete sheaves $\tau_{\leq0}(X\times_YT)\to\tau_{\leq0}(T)$ in $\tau_{\leq0}(\dSt_{/T}).$ Hence, up to choosing an \'etale cover $\{ T_i\to T\}_i$ by derived affines $T_i$, 
the morphism $f_T$ admits a section.
This means that we have a morphism $\coprod_iT_i\to X\times_YT$ whose composition with $f_T$ commutes with the map $\coprod_i T_i\to T.$ The principal $\bB G$-bundle $f_{\coprod_i T_i}$ has a section, so it can be trivialized.
Since $T$ is affine, we can refine the cover to assume that the coproduct $\coprod_i$ is a finite one. 
By applying \cite[Proposition 4.3.6]{DAG-gaitsgory-1} to the Cartesian square
\[
\xymatrix{
\coprod_i T_i \times_S\bB G\ar[r]\ar[d]& X\times_YT\ar[d]\\
\coprod_i T_i\ar[r]&T
}
\]
we see that $X\times_YT$ is 1-Artin. 

As for smoothness and surjectivity, we check the morphism $X\times_YT\to T$ as above has these properties, which by definition is checked on an atlas of the source (\cite[4.1.7]{DAG-gaitsgory-1}). 
Take a smooth atlas $\pi:U\to X\times_YT$. The base change of $\pi$ along the affine morphism $\coprod_iT_i\to T$ is a smooth atlas of $\coprod_iT_i\times_S\bB G$, so Lemma \ref{lemma-BG-qcqs-1-artin-smooth-surj} shows that $U\times_T\coprod_iT_i\to \coprod_iT_i$ is smooth surjective. Since $\coprod_iT_i\to T$ is an \'etale cover, $U\to T$ is smooth surjective. 
\end{proof}
\begin{corollary}\label{corollary-BG-bundle-locally-split-in-topology}
In the situation of Proposition \ref{proposition-BG-bundle-1-Artin-smooth-surj}, we additionally assume that $Y$ is 1-Artin. 
Then,
    {
    \begin{enumerate}
        \item the principal $\bB G$ bundle $f$ is trivialized on a smooth atlas of $Y$, and 
        \item $f$ is quasi-compact.
    \end{enumerate}
    }
\end{corollary}
\begin{proof}
    {The projection $X\times_YX\to X$ has a section $\Delta_{X/Y}$, so this map is equivalent to the projection $X\times\bB G\to X$ (as principal $\bB G$-bundles).
    By the proposition, $f$ is 1-Artin, so $X$ is 1-Artin. 
    Choose a smooth atlas from a 0-geometric $U\to X$. Then, the composition $U\to X\to Y$ is smooth surjective and hence an atlas. The base change of $f$ along this composition is identified with the trivial  principal $\bB G$-bundles $U\times\bB G\to U$, so we obtain the result.
    }
    
    For quasi-compactness, consider any morphism from a derived affine $T\to Y$ and we need to show $X\times_YT$ is quasi-compact as a 1-Artin derived stack. By the previous part, we can find a smooth atlas from a 0-geometric $U\to Y$ for which $f_U:X\times_YU\to U$ can be identified with the projection $\bB G\times_SU\to U$. By Lemma \ref{lemma-BG-qcqs-1-artin-smooth-surj}, this morphism is quasi-compact. $Y$ has a 0-Artin diagonal, so $T\times_YU$ is 0-geometric. Write $T\times_YU=\coprod_{i\in I} V_i$ as a disjoint union of derived affines. The morphism $T\times_YU\to T$ is smooth surjective and $T$ is derived affine, so we can find a finite subset $J\subset I$ such that the composition $V:=\coprod_{i\in J}V_i\to \coprod_{i\in I}V_i=T\times_YU\to T $ surjects. Then, we have the cartesian squares:
    \[
    \xymatrix{
X\times_YV\ar[r]\ar[d]^-{f_{V}}&    X\times_YT\times_YU\ar[r]\ar[d]^-{f_{T\times_YU}}&X\times_YT\ar[d]^-{f_T}\\
V\ar[r]&T\times_YU\ar[r]&T
    }
    \]
    $V$ is affine and $f_U$ is quasi-compact, so $X\times_YV$ is quasi-compact. We can then choose a smooth atlas $W\to X\times_YV$ from a derived affine $W$.
    The composition $W\to X\times_YV\to X\times_YT$ is smooth surjective, so it is an atlas, and we see that $X\times_YT$ is quasi-compact.
\end{proof}
We also note a variant where we  trivialize principal $\bB G$-bundle \'etale locally. 
\begin{proposition}\label{proposition-trivialize-BG-bundle-etale-locally}
Let $G$ be a nice commutative group scheme over an affine $S$. Let $f:X\to Y$ be a principal $\bB G$-bundle over a 1-Artin derived stack $Y$ which admits an \'etale atlas (i.e.\ the atlas morphsm is 0-Artin \'etale). Then, the principal $\bB G$-bundle $f$ can be trivialized after base change along an \'etale atlas $V\to Y$. 
\end{proposition}
\begin{proof}
Let $U\to Y$ be an \'etale atlas. Write $U=\coprod_iU_i$ as a disjoint union of derived affines. 
The base changed morphisms $f_{U_i}:X\times_YU_i\to U_i$ are principal $\bB G$-bundles and in particular an effective epimorphisms. By an argument similar to the proof of Proposition \ref{proposition-BG-bundle-1-Artin-smooth-surj}, for each $i$, we can find an \'etale cover (by derived affines) $\{V_{i,j}\to U_i\}_j$  such that $f_{U_i}$ admits a section after passing to this cover. Then, the morphism $f_U$ admits a section up to passing to the \'etale cover
$V:=\coprod_{i,j}V_{i,j}\to \coprod_iU_i$. 
The morphism $V\to Y$ achieves the result.
\end{proof}
\subsection{Principal $\bB G$-bundles in \'etale and fppf topologies}\label{subsection-principal-BG-bundle-etale-fppf}
In this paper, we will be considering principal $\bB G$-bundles (where $G$ is nice and commutative) in the $\infty$-category of derived stacks in the \'etale topology. One can also consider an analogue in the fppf topology. In this subsection we briefly clarify the relation between these two. 
Let $S$ be an affine, and let $G$ be a commutative nice group over $S$. As before, $\bB G$ is the fppf delooping of $G$, viewed as a 1-Artin derived stack. In this subsection, for clarity, we write it as $i\bB_{\fppf}G$, where $i$ is the inclusion functor from the $\infty$-category of fppf derived stacks to  the $\infty$-category of \'etale derived stacks (which preserves geometric derived stacks by Toen's result, as before). $a$ is the fppf sheafification functor which is the {left} adjoint of $i$ (which also respects geometric derived stacks). 
Note that $i\bB_{\fppf}G$ is still an $\infty$-group as $i$ respects fiber products. 
\begin{proposition}\label{B2G-comparison-etale-fppf}
    We have $\bB_{\et}i\bB_{\fppf} G\simeq i\bB_{\fppf}^2 G$.
\end{proposition}
\begin{proof}
Both of these objects are the colimits of
 \[
    \xymatrix{
    ... \ar@<0.5ex>[r]\ar@<-0.5ex>[r]\ar@<1.5ex>[r]\ar@<-1.5ex>[r]&\bB_{\fppf} G\times\bB_{\fppf} G \ar@<1ex>[r]\ar[r]\ar@<-1ex>[r] & \bB_{\fppf} G \ar@<0.5ex>[r]\ar@<-0.5ex>[r]  & S
    }
    \] 
but computed in different topologies (where we apply $i$ appropriately). The colimit for the left object is computed in the $\infty$-category of \'etale derived stacks, and the right one is computed in its fppf analogue. Then, it's enough to show that the colimit in the \'etale sense is already an fppf derived stack (so that the sheafification yields an equivalent object). 
By \cite{toen-descent}, this follows once we know the colimit in the \'etale sense is already Artin. In the above diagram, the morphism $\bB_{\fppf}G\to S$ is 1-Artin smooth surjective by Proposition \ref{proposition-BG-bundle-1-Artin-smooth-surj}, so by \cite[Proposition 1.3.4.2]{HAG2}, the colimit is 2-Artin. 
\end{proof}
\begin{corollary}
    For any $n$-Artin derived stack $X$, we have \[
    \Bund_{X,\et}(i\bB_{\fppf}G)
    \simeq
    \Bund_{a(X),\fppf}(\bB_{\fppf}G),\]
    where the subscripts $\et,\fppf$ denote which topologies the principal bundles are defined in.
    In particular, the set of equivalence classes of
    principal $\bB G$-bundles (in the convention of this paper) are classified by $\rH^2_{\fppf}(X,G)$.
\end{corollary}
\begin{proof}
    This is immediate from the equivalence $\dArtSt^m\simeq \dGeomSt_{\fppf}^m$ (for each $m$) of \cite{toen-descent}, as it implies
    \[
    \Map_{\dGeomSt_{\fppf}^m}(X, \bB_{\fppf}^2 G)\simeq 
    \Map_{\dArtSt^m}(X, \bB_{\et}i\bB_{\fppf} G)
    \] with $m=\max(n,2).$
\end{proof}
\begin{remark}\label{remark-smooth-group-scheme-Grothendieck}
    By an argument similar to Proposition \ref{B2G-comparison-etale-fppf}, when $G$ is additionally smooth, $\bB_{\fppf}G\simeq\bB_{\et}G$, which is essentially Grothendieck's theorem.
\end{remark}
\section{Scalloped derived stacks}
For a qcqs 1-Artin derived stack $X$, recall that a \textbf{scallop decomposition} (see \cite[Definition 2.7]{Khan-Ravi}) is a finite filtration $\emptyset=U_0\hookrightarrow U_1\hookrightarrow...\hookrightarrow U_n=X$ by quasi-compact open substacks along with ``Nisnevich'' cartesian squares
\[
\xymatrix{
W_i\ar[r]\ar[d]&V_i\ar[d]^-{u_i}\\
U_{i-1}\ar[r]^-{j_i}& U_i
}
\] where $j_i$ are the open immersions, $u_i$ are representable \'etale morphisms of finite presentation inducing isomorphisms away from $U_{i-1}$ (where representability of morphisms in this section means representability by derived algebraic spaces). 
Observe in this definition that the $U_i$'s are automatically qcqs.  
We also note that the morphisms $j_i$ are automatically quasi-compact since the $U_{i-1}$'s are quasi-compact and the $U_i$'s are quasi-separated. 

\begin{definition}[{\cite[Definition 2.9]{Khan-Ravi}}]
A qcqs 1-Artin derived stack is called \textbf{(nicely) scalloped} if it admits a scallop decomposition $(U_i,V_i,u_i)_i$  with (nicely) quasi-fundamental derived stacks $V_i$ (i.e.\ there exists a quasi-affine morphism $V_i\to\bB G_i$ for some nice embeddable group scheme $G_i$ over a classical affine $S_i$). The $\infty$-category of scalloped (qcqs) 1-Artin derived stacks is denoted by $\dArtSt^{\scalloped}$.
\end{definition}
\begin{remark}
    In the above definition, separated diagonal condition from \cite{Khan-Ravi} is dropped as it's automatic in our convention.
\end{remark}
For our purposes, the following characterizations of scalloped derived stacks are useful:
\begin{theorem}[{\cite[Theorem 2.12 (ii),(iii)]{Khan-Ravi}}]\label{theorem-khan-ravi-scalloped-iff-nicely-stabilized}
    Let $X$ be a qcqs 1-Artin derived stack (with an affine diagonal). Then, the following conditions are equivalent.
    \begin{enumerate}
        \item $X$ is scalloped.
        \item $X$ has nice stabilizers.
        \item $X$ admits a relatively affine Nisnevich cover (in the sense of morphisms relatively representable by derived algebraic spaces) $U\to X$ where $U$ is (nicely) fundamental (i.e.\ it admits an affine morphism to $\bB G$ for some nice embeddable group scheme $G\to S$ over a classical affine $S$). 
    \end{enumerate}
\end{theorem}

\begin{proposition}\label{principal-BG-bundle-over-scalloped-is-scalloped}
Let $G$ be a nice commutative group scheme over an affine $S$. 
    Let $Y$ be a scalloped qcqs 1-Artin derived stack over $S$, and let $f:X\to Y$ be a principal $\bB G$-bundle from some derived stack $X$. Then, $X$ is a scalloped qcqs 1-Artin derived stack.
\end{proposition}
\begin{proof}
By Proposition \ref{proposition-BG-bundle-1-Artin-smooth-surj}, $X$ is 1-Artin.
 By Corollary \ref{corollary-BG-bundle-locally-split-in-topology}, it is quasi-compact as well. 
The diagonal of $X$ can be written as the composition $X\to X\times_YX\to X\times X$. 
The second morphism is the base change of $\Delta_Y:Y\to Y\times Y$ along $X\times X\to Y\times Y$, so it's affine as $Y$ has affine diagonals. 
Affineness of morphisms can be checked fppf locally on target, so we may verify the affineness of $\Delta_{X/Y}$ by replacing $Y$ by a smooth atlas where the principal bundle is trivial, using Corollary \ref{corollary-BG-bundle-locally-split-in-topology}.
In this special case, the relative diagonal is affine by the proof of Lemma \ref{lemma-BG-qcqs-1-artin-smooth-surj}. Hence, $X$ has an affine diagonal. For a later use, note that the same argument and Lemma \ref{lemma-BG-qcqs-1-artin-smooth-surj} imply that the morphism $\Delta_{X/Y}$ is fppf.

By Theorem \ref{theorem-khan-ravi-scalloped-iff-nicely-stabilized}, in order to show $X$ is scalloped, we check it has nice stabilizers.
Choose a geometric point $p:\Spec L\to X$. Let $q$ be its target in $Y$. Choose a smooth atlas $U\to Y$ trivializing the principal $\bB G$-bundle (Corollary \ref{corollary-BG-bundle-locally-split-in-topology}). Up to extending the field $L$ by a bigger algebraically closed field (which does not affect niceness of stabilizers by Lemma \ref{lemma-nice-group-field-extension}), we can lift $q$ to a geometric point $r$ on $U.$ Then, we can lift it to a point $s$ on the fiber product $X\times_YU$ as well.
In view of the cartesian square below, the morphism $\Aut_p(X_{\cl})\to\Aut_q(Y_{\cl})$ is fppf, so it's a surjection of group schemes.
\[
\xymatrix{
\Aut_p(X_{\cl})\ar[r]\ar[d]&\Aut_q(Y_{\cl})\ar[d]\\
X_{\cl}\ar[r]^-{\Delta_{X_{\cl}/Y_{\cl}}}&(X\times_YX)_{\cl}
}
\]
Combining this with
\[
\Aut_s((X\times_YU)_{\cl})\simeq\Aut_p(X_{\cl})\times_{\Aut_q(Y_{\cl})}\Aut_r(U_{\cl})\simeq \ker(\Aut_p(X_{\cl})\to\Aut_q(Y_{\cl})),
\] we obtain a short exact sequence
\[
1\to\Aut_s((X\times_YU)_{\cl})\to\Aut_p(X_{\cl})\to \Aut_q(Y_{\cl})\to 1.
\]
We note $\Aut_s((X\times_YU)_{\cl})\simeq\Aut_s(U\times_S\bB G)\simeq G$, so it reads as
\[
1\to G\to\Aut_p(X_{\cl})\to \Aut_q(Y_{\cl})\to 1.
\]
$G$ and $\Aut_q(Y_{\cl})$ are nice, so $\Aut_p(X_{\cl})$ is nice by Remark \ref{remark-nice-group-closed-under-extension-subgroup-quotient}.
\end{proof}
This proposition produces the following important  examples of scalloped derived stacks.
\begin{example}
Let $G$ be a nice commutative group over a field $k$, and let $Y$ be a separated variety. If $f:X\to Y$ is a principal $\bB G$-bundle, $X$ is qcqs scalloped 1-Artin. 
\end{example}

\begin{proposition}[variant of {\cite[Theorem 3.14]{Khan-Ravi}}]\label{proposition-khan-ravi-scalloped-lift-cover}
    Let $Z\to X$ be a closed immersion of qcqs scalloped 1-Artin derived stacks. Consider a smooth morphism $Z'\to Z$ (in the sense of morphisms representable by derived algebraic spaces). Then, up to replacing $Z'$ by a relatively affine Nisnevich cover, the morphism $Z'\to Z$ arises as a base change of some smooth morphism $X'\to X$ (in the sense of morphisms representable by derived algebraic spaces).
\end{proposition}
\begin{proof}
    Our qcqs 1-Artin derived stacks have affine diagonals by convention, so by Theorem \ref{theorem-khan-ravi-scalloped-iff-nicely-stabilized}, we can choose a relatively affine Nisnevich cover $U\to X$ where $U$ is (nicely) fundamental. Note that $U\times_XZ$ is also (nicely) fundamental. After base changing the whole diagram along $U\to X$, the result follows from \cite[Theorem 1.13]{Alper-Hall-HL-Rydh}; we can find a smooth morphism $U'\to U$ lifting $U\times_XZ'\to U\times_XZ$. Thus, the result holds up to replacing $Z'\to Z$ by $U\times_XZ'\to Z'\to Z$.
\end{proof}

\begin{proposition}[{\cite[Theorem 2.24]{Khan-Ravi}}]\label{proposition-khan-ravi-scalloped-perfect}
    Every scalloped derived stack $X$ is perfect (i.e.\ $\Dqc(X)$ is compactly generated by $\Dperf(X)$).
\end{proposition}

\begin{lemma}\label{lemma-scalloped-object-compact-perfect-dualizable}
    Let $X$ be scalloped (qcqs 1-Artin) derived stack. Then, dualizable objects, perfect objects and compact objects all agree.
\end{lemma}
\begin{proof}
The result follows from
    {\cite[Propositions 3.6, 3.9]{bzfn}} as scalloped derived stacks are perfect by  Proposition \ref{proposition-khan-ravi-scalloped-perfect}.
\end{proof}
For later purposes, we record the following variant:
\begin{proposition}\label{proposition-tamely-compactly-generated}
    Let $X$ be a acalloped (qcqs 1-Artin) derived stack. Then, $\Dqc(X)$ is tamely compactly generated (i.e.\ compactly generated, and mapping spectra between compact objects are almost connective).
\end{proposition}
\begin{proof}
We already know compact generation (Proposition \ref{proposition-khan-ravi-scalloped-perfect}), so we choose perfect objects $E,F\in\Dperf(X)$ and show $\RHom_X(E,F)$ is almost connective. 
Since $X$ is scalloped, it admits a scallop decomposition $\emptyset=U_0\hookrightarrow U_1\hookrightarrow...\hookrightarrow U_n=X$ with appropriate properties,
and we prove inductively that $\RHom_{U_i}(E|_{U_i},F|_{U_i})$ is almost connective by induction on $i$. To complete the induction steps, we will show the following statement. 
$X$ and $E,F\in\Dperf(X)$ are as before, and we assume we have a ``Nisnevich'' cartesian square
\[
\xymatrix{
W\ar[r]\ar[d]&V\ar[d]^-p\\
U\ar[r]^-j&X
}
\]
where $j$ is a quasi-compact open immersion, and $p$ is a representable \'etale morphism inducing an isomorphism away from $U$. 
Assume $V$ is (nicely) fundamental and $\RHom_U(E|_U,F|_U)$ is almost connective.
As in \cite[diagram (2.28)]{Khan-Ravi}, we have a pullback square of stable $\infty$-categories
\[
\xymatrix{
\Dperf(X)\ar[r]\ar[d]&\Dperf(U)\ar[d]\\
\Dperf(V)\ar[r]&\Dperf(W)
}
\]
In particular, we have 
\[
\RHom_X(E,F)\simeq\RHom_U(E|_U,F|_U)\times_{\RHom_W(E|_W,F|_W)}\RHom_V(E|_V,F|_V).
\]
By base change, the morphism $W\to V$ is quasi-affine, and $V$ is (nicely) quasi-fundamental, so $W$ is also (nicely) quasi-fundamental.
In the proof of \cite[Lemma 2.26]{Khan-Ravi}, they show that (nicely) quasi-fundamental derived stacks universally have finite cohomological dimensions.
We will use this to show
$\RHom_W(E|_W,F|_W)\simeq\rR\Gamma(W,E|_W^{\vee}\otimes F|_W)$ and $\RHom_V(E|_V, F|_V)\simeq\rR\Gamma
(V,E|_V^{\vee}\otimes F|_V)$ 
are almost connective.
The perfect complexes $E|_W^{\vee}\otimes F|_W$ and $E|_V^{\vee}\otimes F|_V$ are almost connective which we can check after passing to smooth atlases by derived affines using quasi-compactness and
\cite[Proposition 1.5.4 pp.150]{DAG-gaitsgory-1}. 
Then, they remain almost connective after applying $\rR\Gamma(W,-)$ or $\rR\Gamma(V,-)$ by the property of finite cohomological dimensions.
Combining this with the assumption on $\RHom_U(E|_U,F|_U)$, we see $\RHom_X(E,F)$ is also almost connective.
\end{proof}

\begin{proposition}\label{proposition-scalloped-dstack-stable-fiber-product}
    The $\infty$-category $\dArtSt^{\scalloped}$ is closed under fiber products.  
\end{proposition}
\begin{proof}
It's standard that  fiber products exist in $\dArtSt^{\qcqs}$, and we show that this fiber product is scalloped using Theorem \ref{theorem-khan-ravi-scalloped-iff-nicely-stabilized}.
Consider a fiber product $X\times_YZ$ of $X,Y,Z\in\dArtSt^{\scalloped}$ in $\dArtSt^{\qcqs}$.
Given a geometric point $p\in (X\times_YZ)(k)$, the stabilizer at $p$, $\Aut_p(X\times_YZ)_{\cl}$, can be written as $\Aut_p(X_{\cl})\times_{\Aut_p(Y_{\cl})}\Aut_p(Z_{\cl})$, where we are abusing notation to write $p$ for the targets of the geometric point $p$ on $X,Y,Z$. 
By assumption, \[\Aut_p(X_{\cl}),\Aut_p(Y_{\cl}),\Aut_p(Z_{\cl})\] are nice group schemes over $k$, so we are done by Lemma \ref{lemma-nice-groups-fiber-products}.
\end{proof}

\begin{lemma}\label{lemma-nice-group-derived-pushforard-structure-sheaf} Let $S$ be an affine scheme. 
    Let $G$ be a nice commutative group over $S$. Then, any principal $\bB G$-bundle $f:X\to Y$ of qcqs 1-Artin derived stacks over $S$ is perfect as a morphism (\cite[Definition 3.2]{bzfn})
    and 
    satisfies $\cO_Y\simeq f_*\cO_X$ under the natural map.  
\end{lemma}
\begin{proof}
We first prove that $f$ is perfect. This means that for any morphism $U\to Y$ from a derived affine, $U\times_YX$ is a perfect derived stack. Any such fiber product $U\times_YX$ is qcqs 1-Artin and a principal $\bB G$-bundle over a derived affine $U$. Since $G$ is nice, it has nice stabilizers. 
By Theorem \ref{theorem-khan-ravi-scalloped-iff-nicely-stabilized}, $U\times_YX$ is scalloped. By Proposition \ref{proposition-khan-ravi-scalloped-perfect}, $U\times_YX$ is perfect. Hence $f$ is perfect.

By the above perfectness, $f$ satisfies base change for pushforward (\cite[Proposition 3.10]{bzfn}). This allows us to check {$\cO_Y\simeq f_*\cO_X$} smooth-locally, and we reduce to the case where $f$ is a trivial principal $\bB G$-bundle over a derived affine.
As in Remark \cite[Remark 2.3]{Khan-Ravi}, nice groups are linearly reductive (\cite{Alper}), so we are done.
\end{proof}

\section{Geometric pushouts of qcqs 1-Artin derived stacks}
\begin{definition}[{see \cite[Definition 4.1]{Alper-Hall-HL-Rydh}}]\label{definition-geometric-pushout}
The following diagram of qcqs 1-Artin derived stacks is called a \textbf{geometric pushout} if the following conditions hold.
    \[
\xymatrix{
X_{0}\ar[r]^-{a_2}\ar[d]^-{a_1}&X_2\ar[d]^-{b_2}\\
X_1\ar[r]^-{b_1}& X_3
}
    \]
    \begin{enumerate}
        \item All morphisms are affine. 
        \item $a_1,a_2$ are closed immersions.
        \item The morphism \[
        \cO_{X_3}\to (b_1)_*\cO_{X_1}\times_{(b_1a_1)_*\cO_{X_{0}}}(b_2)_*\cO_{X_2}
        \] is an equivalence. 
    \end{enumerate}
    We will also call $X_3$ a \textbf{geometric pushout of }$X_1\leftarrow X_0\to X_2$ in this situation.
    If $a_1,a_2$ are additionally nilpotent closed immersions (on the classical locus), a geometric pushout is called a \textbf{nilpotent geometric pushout.} 
\end{definition}
\begin{remark}\label{remark-geometric-pushout-closed-immersions}
    In the above definition, $b_1,b_2$ are automatically closed immersions. This boils down to the following classical fact: given a diagram of classical rings $B\to A\leftarrow C$ with surjective $C\to A$, the morphism $B\times_AC\to B$ is surjective. 
  \end{remark}
\begin{remark}\label{remark-geometric-pushout-smooth-local-on-X3}
Similarly to  {\cite[Lemma 4.3]{Alper-Hall-HL-Rydh}}, geometric pushouts behave well under base changes. Namely, base changing a geometric pushout as in Definition \ref{definition-geometric-pushout} along $X_3$ gives us another geometric pushout (note that relevant pushforwards  satisfy base change as the morphisms are affine; cf.\ \cite[Proposition 3.10]{bzfn}, Lemma \ref{lemma-nonconnective-base-change-projection-formula}). Also, if a square diagram of qcqs 1-Artin derived stacks is a geometric pushout smooth-locally on the bottom right term, the original diagram is a geometric pushout. 
\end{remark}

\begin{lemma}[see {\cite[Lemma 4.4]{Alper-Hall-HL-Rydh} for classical stacks}]\label{lemma-pushout-cube-faces-pullbacks}
Consider a diagram of qcqs 1-Artin derived stacks below where the top and bottom horizontal faces are geometric pushouts. Assume that the two vertical faces, containing $a_1,b_1$ and $a_2,b_2$ respectively, are cartesian.
    \begin{equation}
\begin{tikzcd}[row sep={3.5em,between origins}, column sep={5em,between origins}]
    &U_0 \arrow[rr, "b_2"] \arrow[dd] \arrow[dl, "b_1"'] & & 
    U_2 \arrow[dd,"f_2"] \arrow[dl] \\
    U_1 \arrow[rr] \arrow[dd,"f_1"] & & 
    U_3 \arrow[dd,"f_3" near end] \\
    & X_0\arrow[rr, "a_2" near end] \arrow[dl, "a_1"'] & & 
    X_2 \arrow[dl] \\
    X_1 \arrow[rr] & & 
    X_3
\end{tikzcd}
\end{equation}
Then, 
\begin{enumerate}
    \item all vertical faces are cartesian, and
    \item $f_3$ is smooth if and only if $f_1$ and $f_2$ are.
\end{enumerate}

\end{lemma}
\begin{proof}
\leavevmode
\begin{enumerate}
    \item 
This is proved in the same as \cite[Lemma 4.4]{Alper-Hall-HL-Rydh}. Since geometric pushouts are stable under base change and cartesianness is smooth-local, we can localize on $U_3$ and $X_3$ and assume these are derived affine. Thus, everything is affine. 
We will need to replace the use of \cite[Thm. 2.2(i)]{Ferrand} in the original proof by its derived analogue. 

We write $X_i=\Spec A_i,U_i=\Spec B_i$ for each $i$. We have $A_3\simeq A_1\times_{A_0}A_2$ and similarly for the $B_i$'s. The goal is to show $B_i\simeq B_3\otimes_{A_3}A_i$ for $i=0,1,2$ under the canonical maps, which we will verify as modules instead of derived rings. 
We recall that we have an equivalence
\[
F:\Mod_{A_3,\geq0}\to \Mod_{A_1,\geq0}\times_{\Mod_{A_0,\geq0}}\Mod_{A_2,\geq0}
\]
as in \cite[Theorem 16.2.0.2, Proposition 16.2.2.1]{SAG}, {\cite[Lemma 3.3]{Moore2}}. 
The functor $F$ is the operation of taking tensor products. The inverse $G$ is the right adjoint of $F$, which takes the object $(M_1,M_2,M_1\otimes_{A_1} A_0\simeq M_0\simeq M_2\otimes_{A_2} A_0)$ with $M_i\in\Mod_{A_i,\geq0}$ in the fiber product category to $M_1\times_{M_0}M_2.$

This formalism applies in particular to 
$(B_1,B_2,B_1\otimes_{A_1}A_0\simeq B_0\simeq B_2\otimes_{A_2}A_0)$ where we view the $B_i$'s as modules. This object is equivalent to
\[(F\circ G)(B_1,B_2,B_1\otimes_{A_1}A_0\simeq B_0\simeq B_2\otimes_{A_2}A_0)
\simeq F(B_1\times_{B_0}B_2)\simeq F(B_3).\]
This says precisely that $B_0,B_1,B_2$ arise from base changing $A_0,A_1,A_2$ along $A_3\to B_3$.  
\item 
Only if direction is clear from the previous part, so we assume $f_1$ and $f_2$ are smooth and prove $f_3$ is. We can prove smoothness smooth-locally on $U_3$ and $X_3$, so again we assume everything is affine. It suffices to show that $f_3$ is locally of finite presentation and the relative cotangent complex is perfect of Tor-amplitude [0,0]. 
By the derived affine case of \cite[Proposition 9.3]{DAG-IX}, 
$f_3$ is locally of finite presentation. Using the base change of cotangent complexes as well as the affine case of \cite[Proposition 7.7]{DAG-IX}, 
we see that the relative cotangent complex of $f_3$ is perfect of the correct Tor-amplitude, so we see $f_3$ is smooth.
\end{enumerate}
\end{proof}
\begin{proposition}\label{proposition-existence-geometric-pushout}
Consider a diagram of scalloped (qcqs 1-Artin) derived stacks:
\[
Y_1\leftarrow Y_0\to Y_2.
\]
Assume the arrows are nilpotent closed immersions. Then, it admits a scalloped (qcqs 1-Artin) nilpotent geometric pushout $Y_3$. 
\end{proposition}
\begin{proof}
    This is clear in the derived affine case (note that when one $Y_i$ is affine, all three derived stacks are affine by Lemma \ref{lemma-affineness-derived-vs-classical}).

     Using quasi-compactness, we first take a ($0$-geometric) smooth atlas of $Y_1$ by a derived affine: $U_1\to Y_1$. By base change along the closed immersion $Y_{0}\to Y_1$, we obtain a smooth atlas with the same properties $U_{0}\to Y_{0}.$
    By Lemma \ref{proposition-khan-ravi-scalloped-lift-cover}, we can find a relatively affine Nisnevich cover $V_0\to U_{0}$ 
    and a cartesian diagram 
    \[
    \xymatrix{
V_0\ar[r]\ar[d]&V_2\ar[dd]^-{f_2}\\
U_{0}\ar[d]&
\\
    Y_{0}\ar[r]&Y_2
    }
    \] where $f_2$ is smooth (in the representable morphism sense). 
    
    Since $U_{0}$ is affine, so is $V_{0}$.
    As the Nisnevich cover $V_0\to U_0$ is in particular \'etale, 
     by Lemma \ref{lemma-infinitesimal-invariance-etale},
    we can descend it to an \'etale morphism of affines $V_1\to U_1$.
    The morphism $V_0\to U_0\to Y_0$ is a relatively 0-geometric morphism which 
    is surjective and smooth in the 0-geometric sense. 
    $V_0\to V_2$ is a nilpotent closed immersion since $Y_0\to Y_2$ is. Then,
    by Lemma \ref{lemma-affineness-derived-vs-classical}, $V_2$ is affine. 
    We also see using nilpotent closed immersions that $f_1$ and $f_2$ are surjective in the 0-geometric sense (note that these morphisms are automatically 0-geometric as $V_1,V_2$ are affine and $Y_1,Y_2$ have 0-geometric diagonals).  
    Thus, $V_i\to Y_i$ is a smooth-atlas for each $i=0,1,2$, and the following diagram has pullback squares: 
   \[
\xymatrix{
V_1\ar[d]&V_{0}\ar[d]\ar[r]\ar[l]&V_2\ar[d]\\
Y_1&Y_{0}\ar[r]\ar[l]&Y_2
}
    \] 
    For $*=0,1,2,$ we consider the Cech nerve $C^\bullet_*:=C_\bullet( V_*/Y_*).$ The $Y_*$'s have affine diagonals, so the Cech nerves $C^\bullet_*$ consist of derived affines for $*=0,1,2$. 
    Levelwise, we can take the geometric pushout of $C^\bullet_1\leftarrow C_{0}^\bullet\to C_2^\bullet$, and we write these pushouts as $C^\bullet$, which forms a simplicial derived affine scheme. We claim that this satisfies the Segal groupoid condition (cf.\ \cite[Definition 1.3.1.6]{HAG2}). For instance, we check $(\sigma_0,\sigma_1):C^2\simeq C^1\times_{C^0}C^1$ (where $\sigma_i$ comes from the map $[1]\to [2]$ sending $0\mapsto i,1\mapsto i+1$). The corresponding statement holds for each $C^\bullet_*$ ($*=0,1,2$) since they are Cech nerves.
    Explicitly, 
    the above isomorphism comes from
    \[
V_*\times_{Y_*}V_*\times_{Y_*}V_*\simeq (V_*\times_{Y_*} V_*)\times_{\pr_2,V_*,\pr_1}(V_*\times_{Y_*} V_*).
    \]
    
    Moreover, by construction of the compatible covers, such equivalences are compatible under pullbacks. 
    We also note by Lemma \ref{lemma-pushout-cube-faces-pullbacks} that all vertical faces in the following diagram are cartesian.
    \begin{equation}
\begin{tikzcd}[row sep={3.5em,between origins}, column sep={5em,between origins}]
    &C_0^2 \arrow[rr, "f_2"] \arrow[dd,shift right]\arrow[dd,shift left]  \arrow[dl, "f_1"'] & & 
    C^2_2  \arrow[dd,shift right]\arrow[dd,shift left] \arrow[dl] \\
    C^2_1 \arrow[rr]  \arrow[dd,shift right]\arrow[dd,shift left] & & 
    C^2  \arrow[dd,shift right]\arrow[dd,shift left] \\
    &C^1_0 \arrow[rr, ""]  \arrow[dd,shift right]\arrow[dd,shift left]\arrow[dl, ""'] & & 
    C^1_2 \arrow[dd,shift right]\arrow[dd,shift left] \arrow[dl] \\
    C^1_1 \arrow[rr]  \arrow[dd,shift right]\arrow[dd,shift left] & & 
    C^1  \arrow[dd,shift right]\arrow[dd,shift left] \\
    & C^0_0\arrow[rr, "" near end] \arrow[dl, ""'] & & 
    C^0_2 \arrow[dl] \\
    C^0_1 \arrow[rr] & & 
    C^0 
\end{tikzcd}
\label{diagram-pushout-cech-nerve-Segal-groupoid}
\end{equation}
Here, vertical maps are given by $\sigma_0,\sigma_1:C^2_*\to C_*^1$ and $d_0,d_1:C_*^1\to C_*^0$.

We can phrase our desired equivalence as 
 \[
(\sigma_0,\sigma_1):\AffPushout_* (V_*\times_{Y_*}V_*\times_{Y_*}V_*)\]
\[\simeq (\AffPushout_* (V_*\times_{Y_*} V_*))\times_{\AffPushout_*( V_*)}(\AffPushout_*( V_*\times_{Y_*} V_*)).
    \]
    where $\AffPushout_*( V_*)$ denotes the pushout of $V_{1}\leftarrow V_{0}\to V_{2}$ in the $\infty$-category of derived affines, and other instances of $\AffPushout_*$ are defined similarly. 
Everything is affine, {so as in the proof of Lemma \ref{lemma-pushout-cube-faces-pullbacks} we can use the equivalence of \cite[Theorem 16.2.0.2]{SAG} and verify that $(\sigma_0,\sigma_1)$ is an equivalence after passing to an appropriate fiber product category of modules. 
Concretely, we only need to show that base changing the diagram $V_{1}\leftarrow V_{0}\to V_{2}$ along $\AffPushout_* (V_*\times_{Y_*}V_*\times_{Y_*}V_*)\to \AffPushout_* (V_*)$ yields an equivalent diagram to the base change of the same diagram along 
\[(\AffPushout_* (V_*\times_{Y_*} V_*))\times_{\AffPushout_*( V_*)}(\AffPushout_*( V_*\times_{Y_*} V_*))\to \AffPushout_* (V_*)\] under the given map.
If we use that the vertical faces in the diagram (\ref{diagram-pushout-cech-nerve-Segal-groupoid}) are cartesian, this follows from the equivalences $(\sigma_0,\sigma_1):C^2_*\simeq C^1_*\times_{C^0_*}C^1_*$ for all $*$.} 
This way, we see that $C^\bullet$ satisfies the Segal groupoid condition. Moreover, the morphism $d_0:C^1\to C^0$ is smooth by Lemma \ref{lemma-pushout-cube-faces-pullbacks}. 
We claim it is also surjective. $d_0:C_*^1\to C_*^0$ is surjective for $*=0,1,2$, so it suffices to show $C_1^0\coprod C_2^0\to C^0$ is surjective. 
This follows from the following easy computation; given surjective maps of ordinary rings $A_1\to A_0\leftarrow A_2$, the morphism $\Spec A_1\coprod \Spec A_2\to \Spec (A_1\times_{A_0}A_2)$ is surjective. 
Then, the colimit of $C^\bullet$ is also a 1-geometric derived stack by \cite[Proposition 1.3.4.2]{HAG2}.
Let us write this colimit as $Y_3$. 
Note that $C^\bullet$ coincides with the Cech nerve of the morphism $C^0\to Y_3$ by an axiom of $\infty$-topos. 
Since the defining properties of geometric pushouts can be checked smooth-locally (Remark \ref{remark-geometric-pushout-smooth-local-on-X3}), we see that the square diagram formed by $Y_0,Y_1,Y_2,Y_3$ is a geometric pushout.

Since $C^0$ is affine, $Y_3$ is quasi-compact. 
Observe that 
the base change of $\Delta_{Y_3}:Y_3\to Y_3\times Y_3$ along $C^0\times C^0\to Y_3 \times Y_3$ is $C^1\simeq C^0\times_{Y_3}C^0\to C^0\times C^0$ which is affine. By \cite[Lemma 3.4.10]{DAG-VIII}, 
we conclude that $Y_3$ has an affine diagonal (hence quasi-separated). 

We lastly check that $Y_3$ is scalloped. We start by observing that the closed immersions $Y_1\to Y_3,Y_2\to Y_3$ are jointly surjective (in the 0-geometric sense) so that they are jointly surjective on classicalizations. 
This surjectivity of $Y_1\coprod Y_2\to Y_3$ follows from the corresponding surjectivity on the Cech nerves proved above. 

We will check that $Y_3$ is scalloped using Theorem \ref{theorem-khan-ravi-scalloped-iff-nicely-stabilized}, so we choose a geometric point on $Y_3$ (or equivalently $(Y_3)_{\cl}$)
     and show the stabilizer group scheme of $(Y_3)_{\cl}$ at this point is nice. By the surjectivity above, the geometric point is in the image of the closed immersion $(Y_1)_{\cl}\to (Y_3)_{\cl}$ or $(Y_2)_{\cl}\to (Y_3)_{\cl}$.
   Relative diagonals of closed immersions of classical Artin stacks are isomorphisms, so closed immersions preserve stabilizer groups at geometric points. Thus, 
     we are done by the assumption that $Y_1,Y_2$ are scalloped.

\end{proof}

\begin{lemma}[{see \cite[Lemma 4.5]{Alper-Hall-HL-Rydh}} for classical stacks]\label{lemma-geometric-pushout-pushout-among-qcqs-Artin-stacks}
Geometric pushouts of qcqs 1-Artin derived stacks along closed immersions are pushouts among qcqs 1-Artin derived stacks.    
\end{lemma}
\begin{proof}
Our proof is inspired by \cite[Lemma 4.5]{Alper-Hall-HL-Rydh}, although our proof is simpler as relevant derived stacks have affine diagonals. Consider a geometric pushout diagram of qcqs 1-Artin derived stacks realizing $X_3$ as the geometric pushout of $X_1\leftarrow X_0\to X_2.$ Let $W$ be a qcqs 1-Artin derived stack. 
We have a natural morphism 
\[
F:\Map(X_3,W)\to \Map(X_1,W)\times_{\Map(X_0,W)}\Map(X_2,W).
\] and we need to show this is an equivalence.
We may prove this smooth-locally on $X_3$ (while base changing $X_0,X_1,X_2$ simultaneously), so we can assume $X_3$ is affine. More precisely, to do this, we take a smooth cover of $X_3$ from a derived affine. The Cech nerve of this map (as well as the base changes to $X_0,X_1,X_2$) consist of derived affines. Thus, if we can prove it for derived affines, by taking the colimit of the Cech nerves, we prove the general statement. 
This way, we assume the $X_i$'s are all affine. 

In this situation, we construct an inverse of $F$ as follows. A point on $\Map(X_1,W)\times_{\Map(X_0,W)}\Map(X_2,W)$  corresponds to the datum of $(\psi_1,\psi_2,\alpha)$ where $\psi_i:X_i\to W$ is a map for $i=0,1,2$ and 
$\alpha:\psi_1\circ j_1\simeq\psi_0\simeq  \psi_2\circ j_2$ is an equivalence if $j_i:X_0\to X_i$ denotes the given closed immersion for $i=1,2$. 
We choose a smooth cover $q:U\to W$ from a derived affine $U$. Since $W$ has an affine diagonal, base changing $q$ along the $\psi_i$'s gives us derived affines $U_i:=U\times_WX_i$. $U_1\leftarrow U_0\to U_2$ has an obvious pushout (among derived affines) which we construct by taking the fiber product of coordinate rings. Since the $U_i,X_i$'s are all affine, we have a map $U_3\to X_3$ making the cubical diagram of the $U_i,X_i$'s commute. By Lemma \ref{lemma-pushout-cube-faces-pullbacks}, the face of this cube containing all the $U_i$'s is the base change of the one with all the $X_i$'s. Also, since $U$ is affine, we obtain a map $\phi^0:U_3\to U$ by universal property among derived affines. The same construction works if we replace $U$ by different terms in the Cech nerve of $U\to W$. The geometric pushouts at different levels of Cech nerves are compatible and in fact they coincide with the Cech nerve of the morphism $U_3\to X_3$, since the cube in the earlier argument have cartesian vertical faces (similarly to the proof of Proposition \ref{proposition-existence-geometric-pushout}). 
Hence, we have a morphism $\phi_\bullet:C_\bullet(U_3/X_3)\to C_\bullet(U/W)$. Taking colimits, we obtain a morphism $\phi:X_3\to W$ satisfying $F(\phi)\simeq (\psi_1,\psi_2,\alpha)$. 
This construction is functorial, and it defines 
\[
G:\Map(X_1,W)\times_{\Map(X_0,W)}\Map(X_2,W)\to \Map(X_3,W).
\] By construction $G$ is an inverse of $F$ up to homotopy.
\end{proof}

\begin{proposition}[Scalloped derived stack analogue of {\cite[Theorem 16.3.0.1]{SAG}}]\label{proposition-schlessinger-criterion-scalloped-d-stack}
    Consider a pushout diagram of derived affine schemes:
    \[
    \xymatrix{
X_{0}\ar[r]\ar[d]&X_1\ar[d]\\
X_2\ar[r]& X_3
    }
    \]
    where $X_{0}\to X_i$  is a nilpotent closed immersion for $i=1,2$ (i.e.\ it is a nilpotent closed immersion on the classical locus).
    Then, pullback induces an equivalence
    \[
   \Phi: \dArtSt^{\scalloped}_{/X_3}\simeq  
\dArtSt^{\scalloped}_{/X_1}
\times_{\dArtSt^{\scalloped}_{/X_{0}}} 
\dArtSt^{\scalloped}_{/X_2}.
    \]
\end{proposition}
\begin{proof}
    We adapt the proofs of \cite[Theorem 16.3.0.1]{SAG} and \cite[Proposition 3.25]{BranTael} to the case of scalloped derived stacks.

    The objects in the codomain are represented by diagrams with cartesian squares
    \[
\xymatrix{
Y_1\ar[d]&Y_{0}\ar[d]\ar[r]\ar[l]&Y_2\ar[d]\\
X_1&X_{0}\ar[r]\ar[l]&X_2
}
    \]
    where $(Y_i\to X_i)\in\dArtSt^{\scalloped}_{/X_i}$ for $i=0,1,2$.
    By Proposition \ref{proposition-existence-geometric-pushout}, the top row of the above diagram admits a geometric pushout $Y_3$ (which is scalloped and qcqs 1-Artin). Moreover, $Y_3$ is a pushout in the category of qcqs 1-Artin derived stacks by Lemma \ref{lemma-geometric-pushout-pushout-among-qcqs-Artin-stacks}, so we have a morphism $Y_3\to X_3$, making the cubical diagram of the $X_i$'s and $Y_i$'s commute. This construction determines a functor 
    \[
\Psi:    \dArtSt^{\scalloped}_{/X_1}
\times_{\dArtSt^{\scalloped}_{/X_{0}}} 
\dArtSt^{\scalloped}_{/X_2}\to \dArtSt^{\scalloped}_{/X_3}
    \]
    as it's constructed using categorical pushouts. 
    In the above construction of $\Psi$, $Y_3$ is the pushout in the slice category $\dArtSt^{
    \scalloped}_{/X_3}$, so one easily checks that $\Psi$ is the left adjoint of $\Phi$.
    The counit $\Psi\circ \Phi\to\id$ is an equivalence by Remark \ref{remark-geometric-pushout-smooth-local-on-X3} and Lemma \ref{lemma-geometric-pushout-pushout-among-qcqs-Artin-stacks}.
    It remains to show that the unit map 
    $\id\to \Phi\circ \Psi$ is an equivalence. This is proved in the same way as \cite[Theorem 16.3.0.1]{SAG}; one can first reduce to the affine case using that geometric pushouts are smooth-local, and the affine computations are identical. 
\end{proof}
\section{Principal $\bB\bG_m$-bundles and twisted derived categories}
\subsection{Banding}\label{subsection-banding}
{Let $S$ be a classical affine, and $G$ denotes either $\bG_m$ or $\mu_n$ for some positive integer $n$ over $S.$
In this subsection, $\bH$ is the $\infty$-topos $\bH=(\dSt_{/S})_{/Y}$ for some Artin derived stack $Y$. Due to our convention on $\bB G$ (Remark \ref{remark-segal-groupoid-assoc-group}), in principle we need to use the fppf topology, but by Remark \ref{remark-smooth-group-scheme-Grothendieck}, the \'etale topology suffices for the case of $G=\bG_m$.}

We start by adapting the theory from \cite{NSS} to our situation. 
\begin{definition}[{\cite[Definition 4.54]{NSS}}]
When $G$ is an $\infty$-group in $\bH$ as above, we define an $\infty$-group
\[\Out(G):=\tau_{\leq0}\sAut_{\bH}(\bB G).\]
\end{definition}
\begin{definition}
    Let $K\to H\xrightarrow[]{q} Q$ be a diagram in $\Grp(\bH)$, and assume that $K$ is the fiber of $q$ (over the unit section). We say that this is a \textbf{short exact sequence of $\infty$-groups in} $\bH$ if $q$ is an effective epimorphism.
\end{definition}
By Remark \ref{remark-truncation-group} and \cite[Proposition 2.10]{NSS}, if $\tau_G:\sAut_{\bH}(\bB G)\to \Out(G)$ denotes the truncation morphism, we have a short exact sequence of $\infty$-groups:
\[
\fib(\tau_G)\to \sAut_{\bH}(\bB G)\to \Out(G).
\] 
\begin{example}\label{example-gerbe-banding-fiber-sequence}
\leavevmode
From Appendix \S\ref{section-infty-grp-dSt}, we have the following examples of fiber sequences constructed as above:
    \begin{enumerate}
        \item When $G=\bG_m$ (in the \'etale topology), the fiber sequence takes the form:
        \[
        \bB\bG_m\to \sAut_{\bH}(\bB\bG_m)\to \ZModTwo
        \]
        \item When $G=\mu_n$ (in the fppf topology), the fiber sequence takes the form:
        \[
        \bB\mu_n\to \sAut_{\bH}(\bB\mu_n)\to \underline{(\bZ/n\bZ)^\times}.
        \]
    \end{enumerate}
\end{example}
\begin{definition}\label{definition-banding-NSS-version}
    Let $X\to Y$ be a $G$-gerbe, which corresponds to a point $[X\to Y]\in\Map(Y,\bB\sAut_{\bH}(\bB G))$ via Proposition \ref{proposition-gerbe-BAutBG}. We define \textbf{the banding associated with $X\to Y$} to be the point of $\Map(Y,\bB\Out(G))$ which is the image of the point 
    \[[X\to Y]\in \Map(Y,\bB\sAut_{\bH}(\bB G))\] in $\Map(Y,\bB\Out(G))$ via $\tau_G.$ 
\end{definition}
\begin{definition}
    We define the \textbf{$\infty$-groupoid $\Gerbe_Y^{\triv}(G)$ of $G$-gerbes over $Y$ with trivialized bandings} as the homotopy fiber of the morphism
    \[
    \Gerbe_Y(G)\simeq \Map(Y,\bB\sAut_{\bH}(\bB G))\to \Map(Y,\bB\Out(G))
    \] over the base point of $\bB\Out(G)$.
\end{definition}
\begin{proposition}[see {\cite[Theorem 2.5]{Nocera-Pernice}}]\label{proposition-principal-BG-bundle-vs-G-gerbe-triv-band}
We have an equivalence of $\infty$-groupoids
\[
\Bund_Y(\fib(\tau_G))\simeq \Gerbe_Y^{\triv}(G)
\]
\end{proposition}
\begin{proof}
{The morphism $\tau_G$ is an effective epimorphism (which can be checked on the truncation), so we have a fiber sequence $\bB\fib(\tau_G)\to \bB\sAut_{\bH}(\bB G)\to \bB\Out(G)$.}
    The result then follows from the fiber sequence
    \[
    \Map(Y,\bB\fib(\tau_G))\to\Map(Y, \bB\sAut_{\bH}(\bB G))\to \Map(Y,\bB\Out(G)),
    \]
\end{proof}
\begin{remark}
For $G=\bG_m,\mu_n$, this shows that principal $\bB G$-bundles over $Y$ are equivalent to $G$-gerbes over $Y$ with trivialized bandings {(in the fppf topology for $\mu_n$).}
\end{remark}
There is an alternative approach to defining bandings in the literature (\cite{lieblich-moduli-twisted}, \cite{Binda-Porta}, \cite{bergh-BS}), but we  check in the appendix (Proposition \ref{proposition-compare-banding}) that this approach coincides with our previous definition for $G=\bG_m,\mu_n$.
\begin{definition}[{\cite[Definition 5.12]{Binda-Porta}}]
Let $f:X\to Y$ be a morphism of derived stacks, and let $H$ be an $\infty$-group in the $\infty$-topos of derived stacks over $Y$. Let $I_{X/Y}$ be the inertia stack with the $\infty$-group structure over $X$. An \textbf{$H$-banding on $X$} is an equivalence $H\simeq f_*I_{X/Y}$ of $\infty$-groups over $Y$.
\end{definition}
\begin{remark}\label{remark-BP-banding-Gm-gerbe}
When $f$ is a $\bG_m$-gerbe, the $\infty$-group $f_*I_{X/Y}$ is a twisted form of $\bG_m$ over $Y$ (Definition \ref{definition-twisted-form}), so we will interpret a $\bG_m$-banding on $X$ as a choice of an equivalence between $f_*I_{X/Y}$ and $(\bG_m)_Y$ in the $\infty$-groupoid of twisted forms of $\bG_m$ over $Y$.
\end{remark}
The following proposition appears known to the experts:
\begin{proposition}\label{proposition-compare-def-banded-gerbe}
There is an equivalence between the $\infty$-groupoid of $\bG_m$-gerbes over $Y$ with $(\bG_m)_Y$-bandings (in the sense of Remark \ref{remark-BP-banding-Gm-gerbe}) and $\Gerbe_Y^{\triv}(\bG_m)$.
\end{proposition}
\begin{proof}
    This follows from Proposition \ref{proposition-compare-banding}.
\end{proof}
\begin{remark}
    Similarly, there is an equivalence between the $\infty$-groupoid of $\mu_n$-gerbes over $Y$ with $(\mu_n)_Y$-bandings and $\Gerbe_Y^{\triv}(\mu_n)$ {in the fppf topology}.
\end{remark}
\subsection{Twisted derived category}
Let $f:X\to Y$ be a $\bG_m$-gerbe with trivialized banding over a scalloped (qcqs 1-Artin) derived stack $Y$.
$X$ is also scalloped (qcqs 1-Artin) by Proposition \ref{proposition-principal-BG-bundle-vs-G-gerbe-triv-band} and Proposition \ref{principal-BG-bundle-over-scalloped-is-scalloped}.
Also, by virtue of Proposition \ref{proposition-compare-def-banded-gerbe}, \cite[Definition 5.18]{Binda-Porta} applies to $f$ 
and for each $d\in\bZ$, we define the full subcategory
\[
\Dqc^d(X/Y)\subset \Dqc(X)
\] spanned by \textbf{$d$-homogeneous} objects. 
The idea is that in this setup we can consider the action map $\act:\bB \bG_m\times X\to X$ as well as the projections $p:\bB \bG_m\times X\to X,q:\bB \bG_m\times X\to \bB\bG_m$.  
There is a projection functor $\delta_d:\Dqc(X)\to\Dqc(X)$ to the $d$-th homogeneous component, which we define as $\delta_d(-):=p_*(\act^*(-)\otimes q^*\cO(-d))$. There is a canonical map $\delta_d\to \id$, and an object $F\in\Dqc(X)$ is $d$-homogeneous when this defines an equivalence  $\delta_d(F)\to F$. 
We will also call $d$-homogeneous objects \textbf{$d$-twisted (quasi-coherent) sheaves.}
When the base $Y$ is clear from context, we also write it as $\Dqc^d(X).$
For $d=1,$ we also refer to $1$-twisted sheaves as \textbf{$X$-twisted shaves}, or simply \textbf{twisted sheaves.} Also see \cite[Remark 2.9]{Nocera-Pernice} for details.
\begin{definition}
We call $\Dqc^1(X/Y)=\Dqc^1(X)$ the \textbf{twisted derived category of $X$.}  
\end{definition}
\begin{remark}
    \cite{Binda-Porta}
 appears to use $\bE_\infty$-ring spectra rather than animated rings, but in this paper, we will stick to derived stacks defined with animated rings. 
    Also, \cite{Binda-Porta} considers  possibly non-geometric derived stacks, but we don't need this generality in this paper. 
\end{remark}
 
By Proposition \ref{proposition-compare-def-banded-gerbe}, we can state \cite[Theorem 5.25 (1),(2)]{Binda-Porta} in the following form:
\begin{proposition}[{\cite[Theorem 5.25(1),(2)]{Binda-Porta},\cite{lieblich-moduli-twisted},\cite{bergh-BS}}]\label{proposition-decomp-Dqc-gerbe}
    Let $Y$ be as above (or more generally a geometric derived stack), and let $f:X\to Y$ be a $\bG_m$-gerbe with trivialized $(\bG_m)_Y$-banding. Then,
    \begin{enumerate}
        \item The pullback $f^*$ defines an equivalence $f^*:\Dqc(Y)\to\Dqc^0(X)$.
        \item The decomposition functor 
        \[
        D=\prod_{d\in\bZ}\delta_d:\Dqc(X)\to \prod_{d\in\bZ}\Dqc^d(X),
        \] is a $t$-exact equivalence. 
    \end{enumerate}
\end{proposition}
    \begin{definition}
Let $f:X\to Y$ be as above.
We define $\Dqc^1(X/Y)^\omega$ to be the full subcategory of $\Dqc^1(X/Y)$ spanned by the compact objects.
\end{definition}
\begin{proposition}\label{proposition-basic-property-Dqc1-presentable-stable-compact}
\leavevmode
Let $f:X\to Y$ be a principal $\bB\bG_m$-bundle of a scalloped derived stack.
    \begin{enumerate}
    \item  $\Dqc^1(X/Y)$ is a presentable stable $\infty$-category.
        \item $\Dqc^1(X/Y)^\omega=\Dqc^1(X/Y)\cap\Dqc(X)^\omega$.
        \item $\Dqc^1(X/Y)$ is tamely compactly generated, so that $\Ind(\Dqc^1(X/Y)^\omega)\simeq\Dqc^1(X/Y)$
    \end{enumerate}
\end{proposition}
\begin{proof}
    \leavevmode
    \begin{enumerate}
        \item We first observe in the product decomposition $\Dqc(X)\simeq\prod_i\Dqc^i(X/Y)$ that mapping spaces in the left hand side are the products of the mapping spaces of components on the right hand side. This immediately implies that $\Dqc^1(X/Y)$ is stable as one can compute limits and colimits of objects in $\Dqc^1(X/Y)$ in the bigger $\infty$-category $\Dqc(X/Y)$ and then project to $\Dqc^1(X/Y)$. 
        We also see that the inclusion $\iota_1:\Dqc^1(X/Y)\to \Dqc(X/Y)$ and the projection $\delta_1:\Dqc(X/Y)\to \Dqc^1(X/Y)$ are the left and the right adjoints of each other. In particular, the endofunctor $\iota_1\circ \delta_1:\Dqc(X/Y)\to\Dqc(X/Y)$ is a cocontinuous localization functor. The essential image is then presentable by \cite[Proposition 5.5.4.2]{HTT}. 
        \item This follows from the fact about mapping spaces in the product above, as well as the fact that $\iota_1,\delta_1$ are cocontinuous. 
        \item Since $\Dqc(X)$ is tamely compactly generated (Proposition \ref{proposition-tamely-compactly-generated}), we can find a small $\infty$-category  $\cC\subset \Dqc(X/Y)$ consisting of compact objects such that its right orthogonal in $\Dqc(X)$ is trivial. 
        We show that the right orthogonal of $\delta_1(\cC)\subset \Dqc^1(X/Y)$ is trivial, and that the objects of $\delta_1(\cC)$ are compact. 
        To see this, pick any $E\in \delta_1(\cC)^\perp\subset \Dqc^1(X/Y)$. For any $C\in\cC$, $\Map_{\Dqc(X)}(C,E)=\Map_{\Dqc^1(X/Y)}(\delta_1(C),E)=0$, so that by choice of $\cC$, $E=0.$ 
        For the latter claim about compactness, we observe again that the mapping space in the product $\infty$-category is computed by the product of mapping spaces of components, as well as that $\iota_1$ is cocontinuous. 

         By part (2), compact objects of $\Dqc^1(X/Y)$ are compact in $\Dqc(X)$, so that ``tameness'' follows from  Proposition \ref{proposition-tamely-compactly-generated}.
    \end{enumerate}
\end{proof}
\begin{remark}\label{remark-Perf1-binda-porta}
Since we assumed $Y$ is scalloped, $\Dqc(X)^\omega\simeq\Dperf(X)$, so by Proposition \ref{proposition-basic-property-Dqc1-presentable-stable-compact}, $\Dqc^1(X/Y)^\omega$ agrees with the notation $\Perf_1(X)$ from \cite[Notation 5.19]{Binda-Porta}.
\end{remark}

    \begin{corollary}\label{corollary-tw-Dqc-dualizable}
        $\Dqc^1(X/Y)$ {(or more generally $\Dqc^i(X/Y)$ for $i\in\bZ$)} is dualizable in $\PrSt$.
    \end{corollary} 
\begin{proof} We simply apply 
    \cite[Proposition D.7.2.3]{SAG} to the proposition.
\end{proof}
\begin{remark}\label{remark-relative-dualizability-transitivity}
{Since $Y$ is scalloped, \cite[Corollary D.7.7.3]{SAG} applies, and $\Dqc^1(X/Y)$ is also dualizable in $\Mod_{\Dqc(Y)}(\PrSt)$. 
} 
    In Proposition \ref{proposition-Dqc1-Dqc(-1)-dual} (or rather its generalization to all weights), in fact we will  see that if $f$ is as before, $\Dqc^i(X/Y)$ is dualizable 
    in $\Mod_{\Dqc(Y)}(\PrSt)$ with
    the dual given by $\Dqc^{-i}(X/Y)$.
\end{remark}

\begin{lemma}[{\cite[Theorem 1.2]{bzfn}}]\label{lemma-kunneth-perfect}
Let $X,Y,X'$ be perfect derived stacks. Then, the obvious functor
    \[
\Dqc(X)\otimes_{\Dqc(Y)}\Dqc(X')\to\Dqc(X\times_YX').
    \] is an equivalence.
\end{lemma}

Next, we study how the symmetric monoidal structure interacts with the decomposition into the homogeneous components. We start with the simplest case.
\begin{example}\label{example-split-BGm-bundle-day-convolution}
Let $Y$ be a scalloped derived stack, and let    $X=Y\times_{\bZ}(\bB\bG_m)_{\bZ}$. 
We know that $\Dqc((\bB\bG_m)_\bZ)\simeq\Fun(\bZ,\Dqc(\bZ))$, where $\bZ$ is a discrete category and the symmetric monoidal structure on the right is given by the Day convolution (see \cite[Theorem 5.2]{Binda-Porta} for example). 
If we ignore the symmetric monoidal structure, we can write the right hand side as 
\[
\Fun(\coprod_{d\in\bZ}\{*\},\Dqc(\bZ))\simeq\prod_{d\in\bZ}\Fun(\{*\},\Dqc(\bZ))\simeq \prod_{d\in\bZ}\Dqc(\bZ).
\] 
The object assigning $\bZ$ to the index $d$ and $0$ for the other indices is denoted by $\cO_{\bZ}(d),$ following \cite[Notation 5.5]{Binda-Porta}. 
Since relevant derived stacks are perfect, we have a symmetric monoidal equivalence by Lemma \ref{lemma-kunneth-perfect}:
\[
\boxtimes:\Dqc(Y)\otimes_{\Dqc(\bZ)}\Dqc((\bB\bG_m)_{\bZ})\simeq \Dqc(Y\times_{\bZ}(\bB\bG_m)_{\bZ})
\]
Since $\Dqc(Y)$ is dualizable, 
the left hand side can be written as
\[
\Dqc(Y)\otimes_{\Dqc(\bZ)}\prod_{d\in\bZ}\Dqc(\bZ)\simeq \prod_{d\in \bZ}\left(\Dqc(Y)\otimes_{\Dqc(\bZ)}\Dqc(\bZ)\right)\simeq \prod_{d\in \bZ}\Dqc(Y).
\]
Combining these, we see that any object $F\in\Dqc(Y\times_{\bZ}(\bB\bG_m)_\bZ)$ can be written as $F=\oplus_{d\in\bZ}F_d\boxtimes\cO_{\bZ}(d)$ for some $F_d\in\Dqc(Y).$ 
Moreover, given elements of form $F_i\boxtimes \cO_{\bZ}(i), F_j\boxtimes \cO_{\bZ}(j)$, their tensor product is 
$(F_i\otimes F_j)\boxtimes \cO_{\bZ}(i+j)
$ since $\cO_\bZ(i)\otimes\cO_\bZ(j)\simeq\cO_\bZ(i+j)$ in view of the Day convolution. These fully describe the symmetric monoidal structure in the trivial principal $\bB\bG_m$-bundle case.

    Also note by projection formula (which holds since relevant derived stacks are perfect as in \cite{bzfn}), the functor $p_*:\Dqc(Y\times_{\bZ}(\bB\bG_m)_{\bZ})\to \Dqc(Y)$ is $\Dqc(Y)$-linear and is describable using the functor $p_*:\Dqc((\bB\bG_m)_\bZ)\to\Dqc(\bZ)$. This implies for instance that $p_*(F_d\boxtimes\cO_{\bZ}(d))=0$ when $d\neq0$, and that $p_*p^*F_0\simeq F_0$ for $F_d,F_0\in\Dqc(Y)$.
\end{example}

\begin{proposition}\label{proposition-symmetric-monoidal-Dqc}
    Let $f:X\to Y$ be a $\bG_m$-gerbe with a trivialized banding over a scalloped derived  stack $Y$. Then, for each $i,j\in\bZ$, the symmetric monoidal structure  $(-)\otimes(-):\Dqc(X)\times \Dqc(X)\to \Dqc(X)$ restricts to a functor 
    \[
    (-)\otimes(-):\Dqc^i(X)\times \Dqc^j(X)\to \Dqc^{i+j}(X).
    \]
\end{proposition}
\begin{proof}
{
When $F\in\Dqc(X)$ is $d$-homogenous, the canonical map $\delta_d(F)\to F$ is an equivalence with $\delta_d(-):=p_*(\act^*(-)\otimes q^*\cO_{\bZ}(-d))$, $p:\bB \bG_m\times X\to X$ and $ q:\bB \bG_m\times X\to \bB\bG_m$. 
As in Proposition \ref{proposition-decomp-Dqc-gerbe}, we have $\oplus_w\delta_w(F)\simeq F$, and we could also phrase $d$-homogeneity as the vanishing of $\delta_w(F)$ for all $w\neq d.$ By  Example \ref{example-split-BGm-bundle-day-convolution}, this means precisely that 
$\act^*F$ can be written as $F_d\boxtimes \cO_{\bZ}(d)$ for some $F_d\in\Dqc(X)$ (where necessarily $F\simeq F_d$ because $p_*p^*F_d\simeq F$).

Now, for objects $F\in\Dqc^i(X),G\in\Dqc^j(X)$, 
we can write $\act^*F\simeq F\boxtimes \cO_{\bZ}(i),\act^*G\simeq G\boxtimes \cO_{\bZ}(j)$. Then,
$\act^*(F\otimes G)\simeq (F\otimes G)\boxtimes\cO_{\bZ}(i+j)$ by Example \ref{example-split-BGm-bundle-day-convolution} again. Then, we have $p_*(\act^*(F\otimes G)\boxtimes\cO_{\bZ}(-i-j))\simeq F\otimes G$, so $F\otimes G\in\Dqc^{i+j}(X).$ 
}
  
\end{proof}
\begin{remark}\label{remark-twisted-Dqc-tensor}
    Similarly, for $\bG_m$-gerbes with trivialized banding $X,X'$ over $Y$, the symmetric monoidal structure on $\Dqc(X\times_YX')$ defines functors
    \[
    \Dqc^{i,i'}(X\times_YX')\times\Dqc^{j,j'}(X\times_YX')\to \Dqc^{i+j,i'+j'}(X\times_YX')
    \] for integers $i,j,i',j'.$
\end{remark}

\begin{lemma}\label{lemma-base-change-tw-Dqc}
Assume $Y'\to Y$ is a morphism of scalloped derived stacks.
Let $X\to Y$ be a $\bG_m$-gerbe with trivialized banding. Then, we have a canonical equivalence
\[
\Dqc^i(X/Y)\otimes_{\Dqc(Y)}\Dqc(Y')\simeq\Dqc^i(X'/Y'),
\]
where we set $X'=X\times_YY'.$
\end{lemma}
\begin{proof}
By Proposition \ref{proposition-khan-ravi-scalloped-perfect}, scalloped derived stacks are perfect, so Lemma \ref{lemma-kunneth-perfect} applies, and the box tensor product defines an equivalence
\[
\phi:\Dqc(X)\otimes_{\Dqc(Y)}\Dqc(Y')\simeq\Dqc(X').
\]
If $g:X'\to X$ denotes the projection, we have $g^*\circ\delta_{d,X}\simeq \delta_{d,X'}\circ g^*$ where $\delta_{d,X}$ is the $d$-th projection operator  on $X$ (and $\delta_{d,X'}$ is similar). This shows the pullback $\Dqc(X)\to\Dqc(X')$ sends $\Dqc^d(X)$ to $\Dqc^d(X').$
Moreover, by Proposition \ref{proposition-decomp-Dqc-gerbe}, the pullback $\Dqc(Y')\to\Dqc(X')$ defines an equivalence $\Dqc(Y')\simeq\Dqc^0(X')$, so by Proposition \ref{proposition-symmetric-monoidal-Dqc}, we see that 
$\phi$ restricts to a functor
\[
\phi_d:\Dqc^d(X)\otimes_{\Dqc(Y)}\Dqc(Y')\to\Dqc^d(X')
\] for each $d\in\bZ$. 
We can put them together to obtain a commutative diagram:
\[
\xymatrix{
\prod_d\left(\Dqc^d(X)\otimes_{\Dqc(Y)}\Dqc(Y')\right)\ar[r]^-{\prod_d\phi_d}&\prod_d\Dqc^d(X')\\
\prod_d\left(\Dqc^d(X)\right)\otimes_{\Dqc(Y)}\Dqc(Y')\ar[u]^-\simeq&\Dqc(X')\ar[u]^\simeq\\
\Dqc(X)\otimes_{\Dqc(Y)}\Dqc(Y')\ar[ur]_-{\phi}^-\simeq\ar[u]^\simeq
}
\]
Here, the top left vertical equivalence comes from the compatibility of product and tensor product (which relies on dualizability of $\Dqc(Y')$). The right vertical arrow and the bottom left vertical arrow are equivalences by Proposition \ref{proposition-decomp-Dqc-gerbe}.
In particular, $\prod_d\phi_d$ is an equivalence, and $\phi_d$ must be an equivalence as well.
\end{proof}
\begin{remark}\label{remark-relative-base-change-tw-Dqc}
    In the same fashion, one has the following relative version of Lemma \ref{lemma-base-change-tw-Dqc}.
    Let $Y\to S\leftarrow S'$ be a diagram of scalloped derived stacks. If $X\to Y$ is a $\bG_m$-gerbe with trivialized banding, for each $i\in\bZ$, we have 
    \[
\Dqc^i(X/Y)\otimes_{\Dqc(S)}\Dqc(S')\simeq\Dqc^i(X\times_SS'/Y\times_SS').
    \]
\end{remark}
\begin{proposition}[see {\cite[Proposition 3.2]{Nocera-Pernice}}]\label{proposition-twisted-kunneth}
    Let $Y$ be a scalloped derived stack, and let $X,X'$ be $\bG_m$-gerbes with trivialized bandings over $Y$. Then, for each integers $i,j\in\bZ$, we have a canonical equivalence defined by the box tensor product:
    \[
\Dqc^i(X)\otimes_{\Dqc(Y)}\Dqc^j(X')\to\Dqc^{i,j}(X\times_YX')
    \]
\end{proposition}
\begin{proof}
We first observe that by the same proof as \cite[Proposition 3.2]{Nocera-Pernice}, the box tensor product functor induces a morphism 
\[
\Dqc^i(X)\otimes_{\Dqc(Y)}\Dqc^j(X')\to\Dqc^{i,j}(X\times_YX')
    \]
    for each $i,j\in\bZ.$

Using Corollary \ref{corollary-BG-bundle-locally-split-in-topology}, we can find a smooth atlas of $Y$ so that the pullback of $X\to Y$ along the cover can be trivialized. 
Since $Y$ is quasi-compact, we can refine the 
atlas if necessary to assume the atlas is a derived affine. The Cech nerve of the atlas consists of derived affines since $Y$ has an affine diagonal. We write 
     $Y\simeq \colim_{m\in I} Y_m$ 
   using this Cech nerve.   
     We write $X_m:=X\times_YY_m,X'_m:=X'\times_YY_m$.
     All relevant derived stacks are scalloped by Proposition \ref{proposition-scalloped-dstack-stable-fiber-product}.
    By \cite[Lemma 5.20]{Binda-Porta}, we have 
    \[
    \Dqc^{i,j}(X\times_YX')\simeq\lim_{m\in I^{\op}}\Dqc^{i,j}(X_m\times_{Y_m}X_m').
    \]
    $\Dqc^j(X')$ is dualizable over $\Dqc(Y)$ by Remark \ref{remark-relative-dualizability-transitivity}. 
    Hence, by \cite[Lemma 5.20]{Binda-Porta}, we have
    \begin{equation*}
\begin{split}
 \Dqc^i(X)\otimes_{\Dqc(Y)}\Dqc^j(X')
 &\simeq
    \left(\lim_{m\in I^{\op}}  \Dqc^i(X_m/Y_m)\right)\otimes_{\Dqc(Y)}\Dqc^j(X')\\
    &\simeq \lim_{m\in I^{\op}}  (\Dqc^i(X_m/Y_m)\otimes_{\Dqc(Y)}\Dqc^j(X')).
\end{split}
\end{equation*}
    Since $Y_m$'s are in particular scalloped, by
    Lemma \ref{lemma-base-change-tw-Dqc}, we can rewrite this as 
    \[
    \lim_{m\in I^{\op}}  (\Dqc^i(X_m/Y_m)\otimes_{\Dqc(Y_m)}\Dqc^j(X'_m/Y_m)).
    \]
    Therefore, in this situation, the functor
    \[
\Dqc^i(X)\otimes_{\Dqc(Y)}\Dqc^j(X')\to\Dqc^{i,j}(X\times_YX')
    \]
    arises as the limit of functors
    \[
\Dqc^i(X_m)\otimes_{\Dqc(Y_m)}\Dqc^j(X_m')\to\Dqc^{i,j}(X_m\times_{Y_m}X_m')
    \] over $m\in I^{\op}$.
    Hence, it suffices to treat the case of trivial principal $\bB\bG_m$-bundles. 
    This case is proved as in \cite[Proposition 3.2]{Nocera-Pernice}.
\end{proof}
\begin{remark}\label{remark-relative-twisted-kunneth}
    Combining Remark \ref{remark-relative-base-change-tw-Dqc} and Proposition \ref{proposition-twisted-kunneth}, one has the following relative version of Proposition \ref{proposition-twisted-kunneth}. If $Y\to S\leftarrow Y'$ is a diagram of scalloped derived stacks, and if $X\to Y,X'\to Y'$ are $\bG_m$-gerbes with trivialized bandings, for each $i,j\in\bZ$, the box tensor product defines an equivalence: 
    \[
    \Dqc^{i}(X/Y)\otimes_{\Dqc(S)}\Dqc^j(X'/Y')\simeq 
    \Dqc^{i,j}(X\times_SX'/Y\times_SY').
    \]
\end{remark}
\begin{proposition}\label{proposition-BGmbundle-pushforward}
    Let $f:X\to Y$ be a $\bG_m$-gerbe with trivialized banding over a scalloped derived stack. Then, 
    \begin{enumerate}
        \item  the canonical morphism $\eta:\id_{\Dqc(Y)}\to f_* f^*$ is an equivalence,
        \item the restricted functor $f_*:\Dqc^0(X)\to\Dqc(Y)$ defines a $t$-exact quasi-inverse of $ f^*:\Dqc(Y)\to\Dqc^0(X)$ from Proposition \ref{proposition-decomp-Dqc-gerbe}(1),
        \item {there is a canonical equivalence $\delta_0\simeq f^*f_*$,} and
        \item the functor $ f_*:\Dqc(X)\to\Dqc(Y)$ is {canonically equivalent} to the composition of $\delta_0:\Dqc(X)\to\Dqc^0(X)$ with the equivalence $f_*:\Dqc^0(X)\simeq \Dqc(Y)$ from (2). In particular, $f_*$ is $t$-exact.
    \end{enumerate}
\end{proposition}
\begin{proof}
\leavevmode
\begin{enumerate}
    \item  $f$ is perfect by Lemma \ref{lemma-nice-group-derived-pushforard-structure-sheaf}, so \cite[Proposition 3.10]{bzfn} implies that projection formula holds, and we can identify the morphism with $(-)\otimes\cO_Y\to(-)\otimes  f_*\cO_X $. The same lemma proves that $\cO_Y\to f_*\cO_X$ is an equivalence, so $\eta$ is an equivalence. 
    \item $f$ is flat, so $f^*$ is $t$-exact. Then, the result is immediate from part (1).
    \item $f$ is perfect, so we can apply base change to the following cartesian diagram by \cite[Proposition 3.10]{bzfn}:
    \[
    \xymatrix{
    X\times\bB\bG_m\ar[r]^-{\act}\ar[d]^-{p}&X\ar[d]^-f\\
    X\ar[r]^-f&Y
    }
    \]
    This shows $f^*f_*\simeq p_*\act^*$, but $\delta_0\simeq p_*\act^*$ by definition, so we obtain the canonical equivalence.
    \item Unit and counit of adjunction give us the following commutative diagram. \[
    \xymatrix{
    f_*\ar[r]^-{\eta\circ f_*}\ar[dr]_-{\id}&f_*f^*f_*\ar[d]^-{f_*\circ\epsilon}\\
    &f_*
    }
    \]
    By part (1), all morphisms in the diagram are equivalences, so by part (3), we have a canonical equivalence
    \[
    f_*\simeq f_*f^*f_*\simeq f_*\delta_0.
    \]
    $\delta_0:\Dqc(X)\to\Dqc^0(X)$ is $t$-exact by \cite[Lemma 5.20]{Binda-Porta}, and $f_*:\Dqc^0(X)\to\Dqc(Y)$ is $t$-exact by part (2), so the proof is complete.
\end{enumerate}
\end{proof}

\setcounter{appendixFix}{\value{section}}
\begin{appendix}
\section{Twisted Fourier-Mukai transform}
The Fourier-Mukai theory for twisted sheaves is well-known in certain settings (see \cite{caldararu} for example), but in this section, we develop the theory in the generality needed for our purposes. 
Let $X\to Y$ be a $\bG_m$-gerbe with trivialized banding over a scalloped derived stack $Y$. Then, $\Delta:X\to X\times_YX$ defines a section of the two projections. 
Let us label the two copies of $X$ appearing in $X\times_YX$ and write $X_1\times_YX_2.$
Since the datum of a $\bG_m$-gerbe with trivialized banding is equivalent to that of a principal $\bB\bG_m$-bundle, the section lets us trivialize the $\bG_m$-gerbes (with trivial banding) $p_i:X_1\times_YX_2\to X_i$ and we can identify $p_i$ with $X_{i}\times\bB\bG_m\to X_i$ for $i=1,2$.
Write $\phi_i:X_i\times\bB\bG_m\to X_1\times_YX_2$ for the trivialization.
\begin{lemma}\label{lemma-trivialize-twist-band}
    In the situation above, the equivalence $\phi_1$ induces an equivalence of banding
    \[
   \widetilde{\phi_1}: (\bG_m)_Y^{(X_1)}\times_Y(\bG_m)_Y^{(\bB\bG_m)}\to  (\bG_m)_Y^{(X_1)}\times_Y(\bG_m)_Y^{(X_2)};(a,b)\mapsto (a,ab),
    \]
    where for each copy of $\bB\bG_m$, the superscript denotes the gerbe whose banding gives rise to that copy for clarity. 
    Similarly, the equivalence $\phi_2$ induces an equivalence of banding
    \[
   \widetilde\phi_2: (\bG_m)_Y^{(X_2)}\times_Y(\bG_m)_Y^{(\bB\bG_m)}\to  (\bG_m)_Y^{(X_1)}\times_Y(\bG_m)_Y^{(X_2)};(a,b)\mapsto (ab,a).
    \]
    In particular, for each $i,j\in\bZ$, $\phi_i$ induce equivalences
    \[
    \phi_1^*:\Dqc^{i,j}(X_1\times_YX_2) \simeq \Dqc^{i+j,j}(X_1\times\bB\bG_m)
    \]
    \[
    \phi_2^*:\Dqc^{i,j}(X_1\times_YX_2) \simeq \Dqc^{i+j,i}(X_2\times\bB\bG_m)
    \]
\end{lemma}
\begin{proof} 
The formulas for the $\widetilde{\phi_i}$'s come from the fact that the maps $\phi_i$'s are constructed using the splittings and the $\bB\bG_m$-actions. 
For the computation of $\phi_1^*,$ for any $E\in\Dqc^{i,j}(X_1\times_YX_2),$ the ``weight'' of $\phi_1^*E$ is given by the character
\[
(\bG_m)_Y\times_Y(\bG_m)_Y\xrightarrow[]{\widetilde\phi_1} (\bG_m)_Y\times_Y(\bG_m)_Y
\xrightarrow[]{(a,b)\mapsto a^ib^j}(\bG_m)_Y;
\]
\[
(a,b)\mapsto(a,ab)\mapsto a^i(ab)^j=a^{i+j}b^j,
\]
so $\phi_1^*E\in\Dqc^{i+j,j}(X_1\times\bB\bG_m).$ The case of $\phi^*_2$ is similar.  
\end{proof}
Using Remark \ref{remark-twisted-Dqc-tensor} and Lemma \ref{lemma-trivialize-twist-band}, we can define a line bundle \[\cP_X:=(\phi_1)_*\pi_1^*\cO(1)\in\Dqc^{-1,1}(X_1\times_YX_2),\]
where 
$\pi_i:X_i\times\bB\bG_m\to\bB\bG_m$ is the projection for $i=1,2.$
By the following lemma, $\cP_X$ can also be defined using $\phi_2$ in the symmetric manner (up to a sign):
\begin{lemma}\label{lemma-symmetry-twisted-diagonal}
    In the situation above, we have $\cP_X=(\phi_1)_*\pi_1^*\cO(1)\simeq (\phi_2)_*\pi_2^*\cO(-1)$.
\end{lemma}
\begin{proof}
By construction of trivializations, the composition $\phi_2^{-1}\circ\phi_1:X_1\times\bB\bG_m\simeq X_2\times\bB\bG_m$ is equivalent to $(\mu,\iota\circ p_2)$ where $\mu:X_1\times\bB\bG_m\to X_2$ is the action map, and $\iota:\bB\bG_m\to\bB\bG_m$ is the inversion map. If we factor this as the composition $ (\mu,\iota\circ p_2)=(\id_X\times\iota)\circ (\mu ,p_2)$, we obtain
\[
(\mu,\iota\circ p_2)^*\pi_2^*\cO(-1)\simeq  (\mu ,p_2)^*(\id_X\times\iota)^*\pi_2^*\cO(-1)\simeq (\mu ,p_2)^*\pi_2^*\cO(1)\simeq \pi_1^*\cO(1),
\] where the last equality holds because $(\mu ,p_2)$ is $\bB\bG_m$-linear.
\end{proof}
\begin{definition}\label{definition-twisted-Fourier-Mukai-transform}
Let $Y\to S$, {$Y'\to S$} be morphisms of scalloped derived stacks. 
    Let $X\to Y,X'\to Y'$ be two $\bG_m$-gerbes with trivialized bandings. Given an integer $n\in\bZ$ and an object \[K\in \Dqc^{-n,n}(X\times_SX'/Y\times_SY'),\] we define \textbf{the ($S$-linear) $n$-twisted Fourier-Mukai transform associated with the kernel $K$} as the functor 
    \[
    \Phi_K^{(n)}=(q_2)_*(q_1^*(-)\otimes K):\Dqc^n(X)\to\Dqc^n(X'),
    \] where $q_1:X\times_SX'\to X,q_2:X\times_SX'\to X'$ are the projections. When $n=1$, we also call it \textbf{the ($S$-linear) twisted Fourier-Mukai transform  associated with the kernel $K$.}
\end{definition}
\begin{proposition}\label{proposition-twisted-id-FM-kernel-1}
   The $Y$-linear functor $\Phi_{\cP_X}^{(1)}:\Dqc^1(X_1)\to\Dqc^1(X_2)$ defined by the kernel $\cP_X\in\Dqc^{-1,1}(X_1\times_YX_2)$ is canonically equivalent to the identity.
\end{proposition}
\begin{proof}
Up to the identification $\phi_2$, the projection map $q_1:X_1\times_YX_2\to X_1$ is identified with the action map $\act:X_2\times \bB\bG_m\to X_1$, 
and the other projection $q_2$ is identified with the projection $\pi_2:X_2\times \bB\bG_m\to X_2$. As a result, the twisted Fourier-Mukai transform $\Phi_{\cP_X}^{(1)}$ is equivalent to the functor
\[
(\pi_2)_*(\act^*(-)\otimes \pi_2^*\cO(-1))
\] by Lemma \ref{lemma-symmetry-twisted-diagonal}.
As in \cite[Remark 2.9]{Nocera-Pernice}, we have a canonical morphism $(\pi_2)_*(\act^*(-)\otimes \pi_2^*\cO(-1))\to\id$ between functors, and by definition of $1$-twisted sheaves, this morphism is an equivalence (when we restrict to $\Dqc^1(X_1)\subset\Dqc(X)$). 
\end{proof}
If we switch the roles of $X_1$ and $X_2$, we can use the kernel $\cP_X$ to define a $(-1)$-twisted Fourier-Mukai transform
\[
\Phi^{(-1)}_{\cP_X}:\Dqc^{-1}(X_2)\to\Dqc^{-1}(X_1).
\] 
In parallel to the previous proposition, we obtain:
\begin{proposition}\label{proposition-twisted-id-FM-kernel-2}
   The $Y$-linear functor $\Phi_{\cP_X}^{(-1)}:\Dqc^{-1}(X_2)\to\Dqc^{-1}(X_1)$ is canonically equivalent to the identity.
\end{proposition}
\begin{proof}
    This is proved as in Proposition \ref{proposition-twisted-id-FM-kernel-1} by symmetry. In this proposition, we trivialize $X_1\times_YX_2$ via $\phi_1$ and then use the description of $\cP_X$ as $(\phi_1)_*\pi_1^*\cO(1)$.
\end{proof}

\begin{proposition}\label{proposition-Dqc1-Dqc(-1)-dual}
    Let $f:X\to Y$ be a $\bG_m$-gerbe with trivialized banding over a scalloped derived stack $Y.$ Then, $\Dqc^1(X)$ is dualizable over $\Dqc(Y)$, and the dual is given by $\Dqc^{-1}(X)$. 
\end{proposition}
\begin{proof}
    We need to define $Y$-linear functors
    \[
    \ev:\Dqc^{-1}(X)\otimes_{\Dqc(Y)}\Dqc^1(X)\to\Dqc(Y),
    \quad
    \coev:\Dqc(Y)\to \Dqc^{1}(X)\otimes_{\Dqc(Y)}\Dqc^{-1}(X)
    \] such that the two compositions
    \[
    F:\Dqc^1(X)\xrightarrow[]{\coev\otimes\id}\Dqc^1(X)\otimes_{\Dqc(Y)}\Dqc^{-1}(X)\otimes_{\Dqc(Y)}\Dqc^1(X)\xrightarrow[]{\id\otimes\ev}\Dqc^1(X),
    \]
    \[
    G:\Dqc^{-1}(X)\xrightarrow[]{\id\otimes\coev}\Dqc^{-1}(X)\otimes_{\Dqc(Y)}\Dqc^{1}(X)\otimes_{\Dqc(Y)}\Dqc^{-1}(X)\xrightarrow[]{\ev\otimes\id}\Dqc^{-1}(X)
    \]
    are equivalent to the identity functors (\cite[Definition 4.6.1.1]{HA}).
    We set
    \[
    \ev:=f_*\circ\Delta_{X/Y}^*:\Dqc^{-1,1}(X\times_YX)\to\Dqc(Y),
    \] up to the equivalence of Proposition \ref{proposition-twisted-kunneth}, where $f:X\to Y$ is the structure morphism of the gerbe.
    Up to Proposition \ref{proposition-twisted-kunneth}, we define 
    \[\coev=\pi^*(-)\otimes\cP_X:\Dqc(Y)\to \Dqc^{1,-1}(X_2\times_YX_1),\] where $\pi:X_2\times_YX_1\to Y$ is the projection, and $\cP_X$ is defined with respect to the switched ordering $X_2\times_YX_1$ instead of $X_1\times_YX_2$ {(so that it is $(1,-1)$-twisted)}.
    All relevant derived stacks are perfect  (Proposition \ref{proposition-khan-ravi-scalloped-perfect}), so that 
    base change holds by \cite[Proposition 3.10, Corollary 3.23]{bzfn}. Then, by
    standard base change argument,
    $F\simeq \Phi_{\cP_X}^{(1)}:\Dqc^1(X)\to\Dqc^1(X)$, and
    $G\simeq\Phi_{\cP_X}^{(-1)}:\Dqc^{-1}(X)\to \Dqc^{-1}(X).$ These two functors are equivalent to the identity by Propositions \ref{proposition-twisted-id-FM-kernel-1} and  \ref{proposition-twisted-id-FM-kernel-2}.
\end{proof}
\begin{remark}\label{remark-relative-Dqc1-Dqc(-1)-dual}
    In the previous proposition, if we assume $Y\to S$ is a 
    morphism of scalloped derived stacks, one can modify the proof to show that $\Dqc^1(X)$ is dualizable over $\Dqc(S)$ (instead of $\Dqc(Y)$) with the dual given by $\Dqc^{-1}(X)$ (see also Remark \ref{remark-relative-dualizability-transitivity}).

To see this, we let $\iota^X:X\times_YX\to X\times_SX$ be the obvious morphism. If $f^{S}:X\to S, \pi^S:X_2\times_SX_1\to S$ are the obvious morphisms, we can set
\[
\ev:=f^S_*\circ\Delta_{X/S}^*:\Dqc^{-1,1}(X\times_SX/Y\times_SY)\to\Dqc(S),
\]
\[
\coev:=(\pi^S)^*(-)\otimes\iota^X_*\cP_X:\Dqc(S)\to\Dqc^{1,-1}(X_2\times_SX_1/Y\times_SY).
\]
Then, we construct $\Dqc(S)$-linear functors $F:\Dqc^1(X)\to\Dqc^1(X), G:\Dqc^{-1}(X)\to\Dqc
^{-1}(X)$ as in the proof of Proposition \ref{proposition-Dqc1-Dqc(-1)-dual}. By the same argument, we see $F\simeq \Phi_{\iota^X_*\cP_X}^{(1)}, G\simeq \Phi^{(-1)}_{\iota^X_*\cP_X}$, where the Fourier-Mukai transforms are defined $S$-linearly rather than $Y$-linearly. By projection formula, we see $\Phi_{\iota^X_*\cP_X}^{(1)}\simeq \Phi_{\cP_X}^{(1)}, \Phi_{\iota^X_*\cP_X}^{(-1)}\simeq \Phi_{\cP_X}^{(-1)},$ where the right hand sides of the equivalences are $Y$-linear Fourier-Mukai transforms. Thus, these functors are equivalent to the identity by Propositions \ref{proposition-twisted-id-FM-kernel-1} and \ref{proposition-twisted-id-FM-kernel-2}.
\end{remark}
\begin{lemma}\label{lemma-dualizability-functor-Fourier-Mukai}
Let $S$ be a derived stack. Let $C$ be a dualizable object of $\cA:=\Mod_{\Dqc(S)}(\PrSt)$ with the dual  given by $C^\vee$ (with $\ev,\coev$ being the evaluation and coevaluation). Let $D$ be another object in $\cA$. Then, we have an equivalence  
\[
\Phi: C^\vee\otimes D\simeq {\Fun}^{\rL}_{S}(C,D),
\]
where $\Phi$ sends an object $K\in C^\vee\otimes D$ to the functor
\[
C\xrightarrow[]{\id\otimes K}C\otimes C^\vee\otimes D\simeq C^\vee\otimes C\otimes D\xrightarrow[]{\ev\otimes \id}D,
\] and $\Phi^{-1}$ sends a functor $F:C\to D$ to the functor
\[
\mathbf{1}\xrightarrow[]{\coev} C\otimes C^\vee\xrightarrow[]{F\otimes\id} D\otimes C^\vee\simeq C^\vee\otimes D,
\] which is identified with an object of $C^\vee\otimes D.$
\end{lemma}
\begin{proof}
    This immediately follows from the compatibility conditions between $\ev$ and $\coev.$ 
\end{proof}
\begin{example}\label{example-twisted-FMT-via-dualizability}
In the situation of Definition \ref{definition-twisted-Fourier-Mukai-transform}, 
we consider the symmetric monoidal $\infty$-category $\cA=\Mod_{\Dqc(S)}(\PrSt)$. Set $C=\Dqc^1(X/Y),D=\Dqc^1(X'/Y')$. By Remark \ref{remark-relative-Dqc1-Dqc(-1)-dual}, we have $C^\vee=\Dqc^{-1}(X/Y),$ 
so using Remark \ref{remark-relative-twisted-kunneth}, we can write the functor from
Lemma \ref{lemma-dualizability-functor-Fourier-Mukai} as
\[
\Phi:\Dqc^{-1,1}(X\times_S X'/Y\times_SY')\simeq\Fun_{\Mod_{\Dqc(S)}}^\rL(\Dqc^1(X/Y),\Dqc^1(X'/Y')).
\]
By the formula of $\Phi$ in the lemma, we see that $\Phi(K)\simeq \Phi_K^{(1)}$ for any kernel $K\in\Dqc^{-1,1}(X\times_S X'/Y\times_SY').$
\end{example}
Using Example \ref{example-twisted-FMT-via-dualizability}, we will next compute some examples of twisted Fourier-Mukai kernels giving rise to some well-known functors between twisted derived categories.
\begin{example}\label{example-twisted-FMT-tensor}
    Let $Y\to S$ be a morphism of scalloped derived stacks, and consider a $\bG_m$-gerbe $X\to Y$ with trivialized banding. For any $E\in\Dqc^0(X/Y)\simeq \Dqc(Y)$, we have a functor $E\otimes(-):\Dqc^1(X)\to\Dqc^1(X)$ by Proposition \ref{proposition-symmetric-monoidal-Dqc}.
    This functor arises as $\Phi(K)$ with
    {$K=\iota^X_*\cP_X\otimes p_2^*E\simeq\iota^X_*(\cP_X\otimes(\iota^X)^* p_2^*E)\simeq \iota_*^X(\cP_X\otimes p_2^*E).$}
    Indeed, by the formula of Lemma \ref{lemma-dualizability-functor-Fourier-Mukai}, the functor $\Phi(K)$ sends $F\in\Dqc^1(X/Y)$ to $(q_2)_*(\Delta_{X/S}\times_S\id_X)^*(p_1^*F\otimes p_{23}^*K)$, where $q_i:X\times_SX\to X$, $p_i:X\times_S X\times_SX\to X,p_{ij}:X\times_SX\times_SX\to X\times_SX$ are the projections to the components specified by the indices $i,j$.  
    $ p_1\circ (\Delta_{X/S}\times_S\id_X)=q_1,  p_{23}\circ (\Delta_{X/S}\times_S\id_X)=\id_{X\times_SX}.$
    Hence, the above expression reduces to 
    \[
    (q_2)_*(q_1^* F\otimes K)\simeq  E\otimes (q_2^Y)_*((q_1^Y)^* F\otimes\cP_{X}),
    \] where the equivalence uses projection formula, and $q^Y_i:X\times_YX\to X$ is the projection specified by the index $i$. This is equivalent to $E\otimes F$ by the argument of Proposition \ref{proposition-twisted-id-FM-kernel-1}.
\end{example}
\begin{example}\label{example-twisted-FMT-push-pull}
Let $S$ be a scalloped derived stack, and consider a morphism $Y\to Y'$ of scalloped derived stacks over $S$. Assume we have a cartesian diagram of derived stacks where the vertical maps are $\bG_m$-gerbes with trivialized banding (which are compatible with pullback):
    \[
\xymatrix{
X\ar[r]^-f\ar[d]& X'\ar[d]\\
Y\ar[r]^-g&Y'
}    
    \]
    In this situation, being a morphism between perfect derived stacks, $g$ is a perfect morphism (\cite[Corollary 3.23]{bzfn}). This allows us to use projection formula and base change as in \cite[Proposition 3.10]{bzfn}. Then, an easy computation shows that functors $f_*:\Dqc(X)\to\Dqc(X')$ and  $f^*:\Dqc(X')\to\Dqc(X)$ respect $n$-twisted sheaves. In particular, we can talk about the restrictions 
    \[
    f_*:\Dqc^1(X)\to\Dqc^1(X'),\quad  f^*:\Dqc^1(X')\to\Dqc^1(X).
    \]
    \begin{itemize}
        \item{$f_*:$} We have $f_*=\Phi^{(1)}_{(\id\times f)_*\iota^X_*\cP_X}\in\Dqc^{-1,1}(X\times_SX'),$ where we recall the notation $\iota^X:X\times_YX\to X\times_SX$. This follows immediately from Remark \ref{remark-relative-Dqc1-Dqc(-1)-dual} and the formula in Lemma \ref{lemma-dualizability-functor-Fourier-Mukai}.
        \item{$f^*:$}
        Again, by Remark \ref{remark-relative-Dqc1-Dqc(-1)-dual} and  Lemma \ref{lemma-dualizability-functor-Fourier-Mukai}, the kernel of interest is given by 
        $K:=
        (\id\times f)^*\iota_*^{X'}\cP_{X'}\in\Dqc^{-1,1}(X'\times_SX)$.
        We can rewrite this kernel to highlight its similarity to $f_*$.
        In the following pullback square, up to the identification of the top left corner $X\times_YX\simeq X'\times_{Y'} X$, the top horizontal arrow and the left vertical arrow can be written as 
         $ f\times f:X\times_{Y}X\to X'\times_{Y'}X'$ and $(f\times\id)\circ\iota^X:
        X\times_YX\to X'\times_SX$: 
        \[
        \xymatrix{
        (X\times_YX\simeq) X'\times_{Y'} X\ar[r]^-{\id\times f}\ar[d]&X'\times_{Y'}X'\ar[d]^-{\iota^{X'}}\\
        X'\times_SX\ar[r]^-{\id\times f}&X'\times_SX'
        }
        \]
        Then, by base change, we can write
        \[
        K=(f\times\id)_*\iota^X_*(f\times f)^*\cP_{X'}.
        \] Since the $\bG_m$-gerbe $X\to Y$ with trivialized banding is the pullback of $X'\to Y'$, we have $(f\times f)^*\cP_{X'}\simeq \cP_{X}$, so that 
        \[
        K\simeq (f\times\id)_*\iota^X_*\cP_{X}.
        \]
    \end{itemize}
\end{example}
\section{$\infty$-groups in derived stacks}\label{section-infty-grp-dSt}
\subsection{Construction of the fiber sequences in {Example \ref{example-gerbe-banding-fiber-sequence}}}
In this subsection, we explain how the first fiber sequence in Example \ref{example-gerbe-banding-fiber-sequence} can be constructed. (The second one is constructed similarly {working with derived stacks in the fppf topology. In particular, the truncation in this subsection is also computed in the fppf sense.}) 
Consider the $\infty$-topos $\bH=\dSt_{/\bZ}$.
For objects $X,Y\in\bH$, we denote the internal mapping stack by $\sMap_{\bH}(X,Y)$. For any $Z\in\bH$, we have (see for example \cite{Dan-HL-Preygel})
\[
\Map_{\bH}(Z,\sMap_{\bH}(X,Y))\simeq\Map_{\bH}(X\times Z,Y)\simeq \Map_{\bH_{/Z}}(X_Z,Y_Z).
\]
We can also define a derived stack $\sIsom_{\bH}(X,Y)$ by the formula
\[
\sIsom_{\bH}(X,Y)(S):=\Map_{(\bH_{/S})^{\simeq}}(X_S,Y_S)
\] for any derived affine $S$ (see \cite{Gepner-Kock} where this is called the ``internal object of equivalences''). When $X=Y$, this is also denoted by $\sAut_{\bH}(X)=\sIsom_{\bH}(X,X).$ 

Let $G,H$ be $\infty$-groups in $\bH$.
We can define a presheaf $\sMap_{\Grp(\bH)}(G,H)$ in $\cS$ by the formula
\[
\sMap_{\Grp(\bH)}(G,H)(S):=\Map_{\Grp(\bH_{/S})}(G_S,H_S)
\]
for any derived affine $S$ over $\bZ$.
Similarly to above, we define $\sIsom_{\Grp(\bH)}(G,H),\sAut_{\Grp(\bH)}(G)$.

As in \cite{NSS}, the looping and delooping functors define an equivalence between $(\bH_{/S})_{_{\geq1}}^{S/}$ and $\Grp(\bH_{/S})$, so using 
$\Map_{\Grp(\bH_{/S})}(G_S,H_S)\simeq\Map_{(\bH_{/S})_{_{\geq1}}^{S/}}(\bB G_S,\bB H_S)$, 
we have a fiber sequence
\[
\Map_{\Grp(\bH_{/S})}(G_S,H_S)\to \Map_{\bH_{/S}}(\bB G_S,\bB H_S)\to \Map_{\bH_{/S}}(S,\bB H_S)),
\] for each $S$. 
By the compatibility of delooping with base change, 
we have a fiber sequence of presheaves (also see \cite[Example 1.2.4]{Dan-HL-Preygel}):
\[
\sMap_{\Grp(\bH_{})}(G,H)\to \sMap_{\bH_{}}(\bB G,\bB H)\to \sMap_{\bH}(\Spec\bZ,\bB H))\simeq \bB H.
\]
This shows in particular that $\sMap_{\Grp(\bH_{})}(G,H)$ is a derived stack.
We can repeat the same argument for $\Map_{(\Grp(\bH_{/S}))^{\simeq}}(G_S,H_S)\simeq\Map_{((\bH_{/S})_{_{\geq1}}^{S/})^{\simeq}}(\bB G_S,\bB H_S)$ to obtain a fiber sequence
\[
\sIsom_{\Grp(\bH_{})}(G,H)\to \sIsom_{\bH_{}}(\bB G,\bB H)\to \sMap_{\bH}(\Spec\bZ,\bB H))\simeq \bB H,
\]
showing that $\sIsom_{\Grp(\bH_{})}(G,H), \sAut_{\Grp(\bH_{})}(G)$ are derived stacks.
\begin{remark}\label{remark-Isom(BG)-Isom-grp(G)-fiber-sequence}
    To emphasize that restricting to the $\sIsom$ still yields a fiber sequence, let us quickly motivate this construction. We explain the following pullback square (and a similar pullback square exists after base change of $\Spec\bZ$):
    \[
    \xymatrix{
\Isom_{\Grp(\bH)}(G,H)\ar[r]\ar[d]&\Map(\Spec\bZ,\Spec\bZ)\ar[d]\\
    \Isom_{\bH}(\bB G,\bB H)\ar[r]&\Map(\Spec\bZ,\bB H)
    }
    \]
    We already know this exists if we work with $\Map$ instead of $\Isom$, and the version for $\Map$ works as follows. Given a point of $\Map_{\Grp(\bH)}(G,H)$, we have a morphism between the simplicial objects of $G$ and $H$ describing the group structures (which is the identity on $\Spec \bZ$), and simply taking the colimit of these simplicial objects gives a diagram: \[
    \xymatrix{
\Spec\bZ\ar[r]^-\id\ar[d]^-e&\Spec\bZ\ar[d]^-e\\
    \bB G\ar[r]&\bB H
    }
    \]
    Here, $e$ denotes the base points of $\bB G$ and $\bB H.$  
    This construction determines a point of $\Map_{\bH}(\bB G,\bB H) \times_{\Map(\Spec\bZ,\bB H)} \Map(\Spec\bZ,\Spec\bZ)$. 
    Conversely a point of \[\Map_{\bH}(\bB G,\bB H)\times_{\Map(\Spec\bZ,\bB H)}\Map(\Spec\bZ,\Spec\bZ)\]
    consists of $\alpha:\bB G \to \bB H$, $\id:\Spec\bZ\to\Spec\bZ$ and a path between 
    $\Spec\bZ\xrightarrow{e}\bB G \xrightarrow{\alpha} \bB H$
    and $\Spec\bZ\xrightarrow{\id}\Spec\bZ\xrightarrow{e} \bB H$, which can be summarized by the (homotopy) commutative diagram:
    \[\xymatrix{
\Spec\bZ\ar@{=}[r]\ar[d]^-e&\Spec\bZ\ar[dd]^-e\\
\bB G\ar[d]^-\alpha&\\
    \bB H\ar@{=}[r]&\bB H
    }\]
    This can be reformatted to resemble the previous square diagram, and taking the Cech nerves of vertical maps $e$ gives a map of simplicial objects, which describes a group homomorphism $G\to H$. 
    Colimits of equivalences are equivalences, and pullbacks of equivalences are equivalences, so indeed we can restrict to $\Isom.$
\end{remark}
\begin{example}
\leavevmode
\begin{enumerate}
    \item When $G=\bG_m$, $\sMap_{\Grp(\bH)}(\bG_m,\bG_m)\simeq \underline{\bZ}_{\bZ},$ and $ \sAut_{\Grp(\bH)}(\bG_m)\simeq \underline{\bZ/2\bZ}_{\bZ}.$
    \item When  $G=\mu_n$, 
$\sMap_{\Grp(\bH)}(\mu_n,\mu_n)\simeq \underline{\bZ/n\bZ}_{\bZ}$, and 
$\sAut_{\Grp(\bH)}(\mu_n)\simeq \underline{(\bZ/n\bZ)^\times}_{\bZ}$.
\end{enumerate}

\end{example}

As a result, we have a fiber sequence of derived stacks
\begin{equation}\label{fiber-sequence-Aut(BGm)-fixing-point}
\underline{\bZ/2\bZ}_{\bZ}\xrightarrow[]{}
\sAut_{\bH}(\bB{\bG_m})\xrightarrow{e^{\#}} 
\bB \bG_m.
\end{equation}
{
 Let $G$ be an $\infty$-group in $\bH$, and choose a point $\phi:\Spec \bZ\to \sAut_{\Grp(\bH)}(G)$.
 Consider the diagram with cartesian squares, where the right square is constructed as in the beginning of the section.
    \begin{equation}\label{diagram-conjugation-action}
\xymatrix{
A\ar[r]\ar[d]&
\sAut_{\Grp(\bH)}(G)\ar[r]\ar[d]&
{\Spec\bZ}\ar[d]
\\
{\Spec\bZ}\ar[r]^-{\bB\phi}&\sAut_{\bH}(\bB G)\ar[r]^-{e^{\#}}&\bB G
}
    \end{equation}
By viewing $\phi$ as an $\infty$-group automorphism $\phi:G\to G$,
we are writing the resulting point on $\sAut_{\bH}(\bB G)$ as $\bB \phi:\Spec\bZ\to \sAut_{\bH}(\bB G)$ (which corresponds to $\bB\phi:\bB G\to\bB G$).
In view of the right pullback square, 
the point $\phi$ determines a path from $\Spec\bZ\xrightarrow{e}\bB G\xrightarrow{\bB \phi}\bB G$ to 
$e:\Spec\bZ\to \bB G$. We write $\widetilde{\phi}$ for the 2-cell:
\[
\xymatrix{
\Spec\bZ\ar@{=}[r]\ar[d]^-e&\Spec\bZ\ar[d]^-e\\
\bB G\ar[r]^-{\bB\phi}\ar@{=>}[ur]^-{\widetilde{\phi}}&\bB G
}
\]
$\widetilde{\phi}$ furnishes a path  between the bottom horizontal composition $e^{\#}\circ\bB\phi$  and $e:\Spec\bZ\to\bB G$, so we obtain an identification $\mu: A\simeq \Spec\bZ\times_{e,\bB G,e}\Spec\bZ\simeq G$.
\begin{lemma}\label{lemma-connecting-map-fiber-sequence-conjugation}
   Up to the identification $\mu:A\simeq G$ from above, the morphism $A\to\sAut_{\Grp(\bH)}(G)$ from (\ref{diagram-conjugation-action}) is the conjugation action map on $\phi:G\to G$.
\end{lemma}
\begin{proof}
We describe this by evaluating on the derived affine $\Spec\bZ$ (evaluations on other derived affines are similar after base changes). Applying $\Map(\Spec\bZ,-)$ to the diagram, we obtain a diagram with cartesian squares: 
        \[
\xymatrix{
\Map(\Spec\bZ,A)\ar[r]^-c\ar[d]&
\Aut_{\Grp(\bH)}(G)\ar[r]\ar[d]&
{*}\simeq\Map({\Spec\bZ},{\Spec\bZ})\ar[d]
\\
{*}\simeq\Map(\Spec\bZ,\Spec\bZ)\ar[r]&\Aut(\bB G)\ar[r]&\Map(\Spec\bZ,\bB G)
}   \]
In view of the big cartesian square, a point $g_0:*\to  \Map(\Spec\bZ,A)$ of the top left corner corresponds to a point of $*\times_{\Map(\Spec\bZ,\bB G)}*$, which records the (homotopy) commutative diagram:
\begin{equation}\label{diagram-g-2-cell}
\xymatrix{
\Spec\bZ\ar@{=}[r]\ar[d]^-e&\Spec\bZ\ar[d]^-e\\
\bB G\ar[r]^-{\bB\phi}\ar@{=>}[ur]^-{\widetilde{g_0}}&\bB G
}
\end{equation}
As in Remark \ref{remark-Isom(BG)-Isom-grp(G)-fiber-sequence}, 
taking the Cech nerves of the vertical maps here yields an $\infty$-group morphism $c(g_0):G\to G$ which is the image of the point $g$ on $\Aut_{\Grp(\bH)}(G)$ via $c$. 
To describe $c(g_0)$, we consider a point $h:\Spec\bZ\to G$ and study how $c(g_0)$ sends this point to another $\bZ$-point $c(g_0)\circ h$ of $G$ (Similar argument works for other points by base change as well, but we focus on $\bZ$ for simplicity.) 
We view $h\in\Map(\Spec\bZ,G)\simeq *\times_{\Map(\Spec\bZ,\bB G)}*$ as defining a (homotopy) commutative diagram:
\[
\xymatrix{
\Spec\bZ\ar@{=}[r]\ar[d]^-e&\Spec\bZ\ar[d]^-e\\
\bB G\ar@{=}[r]\ar@{=>}[ur]^-{\widetilde{h}}&\bB G
}
\]
The morphism $c(g_0):\Spec\bZ\times_{\bB G}\Spec\bZ\to \Spec\bZ\times_{\bB G}\Spec\bZ$ is defined by taking two copies of (\ref{diagram-g-2-cell}) and considering appropriate pullbacks, so
as a point of $*\times_{\Map(\Spec\bZ,\bB G)}*(\simeq \Map(\Spec\bZ,G))$, 
$c(g)\circ h$ corresponds to the big combined 2-cell in:
\begin{equation}\label{diagram-c(g0)-circ-h}
\xymatrix{
&\Spec\bZ\ar@{=}[dr]\ar[dd]^-e\\
\Spec\bZ\ar@{=}[ru]\ar@{=}[dd]&&\Spec\bZ\ar[dd]^-e\\
&\bB G\ar[dr]^-{\bB\phi}\ar@{=>}[dd]^-{\widetilde{g_0}}\ar@{=>}[ru]^-{\widetilde{g_0}}
&\\
\Spec\bZ\ar[ru]^-e\ar@{=}[dr]\ar@{=>}[uuur]^-{\widetilde{h}}&&\bB G\\
& \Spec\bZ\ar[ru]^-e
}
\end{equation}
We need to describe this as a conjugation up to the trivialization $\mu$. 
Let \[g:*\to\Map(\Spec\bZ,G)\simeq\Map(\Spec\bZ,\Spec\bZ)\times_{\Map(\Spec\bZ,\bB G)}\Map(\Spec\bZ,\Spec\bZ)\] be the point identified with $g_0$ under $\mu.$ If we present $g$ as the (homotopy) commutative square
\[
\xymatrix{
\Spec\bZ\ar@{=}[r]\ar[d]^-e&\Spec\bZ\ar[d]^-e\\
\bB G\ar@{=}[r]\ar@{=>}[ur]^-{\widetilde{g}}&\bB G
}
\]
under the identification $\mu$,
the 2-cell $\widetilde{g_0}$ from earlier arises as the big square in
\begin{equation}\label{diagram-mu-identification-squares}
\xymatrix{
\Spec\bZ\ar@{=}[r]\ar[d]&
\Spec\bZ\ar@{=}[r]\ar[d]^-e&\Spec\bZ\ar[d]^-e\\
\bB G\ar[r]^-{\bB \phi}
\ar@{=>}[ur]^-{\widetilde{\phi}}
&
\bB G\ar@{=}[r]\ar@{=>}[ur]^-{\widetilde{g}}&\bB G
}
\end{equation}
With this preparation, the claim is that the big diagram (\ref{diagram-c(g0)-circ-h}) coincides with $\widetilde{(g^{}(\phi\circ h)g^{-1})}$, where $\phi\circ h\in\Map(\Spec\bZ ,G)$.
This is not hard to see if we substitute the instances of $\widetilde{g_0}$ in (\ref{diagram-c(g0)-circ-h}) by (\ref{diagram-mu-identification-squares}) and note that the following big 2-cell represents $\widetilde{\phi\circ h}$.
\[
\xymatrix{
&\Spec\bZ\ar@{=}[dr]\ar[dd]^-e\\
\Spec\bZ\ar@{=}[ru]\ar@{=}[dd]&&\Spec\bZ\ar[dd]^-e\\
&\bB G\ar[dr]^-{\bB\phi}\ar@{=>}[dd]^-{\widetilde{\phi}}\ar@{=>}[ru]^-{\widetilde{\phi}}
&\\
\Spec\bZ\ar[ru]^-e\ar@{=}[dr]\ar@{=>}[uuur]^-{\widetilde{h}}&&\bB G\\
& \Spec\bZ\ar[ru]^-e
}
\]

\end{proof}
\begin{remark}\label{remark-conjugation-action-map-other-points}
    In fact, the same argument works after base change from $\Spec \bZ$ to another derived affine, so that an analogous conclusion holds with respect to $S$-points on $\sAut_{\Grp(\bH)}(G)$ for a derived affine $S$.
\end{remark}
}
\begin{corollary}\label{corollary-truncation-fiber-sequence-still-fiber-sequence}
    The truncation of the fiber sequence (\ref{fiber-sequence-Aut(BGm)-fixing-point}) which we obtain by applying $\tau_{\leq0}$ is still a fiber sequence. Since $\tau_{\leq0}(\bB\bG_m)=\Spec\bZ$, we see in particular that $\tau_{\leq0}(\sAut(\bB\bG_m))\simeq \ZModTwo_{\bZ}$.
\end{corollary}
\begin{proof}
Consider any test derived affine $T$ and look at the fiber sequence of spaces
\[
\Map(T,\ZModTwo_{\bZ})\to\Map(T,\sAut(\bB \bG_m))\to\Map(T,\bB\bG_m).
\]
We can rewrite this as 
$\Map_T(T,\ZModTwo_{T})\to\Map_T(T,\sAut_{T}((\bB \bG_m)_T))\to\Map_T(T,(\bB\bG_m)_T)$
We claim that applying $\pi_0$ to this sequence yields a fiber sequence of sets. By homotopy long exact sequence, it suffices to choose a point  $\phi:*\to \Map_T(T,\ZModTwo_{T})\simeq\Aut_{\Grp(\bH_{/T})}((\bG_m)_T)$ and check that the action of 
\[\pi_1(\Map_T(T,(\bB\bG_m)_T))\simeq\pi_0(\Omega\Map_T(T,(\bB\bG_m)_T))\simeq \pi_0(\Map_T(T,\Omega_e(\bB\bG_m)_T))\]
on the point $\phi$ is trivial. By 
    Lemma \ref{lemma-connecting-map-fiber-sequence-conjugation} and Remark \ref{remark-conjugation-action-map-other-points}, this action is given by the conjugation of $(\bG_m)_T$ on $\phi:(\bG_m)_T\to(\bG_m)_T$, so this justifies the claim. 

    The result so far can be summarized as the fiber sequence of discrete presheaves
    \[
    \tau_{\leq0}^{\PSh}(\ZModTwo_{\bZ})\to \tau_{\leq0}^{\PSh}(\sAut(\bB\bG_m))\to\tau_{\leq0}^{\PSh}(\bB\bG_m),
    \] where $\tau_{\leq0}^{\PSh}(\bB\bG_m)$ denotes the presheaf sending a derived affine $S$ to $\pi_0(\bB\bG_m(S))$ (and we define  $\tau_{\leq0}^{\PSh}(\ZModTwo_{\bZ}),\tau_{\leq0}^{\PSh}(\sAut(\bB\bG_m))$ similarly). $\infty$-sheafification preserves fiber sequences, so the proof is complete.
\end{proof}
Under the equivalence from Corollary \ref{corollary-truncation-fiber-sequence-still-fiber-sequence}, the composition of $\ZModTwo_\bZ\to \sAut(\bB\bG_m)$ with
the
$\infty$-group morphism to the $0$-th truncation (Remark \ref{remark-truncation-group})
\[
T:\sAut_{\bH}(\bB\bG_m)\to \tau_{\leq0}\sAut_{\bH}(\bB\bG_m)\simeq \underline{\bZ/2\bZ}_{\bZ}.
\] is an equivalence.

\begin{lemma}\label{lemma-effective-epi-group-morphism-triv-fiber}
    Let $f:G\to H$ be a morphism of $\infty$-groups with a trivial fiber (i.e.\ the final object). If $f$ is an effective epimorphism, $f$ is an equivalence.
\end{lemma}
\begin{proof}
    If $F$ denotes the fiber, we have $G\times_HG\simeq F\times G$ and a similar computation holds for the higher terms of the Cech nerve of $f$. However, $F$ is trivial by assumption, so the Cech nerve reduces to the constant diagram of $G$. When $f$ is an effective epimorphism, the colimit diagram of the Cech nerve is precisely $f$, so $f$ is an equivalence. 
\end{proof}
\begin{proposition}
    $\fib(T)\simeq\bB\bG_m$ as $\infty$-groups.
\end{proposition}
\begin{proof}
    Using the $\infty$-group structure of $\bB\bG_m$, we have an $\infty$-group morphism which is the translation action map:
    $
    s:\bB\bG_m\to\sAut_{\bH}(\bB\bG_m).
    $ 
    By construction, this is a section of the morphism $e^{\#}:\sAut_{\bH}(\bB\bG_m)\to\sMap_{\bH}(\Spec\bZ,\bB\bG_m)\simeq\bB\bG_m$. 
    Moreover, the composition of $s$ with $T$ is equivalent to the trivial point $\bB\bG_m\to\Spec\bZ\to\underline{\bZ/2\bZ}_{\bZ}$ because $\tau_{\leq0}(\bB\bG_m)\simeq\Spec\bZ$:
    \[
    \xymatrix{
\bB\bG_m\ar[r]\ar[d]^{}&\underline{\bZ/2\bZ}_{\bZ}\ar[d]^-\simeq\\
   \tau_{\leq0}(\bB\bG_m)\simeq\Spec\bZ\ar[r] &\tau_{\leq0}\underline{\bZ/2\bZ}_{\bZ}
    }
    \]
    By universality of fiber sequence,
    $s$ factors through $s':\bB\bG_m\to \fib(T)$.
    The goal is to show $s'$ is an equivalence. 
    Since $s$ is a section, we have a fiber sequence
    \[
    \fib(s)\to \fib(\id_{\bB\bG_m})\to \fib(e^{\#}).
    \]
    This simplifies to 
    \[
    \fib(s)\to \Spec\bZ\to \ZModTwo_{\bZ},
    \] so $\fib(s)\simeq\Spec\bZ$. 
    Next, by the factorization of $s$ via $s'$, we have a fiber sequence of $\infty$-groups
    \[
    \fib(s')\to \fib(s)\to\fib(\fib(T)\to\sAut_{\bH}(\bB\bG_m)).
    \]
    In view of the following diagram of cartesian squares, we see that $\fib(\fib(T)\to\sAut_{\bH}(\bB\bG_m))\simeq \Omega_e\ZModTwo_{\bZ}\simeq\Spec\bZ$.
\[
\xymatrix{
\fib(\fib(T)\to\sAut_{\bH}(\bB\bG_m))\ar[d]\ar[r]&\fib(T)\ar[d]\ar[r]&\Spec\bZ\ar[d]\\
\Spec\bZ\ar[r]&\sAut_{\bH}(\bB\bG_m)\ar[r]^-{T}&\ZModTwo_{\bZ}
}
\]
Thus, the right two terms in the above fiber sequence are trivial. Hence, $\fib(s')=\Spec\bZ$. 
We also claim that $s'$ is an effective epimorphism.
By \cite[Proposition 2.10]{NSS}, we can check this after truncation, and it suffices to check that $\tau_{\leq0}(\fib(T))\simeq \Spec\bZ$. For this, we may show that the truncation $\tau_{\leq0}$ of the fiber sequence $\fib(T)\to \sAut(\bB\bG_m)\xrightarrow{T} \tau_{\leq0}(\sAut(\bB\bG_m))\simeq\ZModTwo_{\bZ}$ is still a fiber sequence. In principle, this is checked similarly to Corollary \ref{corollary-truncation-fiber-sequence-still-fiber-sequence}, but this case is much simpler since $\Omega\ZModTwo_{\bZ}\simeq\Spec
\bZ$.
In summary, $s'$ is a morphism of $\infty$-groups, which is also 
an effective epimorphism with trivial fibers. By Lemma \ref{lemma-effective-epi-group-morphism-triv-fiber}, we are done.

\end{proof}
\begin{corollary}
    We have a fiber sequence of $\infty$-groups induced by the truncation:
    \[
    \bB\bG_m\to \sAut_{\bH}(\bB\bG_m)\xrightarrow[]{T} \ZModTwo
    _{\bZ}.
    \]
\end{corollary}
The case of $G=\mu_n$ is similar. 

\subsection{Clarifications on bandings}
Let $\bH:=(\dSt_{/S})_{/Y}$ in this section, where $Y$ is a derived stack over $S.$ 
Let $G$ be an $\infty$-group in $\bH.$ We will also consider $\sAut_{\Grp(\bH)}(G)$ which is again an $\infty$-group in $\bH$. In this subsection, we briefly clarify that the bandings defined as in \cite{Binda-Porta} and those as in Definition \ref{definition-banding-NSS-version} coincide for $G=\bG_m,\mu_n$, where $n$ is some positive integer. 
\begin{definition}\label{definition-twisted-form}
A \textbf{twisted form} $T$ of $G$ in $\bH$ is an $\infty$-group $T$ in $\bH$ such that for some effective epimorphism $U\to Y$, there exists an equivalence
$T_U\simeq G_U$ of
$\infty$-groups in $\bH_{/U}$.
We let $\TwForm_S(G)$ denote the $\infty$-groupoid of twisted forms in $\bH$ (defined as a full subcategory of the $\infty$-groupoid $\Grp(\bH)^{\simeq}$). 
\end{definition}
Analogously to Proposition \ref{proposition-gerbe-BAutBG}, we have:
\begin{proposition}\label{proposition-twForm-BAutGrp}
    We have an equivalence 
    \[
    \Map(Y,\bB\sAut_{\Grp(\bH)}(G))\simeq \TwForm_Y(G).
    \]
\end{proposition}

Let $G=\bG_m$ or $G=\mu_n$ below. 
Given a $G$-gerbe $f:X\to Y$ {over a derived stack $Y$}, by the same argument as \cite[\S5.2]{Binda-Porta}, we can associate a banding $ f_*I_{X/Y}$ which is a twisted form of $G$ over $Y$. On the other hand, we can also construct a twisted form of $G$ over $Y$ in the abstract fashion as follows. 
By Corollary \ref{corollary-truncation-fiber-sequence-still-fiber-sequence} (and its analogue for $\mu_n$), the composition $\rho:\sAut_{\Grp(\bH)}(\bG)\to \sAut_{\bH}(\bB G)\xrightarrow[]{T}  \Out(G)$ is an equivalence. 
The $G$-gerbe $f$ corresponds to a point of $\Map(Y,\bB \sAut_{\bH}(\bB G))$ (Proposition \ref{proposition-gerbe-BAutBG}), but the truncation $\infty$-group morphism $T:\sAut_{\bH}(\bB G)\to\Out(G)$ and $\rho^{-1}$ give us 
$\bB(\rho^{-1}\circ T)\circ(-):
\Map(Y,\bB \sAut_{\bH}(\bB G))\to
\Map(Y,\bB\sAut_{\Grp(\bH)}(\bG))$.
By Proposition \ref{proposition-twForm-BAutGrp}, the image of the point along this map corresponds to a twisted form of $G$ over $Y$. With somewhat tedious formal computations one can show that these constructions coincide:
\begin{proposition}\label{proposition-compare-banding}
For $G=\bG_m$ or $\mu_n$,
the following diagram in $\cS$ commutes up to homotopy:
    \[
    \xymatrix{
    \Map(Y,\bB \sAut_{\bH}({\bB G}))\ar[rr]^-{\bB(\rho^{-1}\circ \tau)\circ(-)}\ar[d]^-{\simeq}&&\Map(Y,\bB\sAut_{\Grp(\bH)} (G))\ar[d]^-{\simeq}\\
\Gerbe_Y(G)\ar[rr]&&\TwForm_Y(G)\\
    }
    \]
    where the bottom horizontal morphism is constructed by \cite[\S5.2]{Binda-Porta}, and the vertical equivalences come from Proposition \ref{proposition-gerbe-BAutBG} and Proposition \ref{proposition-twForm-BAutGrp}.
    \end{proposition}

\end{appendix}
\part{Atiyah classes}
\setcounter{section}{\value{appendixFix}}
In this chapter, we construct and compute multiple variants of Atiyah classes. 
The most basic version assumes the following setup. $X$ is a derived stack with a cotangent complex over a field $k$, and $E\in\Dqc(X)$ is a complex, and we choose a morphism $d:\bL_{X/k}\to E$. In this situation, we will construct the \textbf{universal $d$-Atiyah transformation of $X$} $\At_d(-):\id\to (-)\otimes E[1]$, which is a natural transformation of functors $\id,(-)\otimes E[1]:\Dqc(X)\to\Dqc(X)$. One can evaluate at each $F\in\Dqc(X)$ to obtain a morphism $\At_d(F):F\to F\otimes E[1]$ in $\Dqc(X)$ which we call the \textbf{$d$-Atiyah class of $F$ on $X$.} When $X$ is a reasonable variety and $d=\id_{\bL_{X/k}}$, this recovers the Atiyah class of $F$ in the classical sense, but we can choose different choices of $d$, allowing for ``derivations'' other than the universal one.

In the literature, Kuhn (\cite{Kuhn}) already constructed Atiyah classes for (classical) algebraic stacks using simplicial methods, 
and for geometric derived stacks in characteristic 0, there exists a construction of   Atiyah classes via the moduli of perfect complexes (see \cite{STV} and \cite{Hennion}).
Our approach is different and is  closer to the viewpoints in \cite{SAG} and \cite{Moore2}. It is worth mentioning that our construction relies on the theory of nonconnective derived stacks (or $\cC$-stacks in the terminology of \cite{blowup-normal-bundle-nonconnective-DAG}). The reason is that the construction analogous to \cite{SAG} and \cite{Moore2} uses the square-zero extension $X[\bL_{X/k}]$ of $X$ by $\bL_{X/k}$. Since such cotangent complexes of stacks are not connective in general, $X[\bL_{X/k}]$ will not be a derived stack. Instead, one can define it as a nonconnective derived stack using the techniques of \cite{blowup-normal-bundle-nonconnective-DAG}. In \cite{Raksit}, Raksit developed the theory of (nonconnective) derived commutative rings, originally due to Bhatt-Mathew, and our stacks are roughly sheaves on such nonconnective affine schemes. 

We will additionally construct a twisted analogue of  Atiyah class which makes sense for twisted sheaves (\S\ref{section-atiyah-class-for-twisted-sheaves}). Furthermore, we will define a noncommutative version of twisted Atiyah classes (Definition \ref{definition-noncommutative-twisted-atiyah-map}), and this is essentially a twisted version of what Moore calls the ``derived Atiyah transformation'' in \cite[\S3.1]{Moore2}.

One key result is that the (commutative) twisted Atiyah class and the noncommutative twisted Atiyah class are compatible (Proposition \ref{proposition-commutative-vs-noncommutative-atiyah-maps}), which is a twisted version of \cite[Proposition 3.12]{Moore2}. This will later be used to compute the twisted-to-noncommutative map on tangent fibres. 

In the last section (\S\ref{section-twisted-atiyah-transformation-trivial-principal-BGm-bundles}), we will compute the twisted Atiyah transformation for trivial $\bB\bG_m$-bundles, and relate it to the HKR morphism. This is analogous to the untwisted result \cite[Theorem 3.1]{Moore2}.

\section{Atiyah class for derived stacks}\label{section-Atiyah-class-for-derived-stacks}
\subsection{Construction}
Let $X$ be a 
derived stack which admits a cotangent complex over a field
$k$ (e.g.\ a 1-Artin derived stack over $k$) 
and let $E\in\Dqc(X)$ be a complex. 
In this section, we will borrow notations from \cite{Raksit} and \cite{blowup-normal-bundle-nonconnective-DAG} such as $\Spec^{\nc}$ and $\QCoh$. 
 As in
 \cite[pp.31]{blowup-normal-bundle-nonconnective-DAG}, we can define the non-connective derived stack $X[E]:=\Spec_X^{\nc}(\cO_X\oplus E)$ using non-connective relative spectrum of a square-zero extension in the $\QAlg$-sense. Below,
\[
i:X\to X[E], \quad u:X[{E[1]}]\to X, \quad\pi:X[E]\to X
\] are the obvious morphisms (which in the classical setup are the closed immersion and projections).
Let us briefly explain how to obtain the morphism $\pi$ for instance.
In the non-connective affine case where $X=\Spec^{\nc}(A)$, we continue to write $E$ to represent the $\Theta(A)$-module corresponding to $E\in\QCoh(X)$ (where we recall the functor $\Theta:\DAlg(\cC)\to\CAlg(\cC)$ from \cite{Raksit} and \cite{blowup-normal-bundle-nonconnective-DAG}). Then, in the notations of \cite[\S4.4]{Raksit}, the square-zero extension $A\oplus E$ arises from $(A,E)\in\DAlgMod(\cC)$, where $\cC=\Mod_{\bZ}$
 is the derived algebraic context of \cite[Example 4.3.1]{Raksit}. The trivial extension $A=A\oplus 0$ arises from $(A,0)\in\DAlgMod(\cC)$, so we can define a morphism $A\to A\oplus E$ from the morphism $(A,0)\to(A,E)$ in $\DAlgMod(\cC)$. Applying $\Spec^{\nc}$ gives us $\pi:\Spec^{\nc}(A\oplus E)\to \Spec^{\nc}(A)$. The non-affine case can be done similarly by defining the map nonconnective-affine-wise.
 
\begin{remark}\label{remark-cotangent-complex-nonconnective-universality}
    Even in the non-connective setting, cotangent complexes enjoy familiar properties, as shown in \cite{blowup-normal-bundle-nonconnective-DAG}. By definition of cotangent complex (\cite[Definition 3.50]{blowup-normal-bundle-nonconnective-DAG}), we have an equivalence 
    $
    \Der_{X/k}\simeq h_{(X,\bL_{X/k})}
    $ of presheaves on $\StMod_k$, where the objects of $\StMod_k$ consist of pairs $(T,M)$ where $T\to\Spec k$ is a non-connective derived stack, and $M\in\QCoh(T)$. Also, by definition, $\Der_{X/k}(T,M):=\Map_k(T[M],X)$.
    For any object $(T,M)\in\StMod_k$, we have a commutative diagram of spaces.
    \[
    \xymatrix{
    \Der_{X/k}(T,M)\ar[r]^-\simeq\ar[d]&h_{(X,\bL_{X/k})}(T,M)\ar[d]\\
    \Der_{X/k}(T,0)\ar[r]^-\simeq&h_{(X,\bL_{X/k})}(T,0)
    }
    \]
    In particular, we have:
    \[
    \xymatrix{
    \Map_k(X[E],X)\ar[r]^-\simeq\ar[d]&\Map_{\StMod_k}((X,E),(X,\bL_{X/k}))\ar[d]\\
    \Map_k(X,X)\ar[r]^-\simeq&\Map_{\StMod_k}((X,0),(X,\bL_{X/k}))
    }
    \]
    By taking the fibers of the vertical maps, we obtain 
    \[
    \Map_{X//k}(X[E],X)\simeq \Map(\bL_{X/k},E),
    \] which is functorial in $E$. 
    We will freely use this fact below. 
\end{remark}
With the above remark in mind, we additionally choose a point \[\alpha\in\Map'_k(X[E],X):=\Map_k(X[E],X)\times_{\Map_k(X,X)}\{\id_X\}\simeq\Map(\bL_{X/k},E) ,\] where the map $\Map_k(X[E],X)\to \Map_k(X,X)$ is the composition with $i$. 

As in \cite[pp.23]{blowup-normal-bundle-nonconnective-DAG}, for a non-connective derived stack $X$, $\QCoh(X)$ is stable presentable and has a natural symmetric monoidal structure. Moreover, we have the usual adjunction of pushforward and pullback for morphisms between such stacks.
\begin{lemma}[{\cite[Lemma 3.25]{blowup-normal-bundle-nonconnective-DAG}. Also see \cite[1.3.10]{Drinfeld-Gaitsgory}}]\label{lemma-nonconnective-base-change-projection-formula}
Nonconnectively affine morphisms satisfy arbitrary base change and projection formula.
More precisely, consider a nonconnectively affine morphism $f:X\to Y$ of nonconnective derived stacks, and a pullback diagram
\[
\xymatrix{
X'\ar[r]^-{g_X}\ar[d]^-{f'}&X\ar[d]^-f\\
Y'\ar[r]^-{g_Y}&Y
}
\]
The canonical morphism
\[
g_Y^*f_*E\to f'_*g_X^* E
\] is an isomorphism. Moreover, for any $E\in\QCoh(X),F\in\QCoh(Y)$, the canonical morphism
\[
f_*E\otimes F\to f_*(E\otimes f^*F)
\] is an isomorphism.
\end{lemma}
\begin{proof}
    The base change statement is proved in \cite[Lemma 3.25]{blowup-normal-bundle-nonconnective-DAG}. Projection formula is clear when $Y$ is nonconnectively affine, and the general case reduces to this by base change (also see \cite[1.3.10]{Drinfeld-Gaitsgory}).
\end{proof}
\begin{remark}\label{remark-square-zero-DAlg-CAlg}
Recall from \cite[Construction 4.4.5]{Raksit} that there is a square-zero extension functor $G:\DAlgMod(\cC)\to \DAlg(\cC)$, where $\cC=\Mod_{\bZ}$ is the derived algebraic context of \cite[Example 4.3.1]{Raksit}.
One can similarly define a $\CAlg$-version of the square-zero extension functor $G_{\CAlg}:\CAlgMod(\cC)\to\CAlg(\cC) $ using the equivalence $\alpha:\Gr^{\{0,1\}}\CAlg(\cC)\simeq\CAlgMod(\cC)$ from \cite[Lemma 4.4.2.]{Raksit}, and by construction it makes the following diagram commute:
\[
\xymatrix{
\DAlgMod(\cC)\ar[r]^-G\ar[d]^-{(\Theta,\id)}&\DAlg(\cC)\ar[d]^-{\Theta}\\
\CAlgMod(\cC)\ar[r]^-{G_{\CAlg}}&\CAlg(\cC)
}
\]
\end{remark}
Assume for a moment that $X=\Spec^{\nc}(A)$. Then, $\QCoh(X[E])\simeq\Mod_{A\oplus E}$ which by definition is $\Mod_{\Theta(A\oplus E)}$. By Remark \ref{remark-square-zero-DAlg-CAlg}, this is the same as $\Mod_{\Theta(A)\oplus E}$, where $\Theta(A)\oplus E:=G_{\CAlg}(\Theta(A),E).$
Then, we obtain a fiber sequence of $\Theta(A)\oplus E$-modules 
\[
E\to \Theta(A)\oplus E\to \Theta(A).
\]
\begin{remark}
While this is not explicitly defined in \cite{blowup-normal-bundle-nonconnective-DAG}, the paper uses the structure sheaf $\cO_X$ which is an object of $\QAlg(X)$. One way to describe this is as follows. Since $\QAlg(X)$ is defined by right Kan extension, we have 
$\QAlg(X)\simeq \lim_{\{\Spec^{\nc}A\to X\}}\DAlg_A$. $\cO_X$ can be described as the compatible family of objects in the right hand side $\{A\}_{\{\Spec^{\nc}A\to X\}}.$ In the sequel, we often abuse notation and also write $\cO_X\in\QCoh(X)$ for the target of $\cO_X$ along the forgetful functor $\QAlg(X)\to\QCoh(X)$.
\end{remark}
\begin{lemma}\label{lemma-nonconnective-sqz-fiber-sequence}
We have a fiber sequence in $\QCoh(X[E])$
\[
i_*E\to\cO_{X[E]}\to i_* \cO_X.
\]
\end{lemma}
\begin{proof}
    By the discussion before the remark above, this holds when $X=\Spec^{\nc}A$ is affine. 
    In general, we write 
    \[
    X[E]\simeq X[E]\times_XX
    \simeq  X[E]\times_X\colim_{\{\Spec^{\nc} A\to X\}}\Spec^{\nc} A,
    \] and this is isomorphic to $ \colim_{\{\Spec^{\nc} A\to X\}}(X[E]\times_X\Spec^{\nc} A)$ since colimits are universal in an ($\infty$-)topos. Nonconnective relative spectra are compatible with base change (\cite[Lemma 3.23]{blowup-normal-bundle-nonconnective-DAG}), so we can write this as $ \colim_{\{\Spec^{\nc} A\to X\}}(\Spec^{\nc}(A)[E|_{\Spec^{\nc} A}])$.
    By \cite[Corollary 3.22]{blowup-normal-bundle-nonconnective-DAG}, we see 
    \[
    \QCoh(X[E])\simeq \lim_{\{\Spec^{\nc} A\to X\}^{\op}}\QCoh(\Spec^{\nc}(A)[E|_{\Spec^{\nc} A}]).
    \]
    On the right hand side, we have a compatible family of (co)fiber sequences in stable infinity categories $\QCoh(\Spec^{\nc}(A)[E|_{\Spec^{\nc} A}])$ by the affine case. Cofibers are computed pointwise in this limit $\infty$-category, 
     so the resulting diagram on the left hand side (which is the correct one in view of Lemma \ref{lemma-nonconnective-base-change-projection-formula})
    \[
i_*E\to\cO_{X[E]}\to i_* \cO_X.
\] is a (co)fiber sequence.
\end{proof}

By rotating the fiber sequence of Lemma \ref{lemma-nonconnective-sqz-fiber-sequence}, we obtain a morphism $\Psi:i_*\cO_X\to i_*E[1]$ in $\QCoh(X[E])$. This can be used to define a morphism
\[
\pi_*(\alpha^*(-)\otimes \Psi):\pi_*(\alpha^*(-)\otimes i_*\cO_X)\to \pi_*(\alpha^*(-)\otimes i_*E[1])
\] in
$\Fun(\QCoh(X),\QCoh(X)).$
We can invoke projection formula by Lemma \ref{lemma-nonconnective-base-change-projection-formula}, so this morphism can be rewritten as 
\[
\pi_*i_*(i^*\alpha^*(-))\to \pi_*i_*(i^*\alpha^*(-)\otimes E[1]).
\]
Using $\pi\circ i\simeq \id_X$ and the identification $\id_X\simeq\alpha\circ  i$, we simplify this as
\[
\pi_*(\alpha^*(-)\otimes\Psi):\id\to(-)\otimes E[1].
\]
Since $X$ is a (connective) derived stack, $\QCoh(X)$ coincides with $\Dqc(X)$ ({because any map from a  nonconnective affine $\Spec^{\nc} A\to X$ factors through the spectrum of its truncation to the connective part $\Spec\tau_{\geq0}A$}),
and this is a morphism in $\Fun_k^{\rL}(\Dqc(X),\Dqc(X))$, and we denote it by $\At_{\alpha}(-):\id\to(-)\otimes E[1]$.
\begin{definition}\label{definition-atiyah-transformation}\leavevmode
\begin{enumerate}
    \item The morphism $\At_{\alpha}(-):\id\to(-)\otimes E[1]$ in $\Fun^{\rL}_k(\Dqc(X),\Dqc(X))$ constructed above is called the \textbf{universal $\alpha$-Atiyah transformation of $X.$} 
    \item For a particular object $F\in\Dqc(X),$ the evaluation of the above transformation at $F$ is denoted by $\At_{\alpha}(F):F\to F\otimes E[1]$, and called the \textbf{$\alpha$-Atiyah class of $F$ on $X$.}
    \item  
    When $X$ is scalloped (qcqs 1-Artin), $\At_{\alpha}(-)$ corresponds to a morphism of Fourier-Mukai kernels via the untwisted analogue of Example \ref{example-twisted-FMT-via-dualizability}, and we denote it by $\At_{\alpha}^{\FM}\in\RHom_{X\times_kX}(\Delta_*\cO_X,\Delta_*E[1])$. This is called the\textbf{ universal $\alpha$-Atiyah class of $X.$} 
    \item If we fix $E$ and vary $\alpha$, via Remark \ref{remark-cotangent-complex-nonconnective-universality} the above construction gives us a map
\[
\bAt_E:\Map_k'(X[E],X)\to \Map_{\Fun^{\rL}_k(\Dqc(X),\Dqc(X))}(\id,(-)\otimes E[1]).
\]
We call this the\textbf{ Atiyah map of $E$.}
\end{enumerate}
\end{definition}

\subsection{Functoriality of Atiyah class}

Next, we check functoriality of universal Atiyah transformations. Let $X,Y$ be derived stacks with cotangent complexes over $k$.  Consider $f:X\to Y$ over $k$, complexes $L_X\in\Dqc(X),L_Y\in\Dqc(Y)$, morphisms $d_X\in\Map'_k(X[L_X],X),d_Y\in\Map'_k(Y[L_Y],Y)$. 
We lastly consider a morphism $\phi:L_Y\to f_*L_X$.
We demand that these data (at the level of morphism of complexes) satisfy $f_*(d_X)\circ [f]\simeq\phi\circ d_Y$ where the composition morphisms are
\[
f_*(d_X)\circ [f]:\bL_{Y/k}\to f_* \bL_{X/k}\xrightarrow{f_*(d_X)} f_*L_X,\quad 
\phi\circ d_Y:\bL_{Y/k}\xrightarrow{d_Y}  L_Y\to f_* L_X.
\]
Geometrically, this condition gives us a commutative diagram compatible with the closed immersions into the trivial square zero thickenings.
\[
\xymatrix{
X[L_X]\ar[r]\ar[d]^-{f_\phi}&X[\bL_{X/k}]\ar[d]^-{f_{[f]}}\ar[r]^-{\delta_X}&X\ar[d]^-f\\
Y[L_{Y}]\ar[r]& Y[\bL_{Y/k}]\ar[r]^-{\delta_Y}&Y
}
\] 
Here, the horizontal compositions are $d_X$ and $d_Y$, and $\delta_X,\delta_Y$ arise from $\id_{\bL_{X/k}},\id_{\bL_{Y/k}}$. It follows that when $\pi^X:X[L_{X}]\to X,\pi^Y:Y[L_{Y}]\to Y$ are the projections and $\iota^X:X\to X[L_{X}],\iota^Y:Y\to Y[L_{Y}]$ are the obvious morphisms, $f_\phi$ and $f$ make obvious square diagrams commute with respect to the pairs of morphisms $(\pi^X,\pi^Y),(\iota^X,\iota^Y),(d_X,d_Y)$. 
Choose $F\in\Dqc(Y)$.
\begin{lemma}\label{lemma-atiyah-functoriality}
Consider a morphism $\psi:A\to (f_\phi)_*B$ in $\Dqc(Y[L_Y])$ for some
objects $A\in\Dqc(Y[L_Y]), B\in\Dqc(X[L_X])$. Then, we obtain a natural morphism associated with $\psi:$
\[
f^*\pi^Y_*(d_Y^* F\otimes A)\to \pi^X_*(d_X^* f^*F\otimes B).
\] 
\end{lemma}
\begin{proof}
We first construct the composition
    \[
    f^*\pi^Y_*(d_Y^* F\otimes A)\to f^*\pi^Y_*(d_Y^* F\otimes (f_\phi)_*B)
    \to  f^*\pi^Y_*(f_\phi)_*(f_\phi^*d_Y^* F\otimes B)
    \simeq f^*f_*\pi^X_* (d_X^* f^*F\otimes B)
    \] where the first arises from $\psi$, the second arises from the projection formula morphism, and the last arises from the commutative squares mentioned right before the lemma. By applying the counit of adjunction, the last object admits a map $f^*f_*\pi^X_* (d_X^* f^*F\otimes B)\to \pi^X_* (d_X^* f^*F\otimes B)$. Composing this map yields the desired map. 
\end{proof}
By using the commutative square diagrams mentioned before the lemma, we will construct a morphism of fiber sequences in $\Dqc(Y[L_{Y}])$
\[
\xymatrix{
\iota^Y_*L_Y\ar[r]\ar[d]^-{\iota^Y_*\phi}&\cO_{Y[L_Y]}\ar[d]\ar[r]&\iota_*^Y\cO_Y\ar[d]\\
(f_\phi)_*\iota^X_*L_X\ar[r]&(f_\phi)_*\cO_{X[L_X]}\ar[r]&(f_\phi)_*\iota_*^X\cO_X
}
\]
To do this, we factor the map $\phi:L_Y\to f_* L_X$ through the unit of adjunction as $L_Y\to f_*f^*L_Y\to f_*L_X$, which geometrically corresponds to the factoring $X[L_X]\to X[f^*L_Y]\to Y[L_Y].$
By looking at one of these two maps at a time, we only need to construct the diagram in the following cases:
\begin{enumerate}
    \item $X=Y$ and $f=\id$. 
    The commutative diagram exists in this case because we can reduce to the nonconnective affine case as in Lemma \ref{lemma-nonconnective-sqz-fiber-sequence}.
    \item $L_X=f^*L_Y$, and $\phi$ is the unit of adjunction $\phi:L_Y\to f_*f^*L_Y.$ In this case, $X[f^*L_Y]$ is the base change $X\times_YY[L_Y]$, and the fiber sequence on $X[f^*L_Y]$ is the pullback of the fiber sequence on $Y[L_Y]$ via Lemma \ref{lemma-nonconnective-base-change-projection-formula}.  In particular, we have the desired diagram by taking the units of adjunction vertically.  
\end{enumerate}
In the diagram we just constructed, applying appropriate functors to the fiber sequences yields the Atiyah classes of $F$ and $f^*F$, but invoking Lemma \ref{lemma-atiyah-functoriality} on the vertical maps, we can relate the resulting fiber sequences as follows:
\[
\xymatrix{
f^*\pi^Y_*(d_Y^* F\otimes\iota^Y_*L_Y)
\ar[r]\ar[d]&
f^*\pi^Y_*(d_Y^* F\otimes\cO_{Y[L_Y]})
\ar[d]\ar[r]&
f^*\pi^Y_*(d_Y^* F\otimes\iota_*^Y\cO_Y)
\ar[d]\\
\pi^X_*(d_X^* f^*F\otimes\iota^X_*L_X)
\ar[r]&
\pi^X_*(d_X^* f^*F\otimes\cO_{X[L_X]})
\ar[r]&
\pi^X_*(d_X^* f^*F\otimes\iota_*^X\cO_X)
}
\]
The rotation of the top fiber sequence yields
$f^*\At_{d_Y}(F):f^*F\to f^*(F\otimes L_Y[1])$, and 
the rotation of the bottom fiber sequence yields
$\At_{d_X}(f^*F):f^*F\to f^*F\otimes L_X[1]$.
The rotated diagram contains the following commutative square where $^A\phi:f^*L_Y\to L_X$ corresponds to $\phi$ by adjunction.  
\[
\xymatrix{
f^*F\ar[rr]^-{f^*\At_{d_Y}(F)}\ar[d]^-{\id}&&f^*F\otimes f^*L_Y[1]\ar[d]^-{f^*F\otimes {}^A\phi[1]}\\
f^*F\ar[rr]^-{\At_{d_X}(f^*F)}&& f^*F\otimes L_X[1]
}
\]
We conclude:
\begin{proposition}
In the above situation, the Atiyah class is functorial in the sense that $\At_{d_X}(f^*F)\simeq (f^*F\otimes {}^A\phi[1]) \circ f^*\At_{d_Y}(F)$.    
\end{proposition}
\begin{remark}\label{remark-special-cases-functorial-atiyah}
    We observe that the above situation applies to the following special cases.
    \begin{enumerate}
        \item When $X=Y$ and $f=\id$, the above proposition describes the effect of changing the ideals of square zero thickenings between $L_X$ and $L_Y$. For each morphism $\phi:L_Y\to L_X$
        satisfying $d_X\simeq \phi\circ d_Y$
        , we obtain a commutative diagram:
        \[
        \xymatrix{
\Map(\bL_{X/k},L_Y)\ar[r]^-{\bAt_{d_Y}}\ar[d]^-{\phi\circ(-)}&\Map_{\Fun^{\rL}(\Dqc(X),\Dqc(X))}(\id,(-)\otimes L_Y[1])\ar[d]^-{(\id\otimes\phi[1])\circ (-)}\\
\Map(\bL_{X/k},L_X)\ar[r]^-{\bAt_{d_X}}&\Map_{\Fun^{\rL}(\Dqc(X),\Dqc(X))}(\id,(-)\otimes L_X[1])
        }
        \]
        As a further special case of this, we can let $L_Y=\bL_{X/k},\phi:\bL_{X/k}\to L_X$ and observe that $\At_{\phi}(-)=(\id\otimes\phi[1])\circ \At_{\id_{\bL_{X/k}}}(-)$. In this sense, $\At_{\id_{\bL_{X/k}}}(-)$ is universal among universal Atiyah transformations on $X$. This shows in particular that the Atiyah map is $k$-linear.
        \item When $d_X=\delta_X, d_Y=\delta_Y$ and $\phi=[f]^A:\bL_{Y/k}\to f_*\bL_{X/k}$, this describes the functoriality of Atiyah classes in the conventional sense.
    \end{enumerate}
\end{remark}
\subsection{Atiyah class on the product}
Let $X,Y$ be scalloped derived stacks over $k$, so that the Kunneth type formula for $\Dqc$ holds over $k$ (Lemma \ref{lemma-kunneth-perfect}). 
We write the projections as  $q_X:X\times_k Y\to X, q_Y:X\times_k Y\to Y$. Write $Z=X\times_k Y$ Take $E\in\Dqc(X),F\in\Dqc(Y).$ 
\begin{proposition}
Let $E\in\Dqc(X),F\in\Dqc(Y)$.
Under the decomposition $\bL_{X\times_k Y/k}\simeq\bL_{X/k}\boxplus\bL_{Y/k}:=q_X^*\bL_{X/k}\oplus q_Y^*\bL_{Y/k}$, the morphism
    \[
    \At_{\delta_{X\times_k Y}}(E\boxtimes F):E\boxtimes F\to (E\boxtimes F)\otimes\bL_{X\times_k Y}[1],
    \]
    is identified with the morphism $\At_{\delta_X}(E)\boxtimes\id_F+ \id_E\boxtimes\At_{\delta_Y}(F)$:
    \[
     E\boxtimes F\to 
    (E\boxtimes F)\otimes q_X^*\bL_{X/k}[1]\oplus (E\boxtimes F)\otimes q_Y^*\bL_{Y/k}[1]\simeq (E\boxtimes F)\otimes\bL_{X\times_k Y/k}[1]
    \]
    functorially in $E$ and $F$. 
\end{proposition}
\begin{proof}
Let $p_X:\bL_{X\times_k Y/k}\to q_X^*\bL_{X/k}$ be the projection. We apply the functoriality framework from the previous subsubsection where $f=\id_{X\times_k Y}$, $L_X=q_X^*\bL_{X/k}, L_Y=\bL_{X\times_k Y},\phi=p_X,$ and ``$d_X$'' and ``$d_Y$'' are respectively given by $\phi$ and $\delta_{X\times Y}.$ As a result, we see
\[
\xymatrix{
E\boxtimes F\ar[rr]^-{\At_{\delta_{X\times_k Y}}(E\boxtimes F)}\ar[drr]_-{\At_{p_X}(E\boxtimes F)}&&E\boxtimes F\otimes \bL_{X\times_k Y}[1]\ar[d]^-{\id\otimes p_X}\\
&& E\boxtimes F\otimes q_X^*\bL_{X/k}[1]
}
\]
This shows that the projection of $\At_{\delta_{X\times_k Y}}(E\boxtimes F)$ to the factor involving $q_X^*\bL_{X/k}$ is given by $\At_{p_X}(E\boxtimes F)$, so we wish to show
$\At_{p_X}(E\boxtimes F)=\At_{\delta_X}(E)\boxtimes\id_F$. 

We observe $(X\times_k Y)[q_X^*\bL_{X/k}]\simeq X[\bL_{X/k}]\times_k Y$ from the functoriality of relative spectrum, and under this isomorphism, the morphism $p_X:(X\times_k Y)[q_X^*\bL_{X/k}]\to X\times_k Y$ is identified with $\delta_X\times\id_Y:X[\bL_{X/k}]\times_k Y\to X\times_k
 Y,$ the projection $\widetilde{\pi}:(X\times_k Y)[q_X^*\bL_{X/k}]\to X\times_k Y$ is identified with $\pi^{X}\times \id_{Y}:X[\bL_{X/k}]\times_k Y\to X\times_k Y$, and the closed immersions from $X\times_k Y$ are identified similarly. Then, by base change (Lemma \ref{lemma-nonconnective-base-change-projection-formula}), we see that $\At_{p_X}(E\boxtimes F)\simeq\At_{\delta_X}(E)\boxtimes\id_F$. By symmetry we have $\At_{p_Y}(E\boxtimes F)\simeq\id_E\boxtimes\At_{\delta_Y}(F)$ when $p_Y:\bL_{X\times_k Y}\to q_Y^*\bL_{Y/k}$ is the other projection. Adding them together, we obtain the desired formula. 
\end{proof}
\begin{corollary}\label{corollary-atiyah-product}
    There is an equivalence of natural transformations
    \[
    {\At_{\delta_{X\times Y}}}\simeq 
    j_X\circ 
    ({\At_{\delta_{X}}}\boxtimes \id_{\Dqc(Y)})+
    j_Y\circ
    (\id_{\Dqc(X)}\boxtimes{\At_{\delta_{ Y}}}):\id\to (-)\otimes \bL_{X\times_k Y/k}[1]
    \] where the functors are $\id,(-)\otimes \bL_{X\times_k Y/k}[1]\in \Fun_k^\rL( \Dqc(X\times_k Y),\Dqc(X\times_k Y))$, and $j_X:(-)\otimes q_X^*\bL_{X/k}[1]\to (-)\otimes \bL_{X\times_k Y/k}[1]$ and $j_Y:(-)\otimes q_Y^*\bL_{Y/k}[1]\to (-)\otimes \bL_{X\times_k Y/k}[1]$ are induced from the inclusions of direct summands.
    Here, box tensor products make sense under the equivalence $\Dqc(X\times_k Y)\simeq \Dqc(X)\otimes_k \Dqc(Y)$. 
\end{corollary}
\begin{proof}
As in \cite[Lemma 2.7]{NCHPD}, it suffices to prove this equivalence after restricting to the full subcategory spanned by box tensor products. Thus, it reduces to the proposition above.
\end{proof}
\subsection{Atiyah class of tensor products}
\begin{proposition}
    Consider $E,F\in\Dqc(X)$, where $X$ is scalloped. Then, we have an equivalence
    \[
    \At_{\delta_X}(E\otimes F)\simeq \At_{\delta_X}(E)\otimes \id_F+\id_E\otimes  \At_{\delta_X}(F):E\otimes F\to E\otimes F\otimes \bL_{X/k}[1]
    \] whose construction is functorial in $E$ and $F$.
\end{proposition}
\begin{proof} 
Applying the functoriality to the morphism $\Delta:X\to X\times_k X$ with respect to the ``universal derivations'' gives us
\[
\xymatrix{
\Delta^*(E\boxtimes F)\ar[rr]^-{\Delta^*\At_{\delta_{X\times_k X}}(E\boxtimes F)}\ar[drr]_-{\At_{\delta_X}(\Delta^*(E\boxtimes F))}&&\Delta^*(E\boxtimes F)\otimes\Delta^*\bL_{X\times_k X/k}[1]\ar[d]^-{\id\otimes \Sigma[1]}\\
&&\Delta^* (E\boxtimes F)\otimes \bL_{X/k}[1]
}
\]
where $\Sigma:\Delta^*\bL_{X\times_k X/k}\simeq \bL_{X/k}\oplus \bL_{X/k}\to \bL_{X/k}$ is the summation map. $\At_{\delta_{X\times_k X}}(E\boxtimes F)$ is computed
in Corollary \ref{corollary-atiyah-product}, and $\At_{\delta_X}(\Delta^*(E\boxtimes F))=\At_{\delta_X}(E\otimes F)$, so we obtain the formula. 
\end{proof}

\subsection{Atiyah class of the universal line bundle on $\bB\bG_m$}
In this subsection, we compute the Atiyah class of $\cO(1)$ on $\bB\bG_m$ (over $k$). This is done in \cite[Example 4.9/4.10]{Kuhn} based on a different construction of Atiyah classes on stacks. Here, we reinterpret Kuhn's computation and present it in a self-contained manner. 
In this section, we write $T$ as the coordinate variable of $\bG_m$, so that $\bG_m=\Spec k[T^{\pm1}].$
In this subsection, for shortness, given a {scalloped derived stack $X$ over $k$, complexes $E\in\Dqc(X),F\in\Dqc(X)$,}
and a morphism $\alpha\in\Map_{X//k}(X[E],X)\simeq \Map(\bL_{X/k},E)$, we write $P(F):=\pi_*\alpha^*(F)$, where $\pi:X[E]\to X$ is the projection. This notation will make $\alpha$ implicit since it should be clear from context in this subsection. We also always work with Atiyah classes with respect to the identity morphisms of cotangent complexes in this subsection, so we will simply write it as $\At(-)$ by omitting subscripts.
\[
\xymatrix{
\bG_m\ar[r]^t\ar[d]^s&\Spec k\ar[d]^x\\
\Spec k\ar[r]^x&\bB\bG_m
}\tag{$\blacksquare$}\label{diagram-blacksquare}
\]
Using classical descent, we can describe the line bundle $\cO(1)$ on $\bB\bG_m$ by the following data: $(x^*\cO(1),x^*\cO(1),\phi: s^*(x^*\cO(1))\simeq t^*(x^*\cO(1)))$. Trivializing $x^*\cO(1)\simeq k$, more explicitly, we can write it as $(k,k,\times T: \cO_{\bG_m}\simeq\cO_{\bG_m})$.
\begin{lemma}\label{lemma-isom-cotanget-complex-point-BGm}
    There is an isomorphism $k\simeq \bL_{\Spec k/\bB\bG_m}$ whose pullback along $t$ is identified with $\cO_{\bG_m}\to \bL_{s}=\bL_{\bG_m/\Spec k}$ sending $1\mapsto d\log T$ via base change $t^*\bL_{\Spec k/\bB\bG_m}\simeq \bL_{s}$. 
\end{lemma}
\begin{proof}
It's well known that $\bL_{\Spec k/\bB\bG_m}$ is a sheaf, so by classical descent, we can describe it by considering the descent data with respect to a cover of $\Spec k$. We use the cover $\bG_m\to\Spec k$ (and the resulting self-product) which we interpret in the following way.
Due to multiple copies of $\Spec k$ we will see below, we write $\Spec k_1\to \bB\bG_m$ for the original morphism which we use to compute $\bL_{\Spec k/\bB\bG_m}$.

    We base change the morphism $\Spec k_1\to \bB\bG_m$ by base-changing $\bB\bG_m$ to $\Spec k$ and $\bG_m$ along the square diagram as in (\ref{diagram-blacksquare}) (but we use labels $s',t'$ to distinguish here).
    \begin{equation}
\begin{tikzcd}[row sep={3.5em,between origins}, column sep={5em,between origins}]
    & \bG_m \times_k \bG_m \arrow[rr, "m"] \arrow[dd, "p_2" near end] \arrow[dl, "p_1"'] & & 
    \bG_m \arrow[dd] \arrow[dl] \\
    \bG_m \arrow[rr,"s" near end,crossing over] \arrow[dd, "t"] & & 
    \Spec k_1 \arrow[dd] \\
    & \bG_m \arrow[rr, "t'" near end] \arrow[dl, "s'"'] & & 
    \Spec k_3 \arrow[dl] \\
    \Spec k_2 \arrow[rr] & & 
    \bB\bG_m 
\end{tikzcd}
\tag{$\square$}\label{diagram-base-change-cube}
\end{equation}
In this diagram, we use the identification \[\Spec k_1\times_{\bB\bG_m}\Spec k_2\times_{\bB\bG_m}\Spec k_3\]
\[\simeq (\Spec k_1\times_{\bB\bG_m}\Spec k_2)\times_{k}(\Spec k_2\times_{\bB\bG_m}\Spec k_3)
\simeq 
\bG_{m,x}\times_k\bG_{m,y},\]
where $x,y$ are the coordinate variables of $\bG_m$. We see in our identifications that the three projection morphisms of form
$\Spec k \times_{\bB\bG_m}\Spec k\times_{\bB\bG_m}\Spec k\to 
\Spec k \times_{\bB\bG_m}\Spec k$ are identified with $p_1,p_2,m:\bG_m\times_k\bG_m\to\bG_m$ as in the diagram.
We describe $\bL_{\Spec k_1/\bG_m}$ by the descent datum $(\bL_{\bG_m/\Spec k_2},\bL_{\bG_m/\Spec k_3}, \phi:p_1^*\bL_{\bG_m/\Spec k_2}\simeq m^*\bL_{\bG_m/\Spec k_3})$.
Here, the gluing datum $\phi$ is the composition $p_1^*\bL_{\bG_m/\Spec k_2}\simeq \bL_{p_2}\simeq m^*\bL_{\bG_m/\Spec k_3}.$ The first map sends the generator $d\log T=\frac{dT}{T}$ to $\frac{dx}{x}=d\log x.$ The inverse of the second map sends the generator $d\log T$ to $\frac{d(xy)}{xy}=\frac{ydx}{xy}=d\log x$, where we used $dy=0$ in $\bL_{p_2}.$ This means that the glueing map $\phi$ sends $d\log T$ to $d\log T$. Thus, this gluing datum is isomorphic to the gluing datum of $k$, which is $(\cO_{\bG_m},\cO_{\bG_m},\id:p_1^*\cO_{\bG_m}\simeq m^*\cO_{\bG_m})$. The isomorphism of gluing data is given by the maps $\cO_{\bG_m}\to \bL_{\bG_m/\Spec k_i};1\mapsto d\log T$ for $i=2,3.$ This isomorphism $k\simeq \bL_{\Spec k/\bB\bG_m}$ has the desired property by construction.
\end{proof}
\begin{construction}\label{construction-psi-trivialize-LBGm}
Using the lemma above we will fix a trivialization $\psi:\cO_{\bB\bG_m}[-1]\simeq \bL_{\bB\bG_m/\Spec k}$. It's well-known that these are canonically isomorphic, but for computational purposes, we construct $\psi$ in the following way.
We consider the chain of canonical isomorphisms
\[
\hspace*{-1.5em}
\Hom(\cO_{\bB\bG_m}[-1], \bL_{\bB\bG_m/\Spec k})\xrightarrow[\simeq]{x^*} \Hom(k[-1],x^*\bL_{\bB\bG_m/\Spec k})\simeq \Hom(k[-1],\bL_{\Spec k/\bB\bG_m}[-1]),
\] where the first map is an isomorphism as we know that $ \bL_{\bB\bG_m/\Spec k}\simeq\cO_{\bB\bG_m}[-1]$. The second isomorphism arises from the transitivity fiber sequence of cotangent complexes
\[
x^*\bL_{\bB\bG_m/\Spec k}\to  \bL_{\Spec k/\Spec k}\to \bL_{\Spec k/\bB\bG_m}.
\]
We define $\psi\in\Hom(\cO_{\bB\bG_m}[-1], \bL_{\bB\bG_m/\Spec k})$ to be the isomorphism corresponding to the $[-1]$-shift of the trivialization of Lemma \ref{lemma-isom-cotanget-complex-point-BGm} via the chain of isomorphisms above. 
\end{construction}

 By the functionality argument of the Atiyah class, we obtain the following commutative diagram, where all the horizontal sequences are fiber sequences: 
 \[
 \begin{tikzcd}[row sep={3.5em,between origins}, column sep={5em,between origins}]
    &  t^*(x^*\cO(1)\otimes\bL_{k/k}) \arrow[rr] \arrow[dd] & & 
    t^*P(x^*\cO(1)) \arrow[rr] \arrow[dd] & & 
   t^*x^*\cO(1) \arrow[dd] 
   \\
    (\cO(1)\otimes\bL_{\bB\bG_m/k})|_{\bG_m} \arrow[rr] \arrow[ur] \arrow[dd] & & 
    P(\mathcal{O}(1))|_{\bG_m} \arrow[rr] \arrow[ur] \arrow[dd] & & 
    \mathcal{O}(1)|_{\bG_m} \arrow[ur] \arrow[dd] 
    \\
    &\cO(1)|_{\bG_m}\otimes\bL_{\bG_m/k}  \arrow[rr] & & 
    P(\cO(1)|_{\bG_m}) \arrow[rr] & & 
    \cO(1)|_{\bG_m} 
    \\
    s^*(x^*\cO(1)\otimes\bL_{k/k}) \arrow[rr] \arrow[ur] & & 
    s^*P(x^*\cO(1)) \arrow[rr] \arrow[ur] & & 
   s^*x^*\cO(1) \arrow[ur]
\end{tikzcd}
\tag{$\square$}\label{diagram-cube}
    \]
    The ``connecting maps'' of these pullbacks are certain pullbacks of Atiyah classes, so we use the following diagram to represent the above commutative diagram of fiber sequences more compactly:
    \[
    \xymatrix{
\At(\cO(1))|_{\bG_m}\ar[r]\ar[d]&t^*\At(x^*\cO(1))\ar[d]\\
 s^*\At(x^*\cO(1))   \ar[r]& \At(\cO(1)|_{\bG_m})
    }
    \]
    For computational purposes, we will use the viewpoint of descent data at the beginning of subsection, so we rewrite this to make the gluing data coming from $\phi$ explicit. 
    \[
    \xymatrix{
s^*x^*\At(\cO(1))\ar[r]_-{\simeq}^-{\phi(\At(\cO(1)))}\ar[d]&s^*x^*\At(\cO(1))\ar[d]\\
s^*\At(x^*\cO(1))\ar[d]^{}&t^*\At(x^*\cO(1))\ar[d]^{}\\
\At(s^*x^*\cO(1))\ar[r]^-{\At(\phi)}_-{\simeq}&\At(t^*x^*\cO(1))
    }
    \]
We temporarily focus on the bottom two rows of this diagram by trivializing the descent data of $\cO(1)$ as $(k,k,\times T:\cO_{\bG_m}\simeq\cO_{\bG_m})$. Writing down the bottom four fiber sequences gives us:
 \[
 \begin{tikzcd}[row sep={3.5em,between origins}, column sep={5em,between origins}]
    & 0 \arrow[rr] \arrow[dd] & & 
    t^*P(k) \arrow[rr] \arrow[dd] & & 
   t^*k \arrow[dd] 
   \\
    0\arrow[rr]  \arrow[dd] & & 
    s^*P(k) \arrow[rr] \arrow[dd] & & 
    s^*k  \arrow[dd] 
    \\
    &t^*k\otimes  \bL_{\bG_m/k}  \arrow[rr] & & 
    t^*k\oplus (t^*k\otimes\bL_{\bG_m/k}) \arrow[rr] & & 
    t^*k
    \\
s^*k\otimes  \bL_{\bG_m/k} \arrow[ur, "{\phi\otimes\id}"] \arrow[rr] & & 
    s^*k\oplus (s^*k\otimes\bL_{\bG_m/k}) \arrow[rr] \arrow[ur, "{\pi_*\delta^*(\phi)}"'] & & 
   s^*k \arrow[ur, "\phi"]
\end{tikzcd}
    \]
    Here, the cotangent complex $\bL_{k/k}$ is trivial, so we have $s^*P(k)\simeq s^*k$ and $t^*P(k)\simeq t^*k$, from which we obtain splittings of the bottom two fiber sequences. We wrote the middle terms of the bottom two fiber sequences as direct sums using these splittings.
    The morphism $\phi:s^*k\to t^*k$ is given by multiplication of $T$, so the morphism $P(s^*k)\to P(t^*k)$ is constructed by applying the functor $\pi_*\delta^*(-)$ to the morphism $\times T:\cO_{\bG_m}\to\cO_{\bG_m}$. Recall that at the ring level, the morphism $\delta:\bG_m[\bL_{\bG_m/k}]\to\bG_m$ sends $T\in \Gamma(\bG_m,\cO_{\bG_m})$ to $(T,dT)\in \Gamma(\bG_m,\cO_{\bG_m})\oplus \Gamma(\bG_m,\Omega_{\bG_m/k}^1)$. Thus, the morphism \[\pi_*\delta^*(\times T):s^*k\oplus(s^*k\otimes\bL_{\bG_m/k})=\cO_{\bG_m}\oplus \bL_{\bG_m/k}\to \cO_{\bG_m}\oplus \bL_{\bG_m/k}=t^*k\oplus(t^*k\otimes\bL_{\bG_m/k})\] is the multiplication map by $(T,dT)$ with respect to the square-zero structure. It is then given by the matrix
    \[
    \begin{pmatrix}
        T&0\\
        dT&T
    \end{pmatrix}.
    \]
    If we trivialize all the instances of $s^*k$ by the standard generator $1\in t^*k$ via the isomorphism $\phi$, this matrix can be rewritten as 
    \[
    \begin{pmatrix}
        1&0\\
        d\log T&1
    \end{pmatrix}.
    \]
    (For example, $1\in \Gamma(\bG_m,t^*k)$ corresponds to the element $T^{-1}\in\Gamma(\bG_m,s^*k)$ via $\phi$, and this element is sent by $\pi_*\delta^*(\phi)$ to $T^{-1}(T+dT)=1+d\log T$ by the original matrix. This describes the left column of the new matrix.)
    We use this computation to rewrite (\ref{diagram-cube}) with respect to the trivialization by $t^*k$ as well as the splitting $t^*k\oplus (t^*k\otimes\bL_{\bG_m/k})$:
    \[
 \begin{tikzcd}[row sep={3.5em,between origins}, column sep={5em,between origins}]
    &0\arrow[rr] \arrow[dd] & & 
    \cO_{\bG_m} \arrow[rr,equal] \arrow[dd, "{(\id,0)}" near start] & & 
   \cO_{\bG_m}\arrow[dd,equal] 
   \\
    (\cO(1)\otimes\bL_{\bB\bG_m/k})|_{\bG_m} \arrow[rr] \arrow[ur] \arrow[dd] & & 
    P(\mathcal{O}(1))|_{\bG_m} \arrow[rr] \arrow[ur] \arrow[dd] & & 
    \mathcal{O}(1)|_{\bG_m} \arrow[ur] \arrow[dd] 
    \\
    &\bL_{\bG_m} \arrow[rr] & & 
    \cO_{\bG_m}\oplus\bL_{\bG_m/k} \arrow[rr] & & 
    \cO_{\bG_m}
    \\
    0\arrow[rr] \arrow[ur] & & 
    \cO_{\bG_m} \arrow[rr,equal] \arrow[ur,"{(\id,\times d\log T)}"'] & & 
   \cO_{\bG_m} \arrow[ur, equal]
\end{tikzcd}
    \]
Here,
the morphism $(\id,0)$ is given by this formula since the we used the splitting $t^*k\oplus (t^*k\otimes\bL_{\bG_m/k})$ instead of $s^*k\oplus (s^*k\otimes\bL_{\bG_m/k})$.
The morphism $(\id,\times d\log T)$ arises as the left column of the matrix we just computed. We extend the diagram by adding two more cubes where the new vertical squares that are not parallel to the paper are cartesian.
    \[
\begin{tikzcd}[row sep={3em,between origins}, column sep={4em,between origins}]
    & & 0 \arrow[rrr] \arrow[ddd] & & & 0 \arrow[rrr] \arrow[ddd] & & & 0 \arrow[ddd] \\
    & 0 \arrow[rrr] \arrow[ur] \arrow[ddd] & & & \cO_{\bG_m} \arrow[rrr,equal] \arrow[ur] \arrow[ddd, "{(\id,0)}" near start] & & & \cO_{\bG_m} \arrow[ur] \arrow[ddd,equal] \\
 (\cO(1)\otimes\bL_{\bB\bG_m/k})|_{\bG_m} \arrow[rrr] \arrow[ur] \arrow[ddd] & & & P(\mathcal{O}(1))|_{\bG_m} \arrow[rrr] \arrow[ur] \arrow[ddd] & & & \mathcal{O}(1)|_{\bG_m} \arrow[ur] \arrow[ddd] \\
    & & \bL_{\bG_m/\Spec k} \arrow[rrr, equal] & & &  \bL_{\bG_m/\Spec k} \arrow[rrr] & & & 0 \\
    & \bL_{\bG_m/\Spec k} \arrow[rrr] \arrow[ur, equal] & & & \cO_{\bG_m}\oplus\bL_{\bG_m/\Spec k} \arrow[rrr] \arrow[ur,"{(0,\id)}"'] & & & \cO_{\bG_m} \arrow[ur] \\
    0 \arrow[rrr] \arrow[ur] & & & \cO_{\bG_m} \arrow[rrr,equal] \arrow[ur,"{(\id,\times d\log T)}"'] & & & \cO_{\bG_m} \arrow[ur, equal]
\end{tikzcd}
\tag{$\sharp$}
\label{diagram-extended-4-cubes}
\]
By the following lemma, in fact all the vertical squares that are not parallel to the paper are cartesian:
\begin{lemma}[{see also \cite[Lemma 2.10/2.11]{Kuhn}}]\label{lemma-cartesian-squares-atiyah-cubes}
    The vertical squares in the diagram (\ref{diagram-cube}) that are not parallel to the paper are cartesian.
\end{lemma}
\begin{proof}
Clearly, the right square is cartesian.
In general, when we have a pullback diagram of Arin derived stacks over the base $S$,
\[\xymatrix{
A\ar[r]\ar[d]&B\ar[d]\\
C\ar[r]&D
}
\] the natural morphisms give us a pullback square of cotangent complexes.
\[\xymatrix{
\bL_{D/S}|_A\ar[r]\ar[d]&\bL_{C/S}|_A\ar[d]\\
\bL_{B/S}|_A\ar[r]&\bL_{A/S}
}
\]
It follows that the left square is cartesian. 
Then, the middle square must be cartesian as well, since we can reconstruct it as the fibers in the rotated fiber sequences.
\end{proof}
\begin{proposition}[{\cite[Example 4.9/4.10]{Kuhn}}]\label{proposition-atiyah-O(1)}
The Atiyah class $\At(\cO(1)):\cO(1)\to \cO(1)\otimes\bL_{\bB\bG_m/\Spec k}[1]$ of the universal line bundle $\cO(1)$ on $\bB\bG_m$ (over $\Spec k$) can be written as 
\[
\At(\cO(1))\simeq\id_{\cO(1)}\otimes\psi[1].
\]
   
\end{proposition}
\begin{proof}
By Lemma \ref{lemma-cartesian-squares-atiyah-cubes}, the right and left faces of the largest cube in diagram (\ref{diagram-extended-4-cubes}) are cartesian squares. The same is true for the big parallel square in the middle. The top vertices of three squares are 0, so they form fiber sequences. Extracting these fiber sequences, we obtain the following commutative diagram whose rows and columns are fiber sequences.
\[
    \xymatrix{
    (\cO(1)\otimes \bL_{\bB\bG_m/\Spec k})|_{\bG_m}\ar[r] \ar[d] & 
    P(\cO(1))|_{\bG_m} \ar[r]^\alpha \ar[d]^\beta & 
    \cO(1)|_{\bG_m} \ar@{=}[d] & 
     \\
    0 \ar[r] \ar[d] & 
    \cO_{\bG_m} \ar@{=}[r] \ar[d]^-{\times d\log T} & 
    \cO_{\bG_m} \ar[d] & \\
    \bL_{\bG_m/\Spec k} \ar@{=}[r] & 
    \bL_{\bG_m/\Spec k}\ar[r] & 
    0 & 
}
    \]
    We emphasize that the ``connecting map'' of the top horizontal fiber sequence is $t^*x^*\At(\cO(1))$.
    {This is the cofiber of $\alpha$, or equivalently $\beta$. Hence, \[t^*x^*\At(\cO(1)):\cO(1)|_{\bG_m}\to(\cO(1)\otimes \bL_{\bB\bG_m/\Spec k})|_{\bG_m}[1]\]
    can be identified with $\times d\log T:\cO_{\bG_m}\to \bL_{\bG_m/\Spec k}$. By an easy homological algebra, this identification is bridged by $\bL_{\bG_m/\Spec k}\simeq (\cO(1)\otimes \bL_{\bB\bG_m/\Spec k})|_{\bG_m}[1]$ which appears as the connecting map of the left vertical fiber sequence.
Up to the identification $\phi\otimes\id: s^*x^*\cO(1)\otimes\bL_{\bG_m/\Spec k}\simeq t^*x^*\cO(1)\otimes\bL_{\bG_m/\Spec k}$, the vertical fiber sequence $A_{1,1}\to A_{1,2}\to A_{1,3}$ in our situation originates from the sequence
\[
t^*x^*(\cO(1)\otimes\bL_{\bB\bG_m/\Spec k})\to 
t^*(x^*\cO(1)\otimes\bL_{\Spec k/\Spec k})
\to 
t^*x^*\cO(1)\otimes\bL_{\bG_m/\Spec k}.
\] 
This is $t^*(-)$ of the transitivity fiber sequence of cotangent complexes for the composition $\Spec k\to\bB\bG_m\to\Spec k$ up to trivializing $t^*x^*\cO(1)$ and base change $t^*\bL_{\Spec k/\bB\bG_m}\simeq \bL_{\bG_m/\Spec k}$. 
On the other hand, by construction of $\psi$, $x^*\psi$ is the $d\log T$ map up to the identification $x^*\bL_{\bB\bG_m/\Spec k}\simeq\bL_{\Spec k/\bB\bG_m}[-1]$ coming from the transitivity fiber sequence. Hence, we see $t^*x^*\At(\cO(1))\simeq t^*x^*(\id_{\cO(1)}\otimes \psi[1]).$
Since the map $t^*x^*:\Hom(\cO(1),\cO(1)\otimes\bL_{\bB\bG_m/\Spec k}[1])\to \Hom(t^*x^*\cO(1),t^*x^*(\cO(1)\otimes\bL_{\bB\bG_m/\Spec k}[1]))$ is injective, we see $\id_{\cO(1)}\otimes\psi[1]\simeq\At(\cO(1))$.
}
\end{proof}
\section{Atiyah class for twisted sheaves}\label{section-atiyah-class-for-twisted-sheaves}
\subsection{Twisted Atiyah class in general}
Assume $f:X\to Y$ is a principal $\bB\bG_m$-bundle over a scalloped derived stack $Y$ over a field $k$.   Choose a not necessarily connective object $E\in\Dqc^0(X/Y)$. 
 
We additionally choose a point $\alpha\in\Map'_k(X[E],X):=\Map_k(X[E],X)\times_{\Map_k(X,X)}\{\id_X\}\simeq\Map(\bL_{X/k},E)$ (Remark \ref{remark-cotangent-complex-nonconnective-universality}). 
By Definition \ref{definition-atiyah-transformation}, we obtain the (untwisted) universal $\alpha$-Atiyah transformation $\At_\alpha(-):\id\to (-)\otimes E[1]$ in $\Fun_k^{\rL}(\Dqc(X),\Dqc(X))$. 
By restricting $\Dqc(X)$ to  the full subcategory $\Dqc^1(X/Y)$, we obtain a morphism in $\Fun_k^{\rL}(\Dqc^1(X/Y),\Dqc(X))$.
$E$ is 0-twisted, so we have $\id,(-)\otimes E[1]\in\Fun_k^\rL(\Dqc^1(X/Y),\Dqc^1(X/Y))$ by Proposition \ref{proposition-symmetric-monoidal-Dqc}.
$\Dqc^1(X/Y)\subset \Dqc(X)$ is a full subcategory, so we have an equivalence
\[\Map_{\Fun^\rL_k(\Dqc^1(X/Y),\Dqc^1(X/Y))}(\id,(-)\otimes E[1])\simeq \Map_{\Fun^\rL_k(\Dqc^1(X/Y),\Dqc(X))}(\id,(-)\otimes E[1]),\]
and we
obtain a morphism $\widetilde{\At}_\alpha(-):\id\to(-)\otimes E[1]$ in $\Fun^\rL_k(\Dqc^1(X/Y),\Dqc^1(X/Y)).$ 
\begin{definition}\label{definition-twisted-atiyah-transformation}

\leavevmode
\begin{enumerate}
    \item The morphism $\widetilde{\At}_\alpha(-):\id\to(-)\otimes E[1]$ in $\Fun^\rL_k(\Dqc^1(X/Y),\Dqc^1(X/Y))$ is called the \textbf{twisted universal $\alpha$-Atiyah transformation of $X/Y$.}
    \item For a particular object $F\in\Dqc^1(X/Y),$ the evaluation of the above transformation at $F$ is denoted by $\widetilde{\At}_{\alpha}(F):F\to F\otimes E[1]$, and called the \textbf{twisted $\alpha$-Atiyah class of $F$ on $X$.} We observe  $\widetilde{\At}_{\alpha}(F)={\At}_{\alpha}(F)$.
    \item  
    By Examples \ref{example-twisted-FMT-via-dualizability} and \ref{example-twisted-FMT-tensor}, $\widetilde{\At}_{\alpha}(-)$ corresponds to a morphism of twisted Fourier-Mukai kernels, and we denote it by $\widetilde{\At}_{\alpha}^{\FM}\in\RHom_{X\times_kX}(\iota^X_*\cP_X,\iota^X_*(\cP_X\otimes p_2^*E[1]))$. This is called the\textbf{ twisted universal $\alpha$-Atiyah class of $X.$} 
    \item If we fix $E$ and vary $\alpha$, via Remark \ref{remark-cotangent-complex-nonconnective-universality} the above construction gives us a map
\[
\widetilde{\bAt}_E:\Map_k'(X[E],X)\to \Map_{\Fun^{\rL}_k(\Dqc^1(X),\Dqc^1(X))}(\id,(-)\otimes E[1]).
\]
We call this the\textbf{ twisted Atiyah map of $E$.}
\end{enumerate}
\end{definition}
{
\begin{remark}
    When $f:X\to Y$ is as above, $\bL_{X/Y}\simeq\cO_{X}[-1]$, so that by the transitivity fiber sequence for the morphism $f$ over $k$, $\bL_{X/k}$ is $0$-twisted. In particular, we may choose $E=\bL_{X/k}$ in the above definition. 
\end{remark}
}
\begin{proposition}\label{proposition-twisted-atiyah-functorial}
The twisted Atiyah map 
$\widetilde{\bAt}_E$ is natural in $E\in\Dqc^0(X)$, so that up to the identifications in Remark \ref{remark-cotangent-complex-nonconnective-universality}, the following diagram commutes for any morphism $\phi:E\to F$ in $\Dqc^0(X)$:
    \[
    \xymatrix{
\Map_k(\bL_{X/k},E)\ar[r]^-{\widetilde{\bAt_E}}\ar[d]^-{\phi\circ(-)}&  \Map_{\Fun_k^{\rL}(\Dqc^1(X),\Dqc^1(X))}(\id,(-)\otimes E[1])\ar[d]^-{(\id\otimes\phi[1])\circ(-)}\\
\Map_k(\bL_{X/k},F)\ar[r]^-{\widetilde{\bAt_F}}&  \Map_{\Fun_k^{\rL}(\Dqc^1(X),\Dqc^1(X))}(\id,(-)\otimes F[1])
    }
    \]
    In particular, for any $\phi:\bL_{X/k}\to E$,  $\widetilde{\At}_{\phi}(-)\simeq (\id\otimes\phi[1])\circ \widetilde{\At}_{\id_{\bL_{X/k}}}(-)$, and $\widetilde{\bAt}_E$ is $k$-linear.
\end{proposition}
\begin{proof}
In the diagram below, the top square commutes by the untwisted analogue Remark \ref{remark-special-cases-functorial-atiyah}, and the remaining squares clearly commute. 
\[
\xymatrix{
\Map(\bL_{X/k},E)\ar[d]^-{\bAt_E}\ar[rr]^-{\phi\circ(-)}
&&
\Map(\bL_{X/k},F)\ar[d]^-{\bAt_F}
\\
  \Map_{\Fun^{\rL}_k(\Dqc(X),\Dqc(X))}(\id,(-)\otimes E[1])\ar[rr]^-{(\id\otimes\phi[1])\circ(-)}\ar[d]
  &&
   \Map_{\Fun^{\rL}_k(\Dqc(X),\Dqc(X))}(\id,(-)\otimes F[1])\ar[d]
  \\
 \Map_{\Fun^{\rL}_k(\Dqc^1(X),\Dqc(X))}(\id,(-)\otimes E[1])\ar[rr]^-{(\id\otimes\phi[1])\circ(-)}
 &&
 \Map_{\Fun^{\rL}_k(\Dqc^1(X),\Dqc(X))}(\id,(-)\otimes F[1])
 \\
\Map_{\Fun^{\rL}_k(\Dqc^1(X),\Dqc^1(X))}(\id,(-)\otimes E[1])\ar[rr]^-{(\id\otimes\phi[1])\circ(-)}\ar[u]^-{\simeq}
&&
 \Map_{\Fun^{\rL}_k(\Dqc^1(X),\Dqc^1(X))}(\id,(-)\otimes F[1])\ar[u]^-{\simeq}
}
\]
The two vertical compositions are $\widetilde{\bAt_E}$ and $\widetilde{\bAt_F}$, so we are done.
\end{proof}
\begin{definition}\label{definition-spectrally-enriched-twisted-atiyah-map}
    By Proposition \ref{proposition-twisted-atiyah-functorial}, the twisted Atiyah map upgrades from mapping spaces to mapping spectra. 
    The resulting map in terms of the twisted Fourier-Mukai kernels is denoted by 
    \[
    \widetilde{\bAt}_E^{\FM}:\RHom(\bL_{X/k},E)\to \RHom_{\Dqc^{-1,1}(X\times_kX/Y\times_kY)}(\iota^X_*\cP_X,\iota^X_*(\cP_X\otimes p_2^* E[1])).
    \]
    We call it the \textbf{spectrally enriched twisted Atiyah map of $E.$}
    We can describe $ \widetilde{\bAt}_E^{\FM}$ as the composition:
\[
\hspace*{-1em}
\RHom(\bL_{X/k},E)\xrightarrow[]{\iota^X_*(\cP_X\otimes p_2^*(-)[1])} \RHom_{\Dqc^{-1,1}(X\times_kX/Y\times_kY)}(\iota_*^X(\cP_X\otimes p_2^*\bL_{X/k}[1]),\iota_*^X(\cP_X\otimes p_2^*E[1]))
\]
\[\xrightarrow[]{(-)\circ \widetilde{\At}_{\id_{\bL_{X/k}}}^{\FM}} \RHom_{\Dqc^{-1,1}(X\times_kX/Y\times_kY)}(\iota^X_*\cP_X,\iota^X_*(\cP_X\otimes p_2^* E[1]))
\]
Here, the morphism $\widetilde{\At}_{\id_{\bL_{X/k}}}^{\FM}$ is of the form $\widetilde{\At}_{\id_{\bL_{X/k}}}^{\FM}:\iota^X_*\cP_X\to \iota^X_*(\cP_X\otimes p_2^* \bL_{X/k}[1])$.
\end{definition}

\subsection{Twisted Atiyah classes for connective complexes}\label{subsection-twisted-atiyah-connective}
In this subsection, we specialize to the case where $E$ is connective. In this case, the construction of the twisted Atiyah class can be viewed as a strict analogue of the untwisted version, where we replace functors by twisted versions of functors in each step. For the rest of this section, we will use several results from \S\ref{subsection-deform-BGm-bundles-sqz-thickening}, but these results will not depend on the current section.

We continue to assume $f:X\to Y$ is a principal $\bB\bG_m$-bundle over a scalloped derived stack $Y$ over a field $k$.  
We choose $E\in\Dqccn^0(X)$ to be connective and {$0$-twisted}, and write
\[
i:X\to X[E], \quad u:X[{E[1]}]\to X, \quad\pi:X[E]\to X
\] for the obvious inclusions and projections.
Additionally, we consider a point \[\alpha\in\Map'_k(X[E],X):=\Map_k(X[E],X)\times_{\Map_k(X,X)}\{\id_X\},\]
where the map $\Map_k(X[E],X)\to \Map_k(X,X)$ is the composition with $i$. By virtue of Proposition \ref{proposition-sqz-BGm-bundle-pullback}, there are morphisms $f^E:X[{E}]\to Y[{f_*E}]$, $f^{E[1]}:X[{E[1]}]\to Y[{f_*E[1]}]$, as well as $i^Y:Y\to Y[{f_*E}],u^Y:Y[{f_*E[1]}]\to Y, \pi^Y:Y[{f_*E}]\to Y,\alpha^Y:Y[{f_*E}]\to Y$ which form pullback squares of principal $\bB\bG_m$-bundle morphisms with respect to $f$, $f^E,$ and $f^{E[1]}$. As a result, by Example \ref{example-twisted-FMT-push-pull}, we can talk about pushforward and pullback of twisted sheaves along $i,u,\pi,\alpha$.

On $X[E]$, we have a fiber sequence of complexes 
$
i_*E\to\cO_{X[E]}\to i_* \cO_X.
$ 
Since these complexes are $0$-twisted, tensoring with objects in the above fiber sequence determine functors in $\Fun_k^\rL(\Dqc^1(X[E]),\Dqc^1(X[E]))$. 
This defines a fiber sequence pointwise, so
this is a fiber sequence of functors.
Composing these functors by $\alpha^*:\Dqc^1(X/Y)\to\Dqc
^1(X[E]/Y[f_*E])$ and $\pi_*:\Dqc^1(X[E]/Y[f_*E])\to\Dqc^1(X/Y)$, we have morphisms of functors from $\Dqc^1(X/Y)$ to itself:
\[
\pi_*(\alpha^*(-)\otimes i_*E)\to \pi_*(\alpha^*(-)\otimes  \cO_{X[E]}) \to \pi_*(\alpha^*(-)\otimes i_*\cO_X)
\] This is still a fiber sequence by the pointwise argument.
As in Lemma \ref{lemma-nonconnective-base-change-projection-formula}, projection formula applies, and this simplifies to 
\[
\pi_*i_*(i^*\alpha^*(-)\otimes E)\to \pi_*\alpha^*(-) \to \pi_*i_*i^*\alpha^*(-).
\]
Using $\pi\circ i\simeq \id_X$ and the identification $\id\simeq\alpha\circ  i$, we simplify this as
\[
(-)\otimes E\to \pi_*\alpha^* \to \id.
\]
By construction of Definition \ref{definition-twisted-atiyah-transformation}, the rotation of this fiber sequence is precisely the twisted universal $\alpha$-Atiyah transformation $\widetilde{\At}_\alpha(-):\id\to(-)\otimes E[1]$.

\subsection{Noncommutative twisted Atiyah map}
We keep the setup of the previous subsection. 
In this subsection, we construct a categorical version of the twisted Atiyah map. In the untwisted case (for derived schemes), this is essentially what Moore calls ``the derived Atiyah transformation'' in \cite[\S3.1]{Moore2}.
\begin{lemma}\label{lemma-pushout-of-sqz-to-pullback-Dqc}
Let $X$ be a qcqs 1-Artin derived stack.
Let $E\in\Dqccn(X)$. Then, pullbacks along inclusions and projections as in the diagram below
\[
\xymatrix{
X[E[1]]\ar[r]^-u\ar[d]^-u&X\ar[d]^-i\\
X\ar[r]^-i& X[E]
}
\] induce a fully faithful functor 
\[
\Phi:\Dqc(X[E])\to \Dqc(X)\times_{\Dqc(X[E[1]])}\Dqc(X)
\]
and this restricts to an equivalence
\[
\Psi:\Dperf(X[E])\simeq \Dperf(X)\times_{\Dperf(X[E[1]])}\Dperf(X).
\]
\end{lemma}
\begin{proof}
In the case where $X$ is a derived affine, this is known as in \cite[Theorem 16.2.0.1]{SAG} or \cite[Lemma 3.3]{Moore2}. 
In general, we write $X=\colim_i U_i$ for  derived affines $U_i$ (by taking the Cech nerve of an atlas).  

We write $E_i:=E|_{U_i}$. By taking the limit of the functors for all $U_i$, we obtain a fully faithful functor
\[
\lim_i\Dqc(U_i[E_i])\to \lim_i\left(\Dqc(U_i)\times_{ \Dqc(U_i[E_i[1]])} \Dqc(U_i)\right).
\]
Since $X$ is the colimit of $U_i$, we have $\Dqc(X)\simeq \lim_i\Dqc(U_i)$ (see \cite[Corollary 3.22]{blowup-normal-bundle-nonconnective-DAG} for example). By the compatibility of base change and relative spectrum, as well as that colimits are universal in $\infty$-topoi, we also see $\Dqc(X[E])\simeq \lim_i\Dqc(U_i[E_i]),$ $\Dqc(X[E[1]])\simeq \lim_i\Dqc(U_i[E_i[1]])$.
By commuting limits, we see that the above functor is precisely $\Phi.$ 
The case of perfect complexes is similar once we show that the equivalence $\Dqc(X)\simeq\lim_i\Dqc(U_i)$ restricts to $\Dperf(X)\simeq\lim_i\Dperf(U_i)$ (and similarly for other derived stacks), but this is clear since perfectness is equivalent to dualizability (\cite[Proposition 3.6]{bzfn}).
\end{proof}
\begin{corollary}[twisted version of {\cite[Lemma 3.4]{Moore2}}]\label{corollary-Map-Dqc1-vs-Dqc1-omega}
\leavevmode
\begin{enumerate}
    \item  We have a pullback functor $i^*:\Dqc^1(X[E]/Y[f_*E])\to\Dqc^1(X/Y)$ which reflects compact objects. 
    \item  Let $\Cat_k^{\perf}$ denote the $\infty$-category of small idempotent-complete $k$-linear stable $\infty$-categories with exact functors as in \cite[\S2.3.1]{Moore2}. Then, \[\Map'_k(\Dqc^1(X/Y),\Dqc^1(X[E]/Y[f_*E]))\simeq \Map_{\Cat_k^{\perf}}'(\Dqc^1(X/Y)^{\omega}, \Dqc^1(X[E]/Y[f_*E])^{\omega}),\] where 
    the left hand side is defined as the fiber of the map
\[i^*\circ (-):\Map_k(\Dqc^1(X/Y),\Dqc^1(X[E]/Y[f_*E]))\to {\Map_k(\Dqc^1(X/Y),\Dqc^1(X/Y))}\] above $\id_{\Dqc^1(X/Y)}$, and the right hand side is defined similarly.
\end{enumerate}
\end{corollary}
\begin{proof}
\leavevmode
\begin{enumerate}
    \item 
    We have commutative diagrams
\[
\xymatrix{
X[E[1]]\ar[r]^-u\ar[d]^-u&X\ar[d]^-i&Y[f_*E[1]]\ar[r]^-{u^Y}\ar[d]^-{u^Y}&Y\ar[d]^-{i^Y}\\
X\ar[r]^-i& X[E]&Y\ar[r]^-{i^Y}& Y[f_*E]
}
\] of derived stacks with respect to the inclusions and projections.
Moreover, using the maps $f$, $f^{E[1]}:X[E[1]]\to Y[f_*E[1]],f^{E}:X[E]\to Y[f_*E],$ we can combine these squares to form a commutative cubic diagram via Proposition \ref{proposition-sqz-BGm-bundle-pullback}. The newly added faces of the cube are cartesian by the proposition. 
Then, the functors $u^*,i^*$ respect the decomposition of $\Dqc$ into homogeneous components. 
    By Lemma \ref{lemma-scalloped-object-compact-perfect-dualizable} and Lemma \ref{lemma-trivial-sqz-thickening-scalloped}, perfect objects and compact objects coincide on relevant derived stacks, so $i^*:\Dqc(X[E])\to\Dqc
    (X)$ reflects compact objects by Lemma \ref{lemma-pushout-of-sqz-to-pullback-Dqc}. 
    By Proposition \ref{proposition-basic-property-Dqc1-presentable-stable-compact}(2), the same holds for $i^*:\Dqc^1(X[E])\to \Dqc^1(X)$.
    \item By part (1), any section of $i^*$ respects compact objects, so \[\Map'_k(\Dqc^1(X/Y),\Dqc^1(X[E]/Y[f_*E]))\simeq \Map'_{\Pr^{\rL,\omega}_k}(\Dqc^1(X/Y),\Dqc^1(X[E]/Y[f_*E])).\]
    By
    Proposition \ref{proposition-basic-property-Dqc1-presentable-stable-compact}, $\Ind(\Dqc^1(X/Y)^\omega)\simeq \Dqc^1(X/Y)$ and $\Ind(\Dqc^1(X[E]/Y[f_*E])^\omega)\simeq \Dqc^1(X[E]/Y[f_*E])$,  
    so recalling the identification $\Ind:\Cat^{\perf}_k\simeq \Pr^{\rL,\omega}_k$ as in \cite[\S2.3.1]{Moore2},  we obtain
    \[\Map'_{\Pr^{\rL,\omega}_k}(\Dqc^1(X/Y),\Dqc^1(X[E]/Y[f_*E]))\simeq \Map'_{\Cat_k^{\perf}}(\Dqc^1(X/Y)^\omega,\Dqc^1(X[E]/Y[f_*E])^\omega).\]
\end{enumerate}
\end{proof}
\begin{proposition}[Twisted analogue of {\cite[Proposition 3.5]{Moore2}}]\label{proposition-twisted-categorical-atiyah}
    There is a canonical equivalence 
    \[
    \Map'_k(\Dqc^1(X/Y),\Dqc^1(X[E]/Y[f_*E]))\simeq 
    \Map_{\Fun_k^{\rL}(\Dqc^1(X/Y),\Dqc^1(X/Y))}(\id,(-)\otimes E[1]).
    \]
    functorial in $E$.
\end{proposition}
\begin{proof}
By Lemma \ref{lemma-pushout-of-sqz-to-pullback-Dqc}, we have a commutative diagram
\[
\xymatrix{
\Dqc(X[E])\ar[r]^-{i^*}\ar[d]^-{i^*}&\Dqc(X)\ar[d]^-{u^*}\\
\Dqc(X)\ar[r]^-{u^*}&\Dqc(X[E[1]])
}
\] inducing a fully faithful functor to pullback (and an equivalence on $\Dperf$).
In view of the cubical diagram described in the proof of Corollary \ref{corollary-Map-Dqc1-vs-Dqc1-omega}(1), the functors $i^*$ and $u^*$ respect the pieces of the decompositions of $\Dqc$ coming from the principal bundle structures.
By Proposition \ref{proposition-basic-property-Dqc1-presentable-stable-compact} and Remark \ref{remark-Perf1-binda-porta}, we obtain a pullback diagram:
\[
\xymatrix{
\Dqc^1(X[E]/Y[f_*E])^\omega\ar[r]^-{i^*}\ar[d]^-{i^*}&\Dqc^1(X/Y)^\omega\ar[d]^-{u^*}\\
\Dqc^1(X/Y)^\omega\ar[r]^-{u^*}&\Dqc^1(X[E[1]]/Y[f_*E[1]])^\omega
}
\]
We apply $\Map_{\Cat_k^{\perf}}(\Dqc^1(X/Y)^\omega,-)$ and take the fiber over $\id$, resulting in the following cartesian square:
\[
\begin{adjustbox}{max width=\linewidth}
\xymatrix{
\Map_{\Cat_k^{\perf}}'(\Dqc^1(X/Y)^\omega,\Dqc^1(X[E]/Y[f_*E])^\omega)\ar[r]\ar[d]&
\Map_{\Cat_k^{\perf}}'(\Dqc^1(X/Y)^\omega,\Dqc^1(X/Y)^\omega)\ar[d]\\
\Map_{\Cat_k^{\perf}}'(\Dqc^1(X/Y)^\omega,\Dqc^1(X/Y)^\omega)\ar[r]&
\Map_{\Cat_k^{\perf}}'(\Dqc^1(X/Y)^\omega,\Dqc^1(X[E[1]]/Y[f_*E[1]])^\omega)
}
\end{adjustbox}
\]
By Corollary \ref{corollary-Map-Dqc1-vs-Dqc1-omega}, we can rewrite this as 
\[
\xymatrix{
\Map_{k}'(\Dqc^1(X/Y),\Dqc^1(X[E]/Y[f_*E]))\ar[r]\ar[d]&
\Map_{k}'(\Dqc^1(X/Y),\Dqc^1(X/Y))\ar[d]\\
\Map_{k}'(\Dqc^1(X/Y),\Dqc^1(X/Y))\ar[r]&
\Map_{k}'(\Dqc^1(X/Y),\Dqc^1(X[E[1]]/Y[f_*E[1]]))
}
\]
The rest of the proof is similar to \cite[Proposition 3.5]{Moore2}, where we remark that $u_*u^*\simeq \id\oplus (-)\otimes E[1]$ by projection formula (which holds for example by Lemma \ref{lemma-nonconnective-base-change-projection-formula}). 
{Let us briefly sketch the intuition. By the cartesian diagram above, a point $\alpha^*:\Dqc^1(X/Y)\to \Dqc^1(X[E]/Y[f_*E])$ at top left corner is roughly equivalent to the datum $(i^*\alpha^*, i^*\alpha^*,\beth(\alpha^*):u^*i^*\alpha^*\simeq u^*i^*\alpha^*)$ (omitting the trivializations by $\id_{\Dqc^1(X/Y)}$). This last gluing datum is equivalent to $\overline{\beth(\alpha^*)}:i^*\alpha^*\to u_*u^*i^*\alpha^*$ by adjunction. If we use the given identifications, the triple can be written as $(\id,\id,\id\to u_*u^*)$ where we remember that last morphism's composition with the projection $u_*u^*\simeq \id\oplus(-)\otimes E[1]\to\id$ is the identity. Then, the data of such triples are the same as morphisms $\id\to (-)\otimes E[1]$. In short, our equivalence sends $\alpha^*$ from the left hand side to the composition
\[
\id\simeq i^*\alpha^*\xrightarrow{\overline{\beth(\alpha^*)}}
u_*u^*i^*\alpha^*\simeq u_*u^*\to (-)\otimes E[1]
\] on the right hand side.}
\end{proof}
\begin{definition}\label{definition-noncommutative-twisted-atiyah-map}
Let $f:X\to Y,E\in\Dqccn^0(X)$ be as before.
 The equivalence from Proposition \ref{proposition-twisted-categorical-atiyah} is denoted as
\[
\widetilde{\bAt}_E^{\nc}: \Map'_k(\Dqc^1(X/Y),\Dqc^1(X[E]/Y[f_*E]))\simeq \Map_{\Fun_k^{\rL}(\Dqc^1(X),\Dqc^1(X))}(\id,(-)\otimes E[1]),
\] and 
we call this the \textbf{noncommutative twisted Atiyah map of $E$.}
\end{definition}
\begin{remark}
    In \cite[Proposition 3.5]{Moore2}, Moore constructed an untwisted version of this equivalence (for qcqs derived $k$-schemes), and used it to construct what he calls the ``derived Atiyah transformation'' in \cite[Construction 3.7]{Moore2}.
\end{remark}
\subsection{Compatibility of commutative and noncommutative twisted Atiyah maps}
We continue to assume that $f:X\to Y$ is a principal $\bB\bG_m$-bundle of scalloped derived stacks over a field $k$, and choose $E\in\Dqccn^0(X)$. 
Given any $\alpha\in \Map'_k(X[E],X)$, by Proposition \ref{proposition-sqz-BGm-bundle-pullback}
this defines pullback squares of principal $\bB\bG_m$-bundles:
\[
\xymatrix{
X\ar[r]^-i\ar[d]^-f&X[E]\ar[r]^-\alpha\ar[d]& X\ar[d]^-f\\
Y\ar[r]&Y[f_*E]\ar[r]& Y
}
\]
This gives us the pullback functor $\alpha^*:\Dqc^1(X/Y)\to\Dqc^1(X[E]/Y[f_*E])$. Also, the trivialization $\alpha\circ i\simeq\id$ gives us a trivialization $i^*\circ \alpha^*\simeq\id$ of the functors between twisted derived categories.
Thus, we have a morphism
\[
(-)^*:\Map'_k(X[E],X)\to \Map'_k(\Dqc^1(X/Y),\Dqc^1(X[E]/Y[f_*E]));\alpha\mapsto\alpha^*
\]
With respect to this morphism, we claim the following:
\begin{proposition}[{twisted analogue of \cite[Proposition 3.12]{Moore2}}]\label{proposition-commutative-vs-noncommutative-atiyah-maps}
In the situation above, the following diagram  commutes.
   \[
\xymatrix{
\Map'_k(X[E],X)\ar[dr]^-{\widetilde{\bAt}_E}\ar[d]^-{(-)^*}\\
\Map'_k(\Dqc^1(X/Y),\Dqc^1(X[E]/Y[f_*E]))\ar[r]^-{\widetilde{\bAt}_E^{\nc}}
&\Map_{\Fun_k^{\rL}(\Dqc^1(X),\Dqc^1(X))}(\id,(-)\otimes E[1]),
}
\] 
\end{proposition}
The rest of the subsection is devoted to the proof of this proposition.
\begin{lemma}[twisted version of {\cite[Lemma 3.14]{Moore2}}]\label{lemma-cartesian-square-pi_*-alpha^*}
For any \[\alpha^*\in\Map'_k(\Dqc^1(X/Y),\Dqc^1(X[E]/Y[f_*E])),\]
we have a pullback diagram in $\Fun_k^{\rL}(\Dqc^1(X/Y),\Dqc^1(X/Y))$ whose construction is functorial in $\alpha^*$:
\[
\xymatrix{
\pi_*\alpha^*\ar[r]\ar[d]& \pi_* i_*i^*\alpha^*\ar[d]^-{\overline{\beth(\alpha^*)}}\\
\pi_*i_*i^*\alpha^*\ar[r]^-\eta& \pi_*i_*u_*u^*i^*\alpha^*
}
\]
or more simply
\[
\xymatrix{
\pi_*\alpha^*\ar[r]\ar[d]& \id\ar[d]^-{\overline{\beth(\alpha^*)}}\\
\id\ar[r]^-\eta& u_*u^*
}
\]
\end{lemma}
\begin{proof}
    {The proof is similar to \cite[Lemma 3.14]{Moore2}, where Moore constructs the second pullback diagram, but in order to clarify functoriality, we sketch the construction while making the instances of identifications $i^*\circ \alpha^*\simeq \id,i^*\circ \pi^*\simeq \id$ explicit. Since pullbacks in the functor category are computed pointwise, for each $F\in\Dqc^1(X/Y)$, we want to establish $\pi_*\alpha^*F\simeq F\times_{u_*u^*F}F$ functorially in $F$. By adjunction, it suffices to show \[\Map_{\Dqc^1(X/Y)}(G,F\times_{u_*u^*F}F)\simeq \Map_{\Dqc^1(X[E]/Y[f_*E])}(\pi^*G,\alpha^*F)\] for any $G\in\Dqc^1(X/Y)$. As in the proof of Proposition \ref{proposition-twisted-categorical-atiyah}, the pullback functor \[\rho:\Dqc^1(X[E]/Y[f_*E])\to\Dqc^1(X/Y)\times_{\Dqc^1(X[E[1]]/Y[f_*E[1]])}\Dqc^1(X/Y)\]
    is fully faithful, so 
    \[
    \Map_{\Dqc^1(X[E]/Y[f_*E])}(\pi^*G,\alpha^*F)\simeq 
    \Map_{\Dqc^1(X/Y)\times_{\Dqc^1(X[E[1]]/Y[f_*E[1]])}\Dqc^1(X/Y)}(\rho(\pi^*G),\rho(\alpha^*F)).
    \]
    $\rho(\pi^*G)$ is the triple $(i^*\pi^*G,i^*\pi^*G,\id:u^*i^*\pi^*G\simeq u^*i^*\pi^*G)$ or more simply $(G,G,\id:u^*G\simeq u^*G)$, and $\rho(\alpha^*F)$ can be written as the triple $(i^*\alpha^*F,i^*\alpha^*F,\beth(\alpha^*)_F:u^*i^*\alpha^*F\simeq u^*i^*\alpha^* F)$, using the symbol from the proof of Proposition \ref{proposition-twisted-categorical-atiyah}. If we use $i^*\alpha^*\simeq\id$, this simplifies to $(F,F,\beth(\alpha^*)_F: u^*F\simeq u^*F)$.
    Without simplificiation, a morphism $\rho(\pi^*G)\to\rho(\alpha^*F)$ consists of two morphisms $f_i:i^*\pi^*G\to i^*\alpha^*F$ ($i=1,2$) along with the following commuting square:
    \[
    \xymatrix{
    u^*i^*\pi^*G\ar@{=}[r]\ar[d]^-{u^*f_1}&u^*i^*\pi^*G\ar[d]^-{u^*f_2}\\
    u^*i^*\alpha^*F\ar[r]^-{\beth(\alpha^*)_F}&u^*i^*\alpha^*F
    }
    \]
    Hence, the mapping space is equivalent to
    \[
    \Map_{\Dqc^1(X/Y)}(i^*\pi^*G,i^*\alpha^*F)
    \times_{\Map_{\Dqc^1(X[E[1]]/Y[f_*E[1]])}    (u^*i^*\pi^*G,u^*i^*\alpha^*F)}
    \Map_{\Dqc^1(X/Y)}(i^*\pi^*G,i^*\alpha^*F)
    \]
    \[
    \simeq 
     \Map_{\Dqc^1(X/Y)}(G,\pi_*i_*i^*\alpha^*F)\times_{\Map_{\Dqc^1(X/Y)}(G,\pi_*i_*u_*u^*i^*\alpha^*F)} \Map_{\Dqc^1(X/Y)}(G,\pi_*i_*i^*\alpha^*F)
    \]
    \[
    \simeq \Map_{\Dqc^1(X/Y)}(G,\pi_*i_*i^*\alpha^*F\times_{\pi_*i_*u_*u^*i^*\alpha^*F} \pi_*i_*i^*\alpha^*F).
    \]
    This shows 
    \[
    \pi_*\alpha^* F\simeq \pi_*i_*i^*\alpha^*F\times_{\pi_*i_*u_*u^*i^*\alpha^*F} \pi_*i_*i^*\alpha^*F,
    \] or more simply
    $\pi_*\alpha^*F\simeq F\times_{u_*u^*F}F$.
    }
\end{proof}
\begin{corollary}\label{corollary-cartesian-square-pi_*-alpha^*}
    For any $\alpha^*\in\Map'_k(\Dqc^1(X/Y),\Dqc^1(X[E]/Y[f_*E]))$, we have a pullback diagram in $\Fun^{\rL}(\Dqc^1(X/Y),\Dqc^1(X/Y))$ whose construction is functorial in $\alpha^*$:
\[
\xymatrix{
\pi_*\alpha^*\ar[r]\ar[d]& \pi_* i_*i^*\alpha^*\simeq \id\ar[d]\\
0\ar[r]& (-)\otimes E[1]
}
\]
where the top horizontal map is induced by the unit of adjunction $\id\to i_*i^*$.
\end{corollary}
\begin{proof}
By projection formula (Lemma \ref{lemma-nonconnective-base-change-projection-formula}), we have an equivalence $u_*u^*\simeq \id\oplus (-)\otimes E[1]$
in $\Fun_k^{\rL}(\Dqc^1(X/Y),\Dqc^1(X/Y))$.
By applying $\pi_*\circ  i_*\circ (-)\circ i^*\circ \alpha^*$ to this decomposition, we obtain a cartesian square
\[
\xymatrix{
\pi_*i_*i^*\alpha^*\ar[r]^-\eta\ar[d]& \pi_*i_*u_*u^*i^*\alpha^*\ar[d]\\
0\ar[r]&\pi_*i_* (i^* \alpha^*(-)\otimes E[1])
}
\]
By pasting this cartesian square with the one from Lemma \ref{lemma-cartesian-square-pi_*-alpha^*}, we obtain the desired cartesian square. By construction, this is functorial in $\alpha^*.$
\end{proof}
\begin{proof}[{Proof of Proposition \ref{proposition-commutative-vs-noncommutative-atiyah-maps}}]
    Consider a morphism $\alpha\in\Map'_k(X[E],X)$. 
    Since $E$ is connective,  $\widetilde{\At}_{E}(\alpha)$ can be described as in \S\ref{subsection-twisted-atiyah-connective}. In particular, $\widetilde{\At}_{E}(\alpha)$ is the cofiber of the morphism $\pi_*\alpha^*(-) \to \pi_*i_*i^*\alpha^*(-)$ in $\Fun^{\rL}(\Dqc^1(X/Y),\Dqc^1(X/Y))$. This is precisely the right vertical morphism in the diagram of  Corollary \ref{corollary-cartesian-square-pi_*-alpha^*}. This morphism is constructed as 
    \[
    \id\simeq \pi_*i_*i^*\alpha^*\xrightarrow{\overline{\beth(\alpha^*)}}\pi_*i_*u_*u^*i^*\alpha^*\simeq u_*u^*\xrightarrow[]{}(-)\otimes E[1],
    \] which is precisely $\widetilde{\bAt}^{\nc}_E(\alpha^*)$ by the last line of the proof of Proposition \ref{proposition-twisted-categorical-atiyah}.

\end{proof}
\section{Twisted Atiyah transformations on trivial principal $\bB\bG_m$-bundles}\label{section-twisted-atiyah-transformation-trivial-principal-BGm-bundles}
Let $Y$ be {smooth proper scheme} over $k$, and let $X=Y\times_k\bB\bG_m$, where $\bB\bG_m$ is understood to be defined over $k$. Let $f:X\to Y$ be the projection, and $X$ has the obvious principal $\bB\bG_m$-bundle structure. $g:X\to\bB\bG_m$ is the other projection.
We compute the twisted universal $\id_{\bL_{X/k}}$-Atiyah class. This transformation takes the form
\[
\widetilde{\At}_{\id_{\bL_{X/k}}}(-):
\id_{\Dqc^1(X/Y)}\to(-)\otimes\bL_{X/k}[1],
\] which lives in $\Fun^{\rL}_k(\Dqc^1(X/Y),\Dqc^1(X/Y))$.
Since the principal $\bB\bG_m$-bundle is trivial in our case, we can identify $\Dqc^1(X/Y)$ with $\Dqc^0(X/Y)\simeq\Dqc(Y).$ Explicitly, this equivalence is given by the functor $(-)\otimes\cO(1)|_X$. Up to this trivialization, the twisted Atiyah transformation looks like 
\[
\left(\widetilde{\At}_{\id_{\bL_{X/k}}}((-)\otimes\cO(1)|_X)\right)\otimes\cO(-1)|_X:\id_{\Dqc^0(X/Y)}\to(-)\otimes\bL_{X/k}[1]
\] which lives in $\Fun^{\rL}(\Dqc^0(X/Y),\Dqc^0(X/Y))$.
For any $E\in\Dqc(Y)$, 
evaluating the above transformation at $f^*E$ yields
\[
\At_{\id_{\bL_{X/k}}}(E\boxtimes\cO(1))\otimes\cO(-1)|_X:f^*E\to f^*E\otimes\bL_{X/k}[1].
\]
By Corollary \ref{corollary-atiyah-product}, this can be written as the sum 
$f^*\At_{\id_{\bL_{Y/k}}}(E)+\id_{E}\boxtimes (\At_{\id_{\bL_{\bB\bG_m/k}}}(\cO(1))\otimes\cO(-1)):E\to E\otimes\bL_{X/k}[1]$. By Proposition \ref{proposition-atiyah-O(1)}, 
this is equal to 
$f^*\At_{\id_{\bL_{Y/k}}}(E)+\id_{f^*E}\otimes g^*\psi[1]$ when we view $f^*\bL_{Y/k},g^*\bL_{\bB\bG_m/k}\subset\bL_{X/k}$ as direct summands, where $\psi:\cO_{\bB\bG_m}[-1]\simeq \bL_{\bB\bG_m/k}$ is defined in Construction \ref{construction-psi-trivialize-LBGm}.
This computation yields the following:
\begin{proposition}\label{proposition-twisted-atiyah-trivial-gerbe-formula}
    Under the identification 
    \[\hspace*{-1em}
    \Map_{\Fun_k^{\rL}(\Dqc^1(X/Y),\Dqc^1(X/Y))}(\id,(-)\otimes\bL_{X/k}[1])
    \simeq 
    \Map_{\Fun^{\rL}_k(\Dqc(Y),\Dqc(Y))}(\id,f_*(f^*(-)\otimes\bL_{X/k}[1]))
    \]
    \[
    \simeq   \Map_{\Fun_k^{\rL}(\Dqc(Y),\Dqc(Y))}(\id,(-)\otimes f_*\bL_{X/k}[1]))
    \]
    induced by the equivalence $f^*(-)\otimes\cO(1):\Dqc(Y)\simeq \Dqc^1(X/Y)$, $\widetilde{\At}_{\id_{\bL_{X/k}}}(-)$ corresponds to
    \[
 \Psi   :=
(j_1\circ \At_{\id_{\bL_{Y/k}}}(-))+(j_2\circ ((-)\otimes f_*g^*\psi[1])):\id_{\Dqc(Y)}\to (-)\otimes f_*\bL_{X/k}[1],
    \] where $j_1: (-)\otimes \bL_{Y/k}[1]\simeq  (-)\otimes f_*f^*\bL_{Y/k}[1]\to (-)\otimes f_*\bL_{X/k}[1]$ and
    $j_2:(-)\otimes f_*g^*\bL_{\bB\bG_m/k}[1]\to(-)\otimes f_*\bL_{X/k}[1]$ are induced by inclusions of direct summands.
\end{proposition}

For any $E\in\Dqc(Y)$, the twisted  Atiyah map of $f^*E$ can be trivialized as:
\[
\xymatrix{
\Map(\bL_{X/k},f^*E)\ar[d]^-{\widetilde{\At}_{f^*E}}\\
\Map_{\Fun_k^{\rL}(\Dqc^1(X/Y),\Dqc^1(X/Y))}(\id,(-)\otimes f^*E[1])\ar[d]^-{\simeq}\ar[r]^-{\simeq}&
\Map_{X\times_k X}(\iota^X_*\cP_X,\iota^X_*(\cP_X\otimes p_2^*f^*E[1]))
\ar[d]^-{\simeq}\\
\Map_{\Fun_k^{\rL}(\Dqc(Y),\Dqc(Y))}(\id,(-)\otimes E[1])\ar[r]^-{\simeq}&\Map_{Y\times_k Y}(\Delta_*\cO_Y,\Delta_*E[1]))
}\tag{$\triangle$}\label{diagram-trivial-gerbe-Atiyah-map}
\]
The horizontal equivalences come from the twisted and untwisted Fourier-Mukai theory as in Lemma \ref{lemma-dualizability-functor-Fourier-Mukai}. 
The top left vertical map sends $\alpha:\bL_{X/k}\to f^*E$ to 
\[ 
\id_{\Dqc^1(X/Y)}\xrightarrow[]{\widetilde{\At}_{\id_{\bL_{X/k}}}}(-)\otimes\bL_{X/k}[1]\xrightarrow[]{(-)\otimes\alpha[1]}(-)\otimes f^*E[1]
\]
Then its target at the bottom left corner can be written as 
\[
\id_{\Dqc(Y)}\xrightarrow[]{\Psi}(-)\otimes f_*\bL_{X/k}[1]\xrightarrow[]{f_*((f^*(-)\otimes \cO(1)|_X)\otimes\alpha[1]\otimes\cO(-1)|_X)=(-)\otimes f_*\alpha[1]}(-)\otimes f_* f^*E[1],
\] where $\Psi$ is defined in
 Proposition \ref{proposition-twisted-atiyah-trivial-gerbe-formula}.
 We can also write this as 
\[
\id\xrightarrow[]{ \At_{\id_{\bL_{Y/k}}}(-)+(-)\otimes f_*g^*\psi[1]} (-)\otimes \bL_{Y/k}[1]\oplus  (-)\otimes f_*g^*\bL_{\bB\bG_m/k}[1]
\]
\[ 
\simeq (-)\otimes f_*\bL_{X/k}[1]\xrightarrow[]{(-)\otimes f_*\alpha[1]}
(-)\otimes f_*f^* E[1]\simeq (-)\otimes E[1].
\]
 By unraveling the bottom horizontal equivalence via Lemma \ref{lemma-dualizability-functor-Fourier-Mukai}, we see its image at the bottom right corner is 
\[
\Delta_*\cO_Y\xrightarrow[]{ (\At_{\id_{\bL_{Y/k}}}^{\FM}, \Delta_*f_*g^*\psi[1])} \Delta_* \bL_{Y/k}[1]\oplus \Delta_* f_*g^*\bL_{\bB\bG_m/k}[1]
\]
\[ \simeq
\Delta_*f_*\bL_{X/k}[1]\xrightarrow[]{\Delta_* f_*\alpha[1]}
\Delta_*f_*f^* E[1]\simeq \Delta_* E[1].
\]
To state this computation cleanly, consider the chain of equivalences
\[\hspace*{-1em}
\Map(\bL_{X/k},f^*E)\simeq \Map(f^*\bL_{Y/k},f^*E)\times \Map(g^*\bL_{\bB\bG_m/k},f^*E)\simeq  \Map(\bL_{Y/k},E)\times\Map(\cO_Y,E[1])
\]
where the second equivalence is constructed using the equivalence $f^*:\Dqc(Y)\simeq \Dqc^0(X/Y)$ and the trivialization $\psi:\cO_{\bB\bG_m}[-1]\simeq \bL_{\bB\bG_m/\Spec k}$ from Construction \ref{construction-psi-trivialize-LBGm}. With these identifications, the map from top left to bottom right  in the diagram (\ref{diagram-trivial-gerbe-Atiyah-map}) takes the form $\Map(\bL_{Y/k},E)\times\Map(\cO_Y,E[1])\to \Map_{Y\times_k Y}(\Delta_*\cO_Y,\Delta_*E[1])$. Our computation can then be summarized as follows:
\begin{proposition}\label{proposition-twisted-atiyah-trivial-gerbe-final-formula-general-E}
Let $E\in\Dqc(Y)$.
Up to the identifications made right before this proposition, the twisted universal $\id_{\bL_{X/k}}$-Atiyah transformation map \[\Map(\bL_{X/k},f^*E)\to\Map_{X\times_k X}(\iota_*^X\cP_X,\iota_*^X(\cP_X\otimes p_2^*f^*E[1]))\] can be written as
\[\hspace*{-1em}
\left((\Delta_*(-)[1])\circ \At^{\FM}_{\id_{\bL_{Y/k}}}
,\Delta_*(-)\right)
:\Map(\bL_{Y/k},E)\times\Map(\cO_Y,E[1])\to \Map_{Y\times_k Y}(\Delta_*\cO_Y,\Delta_*E[1]).\]
\end{proposition}
As a corollary, we relate this to the morphism appearing in the HKR isomorphism.

\begin{remark}\label{remark-recall-HKR}
For a smooth proper variety $Z$ over an algebraically closed field $L$ of characteristic 0, there is a quasi-isomorphism $I:\Delta^*\Delta_*\cO_Z\simeq \bigoplus_i\Omega_{Z/L}^i[i]$, and in \cite[Proposition 4.4]{Caldararu-Mukai-pairing-II} Caldararu showed that the following diagram commutes:
\[
\xymatrix{
\Delta_*\cO_Z\ar[r]^-\eta\ar[dr]_-{\exp\left(\At^{\FM}_{\id_{\Omega^i_{Z/L}}}\right)}&\Delta_*\Delta^*(\Delta_*\cO_Z)\ar[d]^-{\Delta_*I}\\
&\bigoplus_i\Delta_*\Omega_{Z/L}^i[i]
}
\]
Moreover, $I$ induces the HKR isomorphism
\[
\HT^*(Z/L)=\RHom_Z^*(\bigoplus_i\Omega^i_{Z/L}[i],\cO_Z)\simeq \RHom_Z^*(\Delta^*\Delta_*\cO_Z,\cO_Z)
\]
\[\simeq
\RHom_{Z\times_LZ}^*(\Delta_*\cO_Z,\Delta_*\cO_Z) =\HH^*(Z/L).
\]
The above commutative triangle hence relates the HKR isomorphism with the universal Atiyah class. 
Below, we will state a variant that applies to our situation. 
\end{remark}
Working with the variety $Y/k$ now, we consider the projection of $I$ to the first 2 summands (which still exists in positive characteristics; in the notation of \cite[\S3.18]{Moore2}, this is $c_0+c_1$ if we ignore the filtrations): 
\[I_{\leq 1}:\Delta^*\Delta_*\cO_Y\to \bigoplus_{i=0,1}\Omega^i_{Y/k}[i].
\]
The morphism $I_{\leq 1}$ induces a morphism
\[
\hspace*{-1em}
I^{\HKR}_{\leq1}:\RHom_Y(\bigoplus_{i=0,1}\Omega^i_{Y/k}[i],\cO_Y)\to \RHom_{ Y}(\Delta^*\Delta_*\cO_Y,\cO_Y)\simeq \RHom_{Y\times_k Y}(\Delta_*\cO_Y,\Delta_*\cO_Y).\]
We also define $\exp_{\leq1}\left(\At^{\FM}_{\id_{\bL_{Y/k}}}\right):\Delta_*\cO_Y\to \bigoplus_{i=0,1}\Delta_*\Omega^i_{Y/k}[i]$ as the sum of the identity map and $\At^{\FM}_{\id_{\bL_{Y/k}}}$.
\begin{lemma}[{Twisted analogue of \cite[Proposition 4.6(1)]{Moore2}; also see the proof of \cite[Corollary 1.4(1)]{Moore2}}]\label{lemma-injectivity-truncated-HKR}
The morphism
    \[ I^{\HKR}_{\leq1}:\RHom_Y(\bigoplus_{i=0,1}\Omega^i_{Y/k}[i],\cO_Y)\to
     \RHom_{Y\times_k Y}(\Delta_*\cO_Y,\Delta_*\cO_Y)
    \] induces an injective map on the $j$-th cohomology in the following cases:
    \begin{enumerate}
        \item either $\characteristic(k)=0$, or $j\leq p$ where $p=\characteristic(k)>0.$
        \item $p=2, j=3$, and $\rH^0(Y,\bigwedge^2\rT_Y)=0$.
    \end{enumerate}
\end{lemma}
\begin{proof}
\leavevmode
    This is proved in the same way as \cite[Proposition 4.6(1), Corollary 1.4(1)]{Moore2}, so we explain the modifications briefly. We factor the map $I_{\leq 1}$ as \[\Delta^*\Delta_*\cO_Y\to \tau_{<p}(\Delta^*\Delta_*\cO_Y)\to \bigoplus_{i=0,1}\Omega^i_{Y/k}[i].\]
    The complex $\tau_{<p}(\Delta^*\Delta_*\cO_Y)$ splits as a sum of its cohomology sheaves, and the map $I_{\leq 1}$ induces equivalences on $\rH^0$ and $\rH^{-1}$, 
    so the map $\tau_{<p}(\Delta^*\Delta_*\cO_Y)\to \bigoplus_{i=0,1}\Omega^i_{Y/k}[i]$ has a right inverse. It then remains to check the injectivity statement for the map 
    \[
\RHom_{Y}(\tau_{<p}(\Delta^*\Delta_*\cO_Y),\cO_Y)    \to \RHom_{Y}(\Delta^*\Delta_*\cO_Y,\cO_Y)\simeq \RHom_{Y\times_kY}(\Delta_*\cO_Y,\Delta_*\cO_Y).
    \] The rest of the proof is identical to Moore's argument. For instance, in order to check the case (2), 
    it suffices to show that 
    $\RHom^2_{Y}(\tau_{\geq2}(\Delta^*\Delta_*\cO_Y),\cO_Y)=0$, but this is isomorphic to $\RHom^2_{Y}(\cH^{-2}(\Delta^*\Delta_*\cO_Y)[2],\cO_Y)=\RHom^2_{Y}(\Omega_{Y/k}^2[2],\cO_Y)=\rH^0(Y,\bigwedge^2\rT_Y)=0.$
\end{proof}
\begin{proposition}[{Two-degree variant of computations in \cite{Moore2}; see the proof of \cite[Corollary 1.4(1)]{Moore2}}]\label{proposition-injectivity-truncated-HKR}
\leavevmode
\begin{enumerate}
\item The map
 $\rH^{1}(I^{\HKR}_{\leq1}):\rH^{0}(Y,\rT_Y)\oplus \rH^{1}(Y,\cO_Y)\to \HH^1(Y)$ is an isomorphism in any characteristic.
    \item The map
 $\rH^{2}(I^{\HKR}_{\leq1}):\rH^{1}(Y,\rT_Y)\oplus \rH^{2}(Y,\cO_Y)\to \HH^2(Y)$ is an isomorphism in any characteristic if  $\rH^0(Y,\bigwedge^2\rT_Y)=0.$
\item 
 The map
 $\rH^{3}(I^{\HKR}_{\leq1}):\rH^{2}(Y,\rT_Y)\oplus \rH^{3}(Y,\cO_Y)\to \HH^3(Y)$ is injective if $p\neq 2$ or $\rH^0(Y,\bigwedge^2\rT_Y)=0.$
\end{enumerate}
\end{proposition}
\begin{proof}
    By Lemma \ref{lemma-injectivity-truncated-HKR}, the maps are injective in all cases.
    
    For the surjectivity statements, recall that the $(-i)$-th cohomology sheaf of $\Delta^*\Delta_*\cO_Y$ is given by $\Omega^i_{Y/k}.$ 
    Then, applying the hypercohomology spectral sequence to
    $\HH^*(Y/k)\simeq \RHom_Y(\Delta^*\Delta_*\cO_Y,\cO_Y)$, we get a spectral sequence
    \[
    \rE^{p,q}_2=\rH^p(Y,\bigwedge^qT_{Y})\to\HH^{p+q}(Y).
    \]
    Then, we get $\dim\HH^1(Y/k)\leq\dim\rH^1(Y,\cO_Y)+\dim\rH^0(Y,\rT_Y)$, so (1) holds. 
    Also under the vanishing assumption of (2),  $\dim\HH^2(Y/k)\leq\dim\rH^2(Y,\cO_Y)+\dim\rH^1(Y,\rT_Y)$, so (2) holds.
\end{proof}

\begin{lemma}[Caldararu \cite{Caldararu-Mukai-pairing-II}]\label{lemma-caldararu-commutativity-truncated-version}
We have a commutative diagram, where $\eta$ is the unit of adjunction:
    \[
\xymatrix{
\Delta_*\cO_Y\ar[rr]^-\eta\ar[drr]_-{\exp_{\leq1}\left(\At^{\FM}_{\id_{\bL_{Y/k}}}\right)}&&\Delta_*\Delta^*(\Delta_*\cO_Y)\ar[d]^-{\Delta_*I_{\leq 1}}\\
&&\bigoplus_{i=0,1}\Delta_*\Omega^i_{Y/k}[i]
}
\]
\end{lemma}
\begin{proof}
As in Theorem 3.15 in the first version of the preprint \cite{Moore2}, this is an adaptation of Caldararu's result \cite[Proposition 4.4]{Caldararu-Mukai-pairing-II} to positive characteristics where we truncate at $\leq1$.  Caldararu's result requires characteristic assumptions only for constructing higher order terms, so the statement for the truncation $\leq1$ holds in any characteristic.
\end{proof}
\begin{theorem}\label{theorem-twisted-atiyah-trivial-gerbe-formula-HKR}
    Choose a connective object $V\in\Dperf(k)$, and consider the twisted Atiyah map for 
    $E=\cO_Y\otimes_kV\in \Dqc(Y)$: 
    \[
    \widetilde{\bAt}_{f^*E}: \Map(\bL_{X/k},\cO_X\otimes_kV)\to\Map_{X\times_k X}(\iota_*^X\cP_X,\iota_*^X(\cP_X\otimes_kV[1])).\] Under the same identifications as Proposition \ref{proposition-twisted-atiyah-trivial-gerbe-final-formula-general-E}, this map can be written as a morphism
\[
\Map(\bL_{Y/k},\cO_Y\otimes_kV)\times\Map(\cO_Y,\cO_Y\otimes_kV[1])\to \Map_{Y\times_k Y}(\Delta_*\cO_Y,\Delta_*\cO_Y\otimes_kV[1])
\]
This map coincides with $\Omega^\infty(I^{\HKR}_{\leq1}[1]\otimes_kV)$  
(i.e.\ the morphism we obtain by applying the functor $\Omega^\infty(-)$ to the map of $k$-module spectra).
\end{theorem}
\begin{proof} In Proposition \ref{proposition-twisted-atiyah-trivial-gerbe-final-formula-general-E}, we already computed this map as 
    $\left((\Delta_*(-)[1])\circ \At^{\FM}_{\id_{\bL_{Y/k}}}
,\Delta_*(-)\right).
$

On the other hand, 
as a result of Lemma \ref{lemma-caldararu-commutativity-truncated-version}, we obtain a commutative diagram:
\[
\begin{adjustbox}{max width=\linewidth}
\xymatrix{
\RHom_Y(\bigoplus_{i=0,1}\Omega_{Y/k}^i[i],\cO_Y)\ar[r]^-{(-)\circ I_{\leq 1}}\ar[d]^-{\Delta_*}&
\RHom_Y(\Delta^*\Delta_*\cO_Y,\cO_Y)\ar[ddr]^-{\simeq}\ar[d]^-{\Delta_*}\\
\RHom_{Y\times_kY}(\Delta_*\bigoplus_{i=0,1}\Omega^i_{Y/k}[i],\Delta_*\cO_Y)\ar[r]^-{(-)\circ \Delta_*I_{\leq 1}}\ar[drr]_-{(-)\circ \exp_{\leq1}(\At^{\FM}_{\id_{\bL_{Y/k}}})}&
\RHom_{Y\times_kY}(\Delta_*\Delta^*\Delta_*\cO_Y,\Delta_*\cO_Y)\ar[dr]^-{(-)\circ\eta}\\
&&\RHom_{Y\times_k Y}(\Delta_*\cO_Y,\Delta_*\cO_Y)
}
\end{adjustbox}
\]
The composition from top left to bottom right is precisely $I^{\HKR}_{\leq1}$, and the commutative diagram equates this with $(\Delta_*(-))\circ \exp_{\leq1}(\At^{\FM}_{\id_{\bL_{Y/k}}})$, so we are done. 

\end{proof}
As in  Definition \ref{definition-spectrally-enriched-twisted-atiyah-map}, the twisted Atiyah map can be upgraded to a map between mapping spectra. The identifications of mapping spaces used in Theorem \ref{theorem-twisted-atiyah-trivial-gerbe-formula-HKR} 
such as diagram (\ref{diagram-trivial-gerbe-Atiyah-map}) also upgrades to mapping spectra. With this modification, the same proof shows: 
\begin{theorem}[{spectral upgrade of Theorem \ref{theorem-twisted-atiyah-trivial-gerbe-formula-HKR}}]\label{theorem-spectral-variant-twisted-atiyah-trivial-gerbe-formula-HKR}
    Choose a connective object $V\in\Dperf(k)$, and consider the spectrally enriched twisted Atiyah map for 
    $E=\cO_Y\otimes_kV\in \Dqc(Y)$: 
    \[
    \widetilde{\bAt}_{f^*E}^{\FM}: \RHom(\bL_{X/k},\cO_X\otimes_kV)\to\RHom_{\Dqc^{-1,1}(X\times_kX/Y\times_kY)}(\iota_*^X\cP_X,\iota_*^X(\cP_X\otimes_kV[1])).\] We can identify the source and domain similarly to Proposition \ref{proposition-twisted-atiyah-trivial-gerbe-final-formula-general-E}, and we write this map as
\[
\RHom(\bL_{Y/k},\cO_Y\otimes_kV)\times\RHom(\cO_Y,\cO_Y\otimes_kV[1])\to \RHom_{Y\times_k Y}(\Delta_*\cO_Y,\Delta_*\cO_Y\otimes_kV[1])
\]
Then, this map coincides with  the map $I^{\HKR}_{\leq1}[1]\otimes_kV:$ \[
\RHom_Y(\bL_{Y/k},\cO_Y)\otimes_kV
\oplus \RHom_Y(\cO_Y,\cO_Y[1])\otimes_kV
\to 
\RHom_{Y\times_k Y}(\Delta_*\cO_Y,\Delta_*\cO_Y[1])\otimes_kV.\]
\end{theorem}
\part{Derived deformation theory}
\section{Formal moduli problems}
\subsection{Formal moduli problem of scalloped derived stacks}
Throughout the section, let $\Lambda$ be a noetherian complete local ring with residue field $k$. 
The following construction is analogous to \cite[\S3.6]{BranTael}.
Let $f:X\to Y$ be a morphism of scalloped qcqs 1-Artin derived stacks over $k$.

Using the inclusion $\{0\}\to\Delta^2$, we obtain a functor \[\Fun(\Delta^2, (\dArtSt^{\scalloped})^{\op})\to (\dArtSt^{\scalloped})^{\op}.\] Building on this functor and the functor $\Spec:\CAlgan\to(\dArtSt^{\scalloped})^{\op}$ we construct the $\infty$-category
\[
\cC:=\CAlgan\times_{(\dArtSt^{\scalloped})^{\op}}\Fun(\Delta^2, (\dArtSt^{\scalloped})^{\op}),
\] where $\CAlgan$ is the $\infty$-category of derived (or animated) rings. 
$\cC$ comes with a functor 
$p:\cC\to \CAlgan$, 
and this is a cocartesian fibration since $\dArtSt^{\scalloped}$ has fiber products (Proposition \ref{proposition-scalloped-dstack-stable-fiber-product}).
We take $\cC^{\cocart}\subset \cC$ to be the wide subcategory whose edges are the $p$-cocartesian ones.
The restriction $p:\cC^{\cocart}\to \CAlgan$ is a left fibration by \cite[Corollary 2.4.2.5]{HTT}. We view the triangle diagram formed by $f:X\to Y$ and $\Spec k$ as an object $(k,f)\in\cC^{\cocart}$, and we obtain 
the slice $\cC^{\cocart}_{/(k,f)}$ with a left fibration
\[
\cC^{\cocart}_{/(k,f)}\to \CAlgan_{/k}.
\]
By Grothendieck straightening procedure 
this defines a functor $\CAlgan_{/k}\to\widehat{\Cat_{\infty}}.$ Composing this with the forgetful functor $\CAlganArtLambdaK\to\CAlgan_{/k}$, we obtain 
\[
\Def_f:\CAlganArtLambdaK\to\widehat{\Cat_{\infty}},
\] where the superscript $(-)^{\art}$ denotes ``artinian'' (see \cite[Definition 3.1]{BranTael}).
Since endomorphisms of qcqs 1-Artin derived stacks that deform identity morphisms of qcqs 1-Artin derived stacks over $k$ to a base in  $\CAlganArtLambdaK$ are invertible (which can be proved similarly to Lemma \ref{lemma-thickened-end-aut} using Remark \ref{remark-square-zero-thickening-def}), 
$\Def_f$ factors through $\widehat{\cS}.$ 
This functor parametrizes spaces of deformations of the morphism $f$ (between scalloped qcqs 1-Artin derived stacks). 

Similarly to above, given a scalloped qcqs 1-Artin derived stack $X\to \Spec k$, we define a functor \[
\Def_X:\CAlganArtLambdaK\to\widehat{\cS}
\] which parametrizes spaces of deformations of the scalloped qcqs 1-Artin derived stack $X$.
\begin{lemma}\label{lemma-loopspace-Def-X-cotangent-complex}
Let $V\in\Mod_{k,\geq0}$ be perfect. We have
\[
\Omega\Def_X(\sqz_k(V))\simeq\Map(\bL_{X/k},\cO_X\otimes_kV),
\]
where the base point of the loopspace is at the trivial deformation. 
\end{lemma}
\begin{proof}
Since $X$ is a 1-Artin derived stack, it admits a cotangent complex. As in  
    Remark \ref{remark-cotangent-complex-nonconnective-universality}, by definition of cotangent complex, we have
    \[
    \Map_{X//k}(X[\cO_X\otimes_kV],X)\simeq \Map_X(\bL_{X/k},\cO_X\otimes_kV).
    \]
    The left hand side is identified with $\Omega\Def_X(\sqz_k(V))$ similarly to \cite[\S2.2.1]{Moore2}
\end{proof}
\begin{proposition}[See {\cite[Proposition 3.25]{BranTael}, \cite[Theorem 16.3.3.1]{SAG} for other situations}]\label{proposition-Def-X-FMP}
The functor $\Def_X$ is a formal moduli problem, so it is valued in $\cS$ in particular.
\end{proposition}
\begin{proof}
The proof is similar to \cite[Proposition 3.25]{BranTael}. 
The Schlessinger criterion for derived deformations of
scalloped qcqs 1-Artin derived stacks holds by Proposition \ref{proposition-schlessinger-criterion-scalloped-d-stack}.

For any perfect complex $V\in\Mod_{k,\geq0}$, the combination of the Schlessinger criterion and Lemma \ref{lemma-loopspace-Def-X-cotangent-complex} imply
\[
\Def_X(\sqz_k(V))\simeq \Omega\Def_X(\sqz_k(V[1]))\simeq \Map_X(\bL_{X/k},\cO_X\otimes_kV[1]).
\]
Thus, $\Def_X$ is valued in essentially small spaces if the inputs are of the form $\sqz_k(V)$. This can be generalized to other inputs by the argument of \cite[Proposition 3.25]{BranTael}.
\end{proof}
\begin{proposition}\label{proposition-omega-Def-f}
    Let $V\in\Mod_{k,{\geq0}}$ be perfect. Then, we have \[\Omega_{}\Def_f(\sqz_k(V))\simeq
    \Map(\bL_{X/k},\cO_X\otimes_kV)\times_{ \Map(f^*\bL_{Y/k},\cO_X\otimes_kV)}\Map(\bL_{Y/k},\cO_Y\otimes_kV)
    \] where the fiber product is defined using the maps
    \[
     f^*: \Map(\bL_{Y/k},\cO_Y\otimes_kV)\to \Map(f^*\bL_{Y/k},\cO_X\otimes_kV)),
    \]
    \[
   (-)\circ [f]:\Map(\bL_{X/k},\cO_X\otimes_kV)\to
    \Map(f^*\bL_{Y/k},\cO_X\otimes_kV)).
    \]
\end{proposition}
\begin{proof}
For shortness, we write $X^V:=X[\cO_X\otimes_kV],Y^V:=Y[\cO_Y\otimes_kV]$.
By definition of cotangent complex (see Remark \ref{remark-cotangent-complex-nonconnective-universality}), as in Lemma \ref{lemma-loopspace-Def-X-cotangent-complex} we have
\[
\Omega_{}\Def_X(\sqz_k(V))\simeq \Map(\bL_{X/k},\cO_X\otimes_kV),\quad
\Omega_{}\Def_Y(\sqz_k(V))\simeq \Map(\bL_{Y/k},\cO_Y\otimes_kV),\]
\[
\Map_{\sqz_k(V)}(X^V,Y^V)\times_{\Map_k(X,Y)}\{f\}\simeq 
\Map_k'(X^V,Y)\simeq \Map(f^*\bL_{Y/k},\cO_X\otimes_kV),
\] where $\Map_k'(X^V,Y)$ denotes the space of morphisms whose compositions with $X\to X^V$ come with identifications with $f:X\to Y.$

Base-changing along  $\Spec k\to\Spec \sqz_k(V)$ gives us a functor from \newline $\Fun(\Delta^2,\dArtSt^{\scalloped})\times_{\dArtSt^{\scalloped}}\{\Spec\sqz_k(V)\}$ to $ \Fun(\Delta^2,\dArtSt^{\scalloped})\times_{\dArtSt^{\scalloped}}\{\Spec k\}$. 
We can write $\Def_f(\sqz_k(V))$ as
the fiber of this functor over $\{(f,\Spec k)\}$.

Then, one can compute $\Omega_{}\Def_f(\sqz_k(V))$ by substituting the above expression into 
\[\{\triv\}\times_{\Def_f(\sqz_k(V))}\Fun(\Delta^1,\Def_f(\sqz_k(V)))\times_{\Def_f(\sqz_k(V))}\{\triv\}.\] 
After a series of simplifications, we see that a point of
$\Omega_{}\Def_f(\sqz_k(V))$ corresponds to a commutative diagram
\[
\xymatrix{
X^V\ar[r]^-\alpha\ar[d]^-{f^V}&X^V\ar[d]^-{f^V}\\
Y^V\ar[r]^-\beta& Y^V
}
\] over $\sqz_k(V)$ such that $f^V$ is the trivial deformation of $f,$ and its base change to $k$ is identified with the diagram
\[
\xymatrix{
X\ar[r]^-\id\ar[d]^-{f}&X\ar[d]^-{f}\\
Y\ar[r]^-\id& Y
}
\]
More precisely, we see that $\Omega_{\triv}\Def_f(\sqz_k(V))$ is equivalent to the fiber of the map
\[\hspace*{-3em}
\Map_{\sqz_k(V)}(X^V,X^V) \times_{\Map_{\sqz_k(V)} (X^V,Y^V)} \Map_{\sqz_k(V)}(Y^V,Y^V) 
\to 
{\Map_{k}(X,X) \times_{\Map_{k}(X,Y)} \Map_{k}(Y,Y)}\] above the point $(\id_X,\id_Y,\id_f:f\implies f)$, 
where the map 
\[\Map_{\sqz_k(V)}(X^V,X^V)\to {\Map_{\sqz_k(V)}(X^V,Y^V)}\]
is the post-composition map with $f^V$, and the map \[\Map_{\sqz_k(V)}(Y^V,Y^V)\to {\Map_{\sqz_k(V)}(X^V,Y^V)}\] is the pre-composition map with $f^V$.
Commuting pullbacks, we rewrite this as the fiber product of $\Map_{\sqz_k(V)}(X^V,X^V)\times_{\Map_{k}(X,X)}\{\id_X\}$ and $\Map_{\sqz_k(V)}(Y^V,Y^V)\times_{\Map_{k}(Y,Y)}\{\id_Y\}$ over $\Map_{\sqz_k(V)}(X^V,Y^V)\times_{\Map_{k}(X,Y)}\{f\}$.
In other words, 
\[
\Omega_{}\Def_f(\sqz_k(V))\simeq \Omega_{}\Def_X(\sqz_k(V))
\times_{\Map_{\sqz_k(V)}(X^V,Y^V)\times_{\Map_k(X,Y)}\{f\}
} \Omega_{}\Def_Y(\sqz_k(V)).
\]
By the comments at the beginning of the proof, this is identified with 
\[
\Map(\bL_{X/k},\cO_X\otimes_kV)\times_{\Map(f^*\bL_{Y/k},\cO_X\otimes_kV)} \Map(\bL_{Y/k},\cO_Y\otimes_kV).
\]
Let us describe the two structure maps in this fiber product to finish the proof. If we unravel the identifications, the left structure morphism
\[
A:\Map(\bL_{X/k},\cO_X\otimes_kV)\to \Map(f^*\bL_{Y/k},\cO_X\otimes_kV)
\]
comes from the post-composition morphism
\[
f^V\circ(-):\Omega\Def_{X}(\sqz_k(V))\to \Map_{\sqz_k(V)}(X^V,Y^V)\times_{\Map_k(X,Y)}\{f\},
\] so $A$ is the precomposition map by
 $[f]:f^*\bL_{Y/k}\to\bL_{X/k}$. 
The right structure morphism
\[B:
\Map(\bL_{Y/k},\cO_Y\otimes_kV)\to \Map(f^*\bL_{Y/k},f^*\cO_Y\otimes_kV).
\]comes from the pre-composition morphism
\[
(-)\circ f^V:\Omega_{}\Def_{Y}(\sqz_k(V))\to \Map_{\sqz_k(V)}(X^V,Y^V)\times_{\Map_k(X,Y)}\{f\},
\] so $B$ is the pullback map 
\[
f^*:\Map(\bL_{Y/k},\cO_Y\otimes_kV)\to \Map(f^*\bL_{Y/k},f^*\cO_Y\otimes_kV).
\] 
This concludes the proof.
\end{proof}

\begin{proposition}
The functor $\Def_f$ is a formal moduli problem, so it is valued in $\cS$ in particular.
\end{proposition}
\begin{proof}
    We first check the Schlessinger criterion. Consider a 
    pullback diagram
    \[
    \xymatrix{
A\ar[r]\ar[d]&B_1\ar[d]\\
B_2\ar[r]& C    
    }
    \] in $\CAlganArtLambdaK$ with surjective $\pi_0(B_i)\to \pi_0(C)$.
    For each object $R\in\CAlganArtLambdaK$, we 
    define an $\infty$-category 
    \[
\Mor_R:=\Fun(\Delta^2,(\dArtSt^{\scalloped})^{\op})\times_{(\dArtSt^{\scalloped})^{\op}} \{R\},
    \] where the fiber products are defined via the functors 
    \[\hspace*{-3em}
    (\{0\}\to\Delta^2):\Fun(\Delta^2,(\dArtSt^{\scalloped})^{\op})\to(\dArtSt^{\scalloped})^{\op}
    , 
    \Spec(-):\{R\}\to(\dArtSt^{\scalloped})^{\op}.
    \]
    Observe that $\Mor_R$ is equivalent to 
    $\Fun(\Delta^1,\dArtSt^{\scalloped}_{/R})^{\op}$. 
    It will suffice to check an equivalence of $\infty$-categories
    \[
\Mor_A\simeq \Mor_{B_1}\times_{\Mor_C}\Mor_{B_2}
    \] induced by base change. 
    This follows immediately by applying $\Fun(\Delta^1,-)^{\op}$ to the diagram as in Proposition \ref{proposition-schlessinger-criterion-scalloped-d-stack}:
    \[
    \Mor_A\simeq 
    \Fun(\Delta^1,\dArtSt^{\scalloped}_{/A})^{\op}\simeq
    \]\[
    \Fun\left(\Delta^1,\dArtSt^{\scalloped}_{/B_1}\times_{\dArtSt^{\scalloped}_{/C}}\dArtSt^{\scalloped}_{/B_2}\right)^{\op}\]
    \[\simeq 
\Fun(\Delta^1,\dArtSt^{\scalloped}_{/B_1})^{\op}\times_{\Fun(\Delta^1,\dArtSt^{\scalloped}_{/C})^{\op}} \Fun(\Delta^1,\dArtSt^{\scalloped}_{/B_2})^{\op}
    \simeq \]
    \[\Mor_{B_1}\times_{\Mor_C}\Mor_{B_2}.
    \]
    Essential smallness follows from Lemma \ref{proposition-omega-Def-f} similarly to Proposition \ref{proposition-Def-X-FMP}.
\end{proof}
\subsection{Deformations of $G$-gerbes}
{We fix a nice commutative group scheme $G$ over $\Lambda$ and follow our convention for principal $\bB G$-bundles and $G$-gerbes from Remark \ref{remark-segal-groupoid-assoc-group}. In particular, when $G$ is not smooth, $G$-gerbes should be interpreted in the fppf sense.} 
Assume in the situation of the previous subsection that $f:X\to Y$ is a principal $\bB G$-bundle. 
We define a subfunctor 
$\Def_f^{\Gerbe(G)}\subset \Def_f$ so that $\Def_f^{\Gerbe(G)}(A)\subset \Def_f(A)$ is spanned by the deformations of $f$ over $A$ that are $G$-gerbes for each $A\in\CAlganArtLambdaK$ (without trivialization of the banding). Note that this is well-defined as a functor since base changes of $G$-gerbes are still $G$-gerbes. 
Given the above setup, we have the classifying $k$-linear morphism $[X]:Y\to(\bB^2 G)_k$ for $f$. 
{We will use it to define a formal moduli problem $\Def_f^{\Bund(\bB G)}$ which parametrizes derived deformations of $f$ as principal $\bB G$-bundles. Informally, given a $A\in\CAlganArtLambdaK$, $\Def_f^{\Bund(\bB G)}(A)$ is the space of $A$-linear morphisms $\cY\to (\bB^2 G)_A$ from $1$-Artin scalloped derived stacks which deform the $k$-linear map $[X]:Y\to (\bB^2G)_k$. As an $\widehat{\cS}$-valued functor, this can be constructed similarly to $\Def_f$ from the previous section, by starting with the $\infty$-category $\cC$ define by the pullback 
\[
\xymatrix{
\cC\ar[r]\ar[d]&\Fun(\Delta^2,(\dArtSt^{2,\scalloped})^{\op})\ar[d]\\
\CAlg^{\an}_\Lambda\times(\dArtSt^{\scalloped})^{\op}\ar[r]^-{}&
\Fun(\Delta^1,(\dArtSt^{2,\scalloped})^{\op})\times{(\dArtSt^{2,\scalloped})^{\op}}
}
\] 
Here, $\dArtSt^{2,\scalloped}$ refers to the $\infty$-category of scalloped 2-Artin derived stacks. The right vertical functor comes from the inclusion $\Delta^{\{0,1\}}\coprod\{2\}\hookrightarrow\Delta^2$ of a vertex and an edge. The bottom horizontal functor arises from the combination of the functors 
\[\CAlgan_{\Lambda}\to
\Fun(\Delta^1,
(\dArtSt^{2,\scalloped})^{\op});A\mapsto (\Spec A\to (\bB^2G)_A)\]
(note the direction of the arrow due to taking the opposite category) 
and the inclusion $(\dArtSt^{\scalloped})^{\op}\to (\dArtSt^{2,\scalloped})^{\op}$. We will prove that the resulting functor is a formal moduli problem.
}
\begin{proposition}\label{proposition-Omega-Def-f-as-BGm-bundle} In the situation above, for any perfect $V\in\Mod_{k,\geq0}$, we have
    \[
    \Omega_{}\Def_f^{\Bund(\bB G)}(\sqz_k(V))\simeq \Map(\bL_{Y/(\bB^2 G)_k},\cO_Y\otimes_kV).
    \]
\end{proposition}
\begin{proof}

    $\Omega_{}\Def_f^{\Bund(\bB G)}(\sqz_k(V))$ is the space of automorphisms of the morphism $f^V:Y^V\to (\bB^2 G)_k^V$ with the identity on $(\bB^2 G)_k^V.$ Such automorphisms are the automorphisms of $Y^V$ respecting the morphism $f^V$ whose base changes along $\Spec k\to \Spec \sqz_k(V)$ are identified with the identity of $Y$. In other words, this space is $\Map_{(\bB^2 G)_k^V}(Y^V,Y^V)\times_{\Map_{(\bB^2 G)_k}(Y,Y)}\{\id_Y\}$. Universality of fiber product says 
    $\Map_{(\bB^2 G)_k^V}(Y^V,Y^V)\simeq \Map_{(\bB^2 G)_k}(Y^V,Y)$, so this space is the same as $\Map_{Y//(\bB^2 G)_k}(Y^V,Y)$ which by Remark \ref{remark-cotangent-complex-nonconnective-universality} is the same as $\Map(\bL_{Y/(\bB^2 G)_k},\cO_Y\otimes_kV)$. 
\end{proof}
{
\begin{proposition}
    Given the following pullback diagram in $\CAlganArtLambdaK$ with surjective $\pi_0(A_i)\to \pi_0(A_0)$ ($i=1,2$):
    \[
    \xymatrix{
A_3\ar[r]\ar[d]&A_1\ar[d]\\
    A_2\ar[r]& A_0
    }
    \]  applying the functor $\Def_f^{\Bund(\bB G)}$ gives us a  pullback diagram in $\widehat{\cS}$:
    \[
    \xymatrix{
\Def_f^{\Bund(\bB G)}(A_3)\ar[r]\ar[d]&\Def_f^{\Bund(\bB G)}(A_1)\ar[d]\\
   \Def_f^{\Bund(\bB G)}( A_2)\ar[r]& \Def_f^{\Bund(\bB G)}(A_0)
    }
    \]
\end{proposition}
\begin{proof}
   We have an obvious morphism from $\Def_f^{\Bund(\bB G)}(A_3)$ to the pullback of the remaining three terms by base change, and we describe the inverse construction. Thus, start with $[f_i:Y_i\to (\bB^2 G)_{A_i}]\in \Def_f^{\Bund(\bB G)}( A_i)$ for $i=0,1,2$ along with an identification $f_1\times_{A_1}A_0\simeq f_0\simeq f_2\times_{A_2}A_0$. First, by the proof of Proposition \ref{proposition-schlessinger-criterion-scalloped-d-stack}, we can find a 1-Artin scalloped derived stack $Y_3$ over $A_3$ which is the geometric pushout of $Y_1\leftarrow{Y_0}\to Y_2$ (and the $Y_i$'s for $i=0,1,2$ are the base changes of $Y_3$ along the maps $A_3\to A_i$).
   We observe that 
   \[
   \Map_{A_3}(Y_3,(\bB^2G)_{A_3})\simeq \Map_{A_3}(Y_1,(\bB^2G)_{A_3})\times_{\Map_{A_3}(Y_0,(\bB^2G)_{A_3})} \Map_{A_3}(Y_2,(\bB^2G)_{A_3}).
   \] Indeed, geometric derived stacks are infinitesimally cohesive by Artin-Lurie representability theorem (see for example \cite[Proposition 5.3.7]{DAG}), so this holds if we replace $Y_3$ by a smooth atlas which is a 0-geometric derived stack (and  replace the other $Y_i$'s by base changes of this atlas). The same conclusion also holds for other terms of the Cech nerve of the smooth atlas. By taking the colimits of the Cech nerves, we obtain the version for the $Y_i$'s.

   By base change, the above equivalence can be rewritten as 
   \[
   \Map_{A_3}(Y_3,(\bB^2G)_{A_3})\simeq \Map_{A_1}(Y_1,(\bB^2G)_{A_1})\times_{\Map_{A_0}(Y_0,(\bB^2G)_{A_0})} \Map_{A_2}(Y_2,(\bB^2G)_{A_2}),
   \]
   where the morphism from left to right is interpreted as base changes from $A_3$ to $A_i (i=0,1,2)$.
   This lets us construct $f_3:Y_3\to (\bB^2G)_{A_3}$, whose base change recovers the given data. Moreover, by construction, this is inverse to the base change morphism $\Def_f^{\Bund(\bB G)}(A_3)\to \Def_f^{\Bund(\bB G)}(A_1)\times_{\Def_f^{\Bund(\bB G)}(A_0)} \Def_f^{\Bund(\bB G)}(A_2).
   $
\end{proof}
\begin{corollary}\label{corollary-BG-bundle-formal-moduli-problem}
     In the situation above, $\Def_f^{\Bund(\bB G)}$ is an (${\cS}$-valued) formal moduli problem.
\end{corollary}
\begin{proof}
    This follows from the last two propositions. 
\end{proof}
}
\begin{lemma}\label{lemma-compute-fiber-Def-f-for-Gm-gerbes}
    In the situation above,
    \begin{enumerate}
        \item pulling back along $f$ induces equivalences \[
         \Map(\bL_{Y/(\bB^2G)_k},\cO_Y)\to \Map(\bL_{X/k},\cO_X),\quad 
        \Map(\bL_{Y/k},\cO_Y)\simeq \Map(f^*\bL_{Y/k},\cO_X),
        \] and
        \item we have a cartesian square \[
        \xymatrix{
\Map(\bL_{X/k},\cO_X)\ar[r]\ar[d]^-{\id}&\Map(\bL_{Y/k},\cO_Y)\ar[d]^-{f^*}\\
\Map(\bL_{X/k},\cO_X)\ar[r]^-{(-)\circ[f]}&\Map(f^*\bL_{Y/k},\cO_X)
        }
        \]
        where the top horizontal map is the composition
        \[
\Map(\bL_{X/k},\cO_X)\xrightarrow[]{(-)\circ[f]}\Map(f^*\bL_{Y/k},\cO_X)\simeq \Map(\bL_{Y/k},f_*\cO_X)\simeq \Map(\bL_{Y/k},\cO_Y).
        \]
        In particular, the projection map
        \[
\Map(\bL_{X/k},\cO_X)\times_{\Map(f^*\bL_{Y/k},\cO_X)}\Map(\bL_{Y/k},\cO_Y)\to \Map(\bL_{X/k},\cO_X)
        \] is an equivalence.
    \end{enumerate}
\end{lemma}
\begin{proof}\leavevmode
    \begin{enumerate}
        \item 
        We have $X\simeq Y\times_{(\bB^2 G)_k}\Spec k$ with respect to the classifying morphism $[X]:Y\to(\bB^2 G)_k$ and the universal principal $\bB G$-bundle $\pi:\Spec k\to(\bB^2 G)_k$.
        The morphism $\cO_Y\to f_*\cO_X$ is an equivalence by Lemma \ref{lemma-nice-group-derived-pushforard-structure-sheaf}. The results follow from the adjunction $(f^*,f_*)$  
         once we observe that $f^*\bL_{Y/(\bB^2G)_k}\simeq\bL_{X/k}.$
        \item 
         If we replace the top right vertex by $\Map(\bL_{Y/k},f_*\cO_X)$, 
        the top horizontal map  by $\Map(\bL_{X/k},\cO_X)\xrightarrow[]{(-)\circ[f]}\Map(f^*\bL_{Y/k},\cO_X)\simeq \Map(\bL_{Y/k},f_*\cO_X)$ and the right vertical map by $\Map(\bL_{Y/k},f_*\cO_X)\to \Map(f^*\bL_{Y/k},f^*f_*\cO_X)\simeq \Map(f^*\bL_{Y/k},\cO_X)$, then it's easy to see this commutes by formal computations.
        Up to the equivalence \newline $\Map(\bL_{Y/k},\cO_Y)\simeq \Map(\bL_{Y/k},f_*\cO_X)$, we obtain the desired commutative diagram.
        The vertical maps are equivalences by the previous part, so this is cartesian.
    \end{enumerate}
\end{proof}
\begin{proposition}\label{proposition-equiv-FMP-f-X-Gmgerbe-BGmbundle} The morphisms of functors
    \[
    \Def^{\Bund(\bB G)}_f\to \Def_f^{\Gerbe(G)}\to \Def_f\to \Def_X
    \]are all equivalences. In particular, $\Def_f^{\Gerbe(G)}$ is a formal moduli problem.
\end{proposition}
\begin{proof}
\leavevmode
We will verify the claim by checking certain compositions are equivalences as follows:
    \begin{enumerate}
        \item $\underline{\Def_f\to \Def_X}:$ By construction of the morphism and Proposition \ref{proposition-omega-Def-f}, the morphism $\Omega_{}\Def_f(\sqz_k(V))\to\Omega_{} \Def_X(\sqz_k(V))$ is identified with the projection map
        \[
        \Map(\bL_{X/k},\cO_X\otimes_kV)\times_{\Map(f^*\bL_{Y/k},\cO_X\otimes_kV)}\Map(\bL_{Y/k},\cO_Y\otimes_kV)\to \Map(\bL_{X/k},\cO_X\otimes_kV).
        \] This is an equivalence by Lemma \ref{lemma-compute-fiber-Def-f-for-Gm-gerbes} (2), so we obtain an equivalence of tangent fibres $\rT_{\Def_f}\simeq \rT_{\Def_X}$ and formal moduli problems $\Def_f\simeq \Def_X$ (\cite[Proposition 3.9]{BranTael}).
        \item $\underline{\Def_f^{\Bund(\bB G)}\to \Def_X}:$ By unraveling the relation between principal bundles and classifying morphisms, we characterize the morphism 
        \[\Omega_{}\Def_f^{\Bund(\bB G)}(\sqz_k(V))\to \Omega_{}\Def_f(\sqz_k(V))\] as follows.
        Given an automorphism $\alpha:Y^V\to Y^V$ of $Y^V$ which deforms $\id_Y$ compatibly with $[X]^V:Y^V\to(\bB^2 G)_{\sqz_k(V)}$ (which is the trivial deformation of $[X]$), this is sent to the datum of the following commutative diagram  over $\sqz_k(V)$:
        \[
        \xymatrix{
        X^V\ar[r]^-{\alpha_X}\ar[d]^-{f^V}&X^V\ar[d]^-{f^V}\\
        Y^V\ar[r]^-{\alpha}&Y^V
        }
        \] Here, $\alpha_X:Y^V\times_{(\bB^2 G)_{\sqz_k(V)}}\Spec (\sqz_k(V))\to Y^V\times_{(\bB^2 G)_{\sqz_k(V)}}\Spec (\sqz_k(V))$ is the base change of $\alpha$ along $\pi^V:\Spec (\sqz_k(V))\to (\bB^2 G)_{\sqz_k(V)}$. Further applying the morphism $\Omega_{}\Def_f(\sqz_k(V))\to \Omega_{}\Def_X(\sqz_k(V))$ will send this diagram to $\alpha_X:X^V\to X^V$ over $\sqz_k(V)$.
        This description and Proposition \ref{proposition-Omega-Def-f-as-BGm-bundle} show that the morphism $\Omega_{}\Def_f^{\Bund(\bB G)}(\sqz_k(V))\to \Omega_{} \Def_X(\sqz_k(V))$ can be identified with the morphism induced by $f^*$ (as it comes from the base change along $\pi:\Spec k\to (\bB^2 G)_k$):
        \[
\Map(\bL_{Y/(\bB^2 G)_k},\cO_Y\otimes_kV)\to\Map(\bL_{X/k},\cO_X\otimes_kV).
        \]
        By Lemma \ref{lemma-compute-fiber-Def-f-for-Gm-gerbes}(1), this is an equivalence, and so is the map $\Def_f^{\Bund(\bB G)}\to \Def_X$.
        \item $\underline{\Def_f^{\Gerbe(G)}\to \Def_f}:$ For each $A\in\CAlganArtLambdaK$, $\Def_f^{\Gerbe(G)}(A)\subset \Def_f(A)$ is a full subcategory, so we only need to check essential surjectivity of the inclusion. This follows since the equivalence $\Def_f^{\Bund(\bB G)}\to \Def_X$ factors through $\Def_f^{\Gerbe(G)}$.
    \end{enumerate}
\end{proof}
\begin{remark}
    Proposition \ref{proposition-equiv-FMP-f-X-Gmgerbe-BGmbundle} is related to the result \cite[Theorem 2.11]{AHPS} where they deform good moduli space morphisms along square zero extensions of the source. In its proof, they consider a good moduli space morphism $q:X\to Y$ from an Artin derived stack to a derived algebraic space, and assume they have a square zero extension $X'$ of $X$. Then, they construct a deformation $q':X'\to Y'$ as a good moduli space morphism in a natural way. 
\end{remark}
\subsection{Deforming principal $\bB\bG_m$-bundles along square-zero thickenings}\label{subsection-deform-BGm-bundles-sqz-thickening}
Let $f:X\to Y$ be a principal $\bB\bG_m$-bundle over a scalloped derived stack $Y$ over $k$. 
Let $E\in\Dqccn^0(X)$ be a $0$-twisted connective object, and let $X':=X[E]$ be the trivial square-zero thickening by $E.$ For instance, $X'$ can be a derived deformation of $X$ with respect to some base derived ring extension (as treated in the previous subsection), but we need this more general version.
Let $i^X:X\to X[E],i^Y:Y\to Y[f_*E]$ be the inclusions and let $\pi^X:X[E]\to X, \pi^Y:Y[f_*E]\to Y$ be the projections. 
In this subsection, we construct a thickening $f':X'\to Y'$ of the morphism $f$ which is a  principal $\bB\bG_m$-bundle and satisfies some desirable properties. Our construction is inspired by \cite{AHPS}.

\begin{lemma}\label{lemma-trivial-sqz-thickening-scalloped}
 $X[E]$ is a scalloped qcqs 1-Artin derived stack.
\end{lemma}
\begin{proof}
    The structure morphism $X[E]\to X$ is affine, so when $X$ is qcqs, $X[E]$ is as well. 
    We observe that $(X[E])_{\cl}\simeq X_{\cl}[\pi_0(E)]$, which can be checked affine-locally. 
    This shows $(X[E])_{\cl}$ is a thickening of $X_{\cl},$ and Theorem \ref{theorem-khan-ravi-scalloped-iff-nicely-stabilized} ensures that $X[E]$ is scalloped.
\end{proof}

For our purposes, we will define the following notion:

\begin{definition}
\label{definition-square-zero-thickeninig-stack}
    A morphism $f:X\to Y$ of 1-Artin derived stacks is called \textbf{a square-zero thickening} if it is
    of form $\Spec_Y(\cA)\to Y$ where $\cA$ is a connective $\cO_Y$-algebra with connective $\cJ:=\fib(\cO_Y\to\cA)$, and  the morphism $\cJ\otimes\cJ\to\cO_Y\otimes\cO_Y\simeq\cO_Y$ is trivial in $\Dqc(Y)$. 
\end{definition}
\begin{remark}\label{remark-square-zero-thickening-def}
\leavevmode
\begin{enumerate}
    \item When $X$ is 1-Artin and $E\in\Dqccn(X)$ is connective,
    the closed immersion $j:X\to X[E]$ is a square-zero thickening in the above sense. 
    \item Square-zero thickenings are clearly stable under base change. 
    \item This differs from the standard notion of ``square-zero extension'' in the literature (\cite[Definition 7.4.1.6]{HA}). However, these notions are closely related by \cite[Proposition 7.4.1.14]{HA}.
\end{enumerate}

\end{remark}
\begin{lemma}\label{lemma-sqz-thickening-pullback-conservative}
    Let $\Spec B\to \Spec A$ be a square-zero thickening. Then, the base change functor $(-)\otimes_AB:\Mod_A\to\Mod_B$ is conservative.
\end{lemma}
\begin{proof}
Let $I=\fib(A\to B)$. 
    Choose a morphism $f:M\to N$ in $\Mod_A$ which induces $M\otimes_AB\simeq N\otimes_AB$. We need to show $f:M\simeq N$. By taking the fiber of $f$, we reduce to showing that $K\simeq0$ when $K\otimes_AB\simeq 0$. Suppose $K\otimes_AB\simeq 0$.
    Applying $(-)\otimes_AB$ to the fiber sequence $I\to A\to B$, we obtain the equivalence $K\otimes_AI\simeq K.$ Then, we obtain $K\otimes_AI\otimes_AI\simeq K\otimes_AI\simeq K$, but this map is obtained by applying $K\otimes_A(-)$ to the nullhomotopic map $I\otimes_AI\to A$, so we see $K\simeq 0.$
\end{proof}

\begin{lemma}\label{lemma-affineness-derived-vs-classical}
    Let $X$ be a 1-Artin derived stack. 
    Then $X$ is derived affine if and only if $X_{\cl}$ is affine.
\end{lemma}
\begin{proof}
    Only if direction is obvious so assume $X_{\cl}$ is affine.
     We observe that $X$ is nilcomplete as it's geometric. 
    Hence, {for any test derived affine $B$, we have $X(B)\simeq \lim_nX(B_n)$, where $B_n:=\tau_{\leq n}(B)$.
    We write $X_n$ for the $n$-th truncation of $X$, so that $X_{\cl}=X_0$ in particular. 
    }  
    Since $X_n$ is the pushout of the diagram (by a standard fact about square-zero Postnikov extensions)
    \[X_{n-1}\xleftarrow[]{\pi} X_{n-1}[\pi_n(\cO_X)[n+1]]\xrightarrow[]{d} X_{n-1}\]
    in the $\infty$-category of 1-Artin derived stacks,  
    when $X_{n-1}$ is derived affine, so is $X_n$. Thus, each $X_n$ is derived affine by induction.
  {  We can then write $X_n=\Spec A_n$, compatibly for various $n$.
    We set $A:=\lim_nA_n$
    We see that for any derived affine $\Spec B$,
\[
X(B)\simeq \lim_nX(B_n)
\simeq \lim_n\Map(\Spec B_n,X)
\simeq \lim_n\Map(\Spec B_n,X_n)
\simeq \lim_n\Map(A_n, B_n)\]
\[
\simeq \lim_n\Map( A,B_n)
\simeq \Map(A,B)
\simeq (\Spec A)(B).
\]
   Thus, $X\simeq\Spec A$ is affine.
   }
\end{proof}
\begin{lemma}\label{lemma-sqz-thickening-derived-vs-classical}
    When $f:X\to Y$ is a square-zero thickening of 1-Artin derived stacks (Definition \ref{definition-square-zero-thickeninig-stack}), 
    the classicalization $X_{\cl}\to Y_{\cl}$ is a square-zero thickening in the traditional sense. 
\end{lemma}
\begin{proof}
The classical locus functor is compatible with pullbacks and colimits along Cech nerves. 
Also, classicalization of a smooth cover morphism is a smooth cover in the classical sense,
so we may check this locally on $Y$. Then, we can assume $Y$ is affine, and write $f$ as $\Spec A\to \Spec B$ for derived rings $A,B$. Let $I=\fib(B\to A).$
Then, $\pi_0(B)\to\pi_0(A)$ is surjective with kernel $\text{Im}(\pi_0(I)\to\pi_0(B))$.
By assumption, the map \[\pi_0(I)\otimes_{\pi_0(B)}\pi_0(I)\simeq \pi_0(I\otimes_B I)\to\pi_0(B\otimes_BB)\simeq\pi_0(B)\]
is trivial, so $\text{Im}(\pi_0(I)\to\pi_0(B))$ is square-zero in the traditional sense.

\end{proof}
\begin{corollary}\label{corollary-sqz-thickening-derived-affineness}
    Let $X\to Y$ be a square-zero thickening of 1-Artin derived stacks (Definition \ref{definition-square-zero-thickeninig-stack}).
  If $X$ is derived affine, then so is $Y$.
\end{corollary}
\begin{proof}
 By Lemma \ref{lemma-affineness-derived-vs-classical},
 $X_{\cl}$ is affine. By Lemma \ref{lemma-sqz-thickening-derived-vs-classical}, $X_{\cl}\to Y_{\cl}$ is a square-zero thickening in the classical sense and hence $Y_{\cl}$ is affine. Then, using Lemma  \ref{lemma-affineness-derived-vs-classical} again, we see that $Y$ is derived affine.
\end{proof}

\begin{lemma}\label{lemma-thickened-end-aut}
    Let $X,Y,Z,W$ be 1-Artin derived stacks. Suppose that we have a cartesian square where the vertical maps are square-zero thickenings:
    \[
    \xymatrix{
    X\ar[r]^\psi\ar[d]&Y\ar[d]\\
    Z\ar[r]^-\phi& W
    }
    \]
    If $\psi$ is an equivalence, so is $\phi$.
\end{lemma}
\begin{proof}
Since we can check equivalence locally and square-zero thickenings are preserved under base change, we may assume $W$ is derived affine. Square-zero thickenings are affine, and $\psi$ is an equivalence, so $X$ and $Y$ are also affine. By Corollary
\ref{corollary-sqz-thickening-derived-affineness}, $Z$ is also derived affine. At this point, the result follows from Lemma \ref{lemma-sqz-thickening-pullback-conservative} since the forgetful functor sending algebras over a derived ring to modules is conservative.  
\end{proof}
By Proposition \ref{proposition-BGmbundle-pushforward}, $E\simeq f^*f_*E$, and $f_*$ is $t$-exact, so that $f_*E$ is connective. We define $Y':=Y[f_*E]$. 
By the compatibility of relative spectrum with base change, we have a morphism
\[
f':X'=X[E]\simeq X[f^*f_*E]\simeq Y[f_*E]\times_YX\to Y[f_*E]=Y'.
\] In other words, $f'$ is the base change of $f$ along the projection $\pi_Y:Y'\to Y$. $f'$ carries a principal $\bB\bG_m$-bundle structure from $f.$ 
\begin{lemma}\label{lemma-invariance-etale-topos-thickening-bZ2Z}
    The  composition with $i_Y$ induces an equivalence
    \[
    \Map(Y',\bB\ZModTwo)\simeq  \Map(Y,\bB\ZModTwo).
    \]
\end{lemma}
\begin{proof} 
    We can rewrite this morphism as 
    \[
    ({i^Y})^*:\Map_{Y'}(Y',(\bB\ZModTwo)_{Y'})\to \Map_Y(Y,(\bB\ZModTwo)_{Y}).
    \]
    Since $\ZModTwo$ is \'etale over the base, we can interpret these mapping spaces in the underlying \'etale topoi of $Y$ and $Y'$. By the invariance of \'etale topoi under thickenings (see for example \cite[Reminders 1.6.6]{HHR}), this is an equivalence.
\end{proof}

\begin{proposition}\label{proposition-sqz-BGm-bundle-pullback}
\leavevmode
Consider the morphisms $f:X\to Y,f':X'\to Y'$ defined at the beginning of the subsection. 
In this situation the following statements hold. 
    \begin{enumerate}
    \item We have a pullback diagram
    \[
     \xymatrix{
        X\ar[r]^-{i^X}\ar[d]^-{f}&X'\ar[d]^-{f'}\\
        Y\ar[r]^-{i^Y}&Y'
        }
    \] of principal $\bB\bG_m$-bundles.
        \item Let $[f]^A:\bL_{Y/k}\to f_*\bL_{X/k}$ denote the morphism associated with the $k$-linear morphism $f.$
        For any morphism $d\in\Map_{X//k}(X',X)\simeq\Map(\bL_{X/k},E)$, there is a canonical pullback diagram of principal $\bB\bG_m$-bundles 
        \[
        \xymatrix{
        X'\ar[r]^-{d}\ar[d]^-{f'}&X\ar[d]^-{f}\\
        Y'\ar[r]^-{d^Y}&Y
        }
        \] whose construction is functorial in $E$, where $d^Y\in\Map_{Y//k}(Y',Y)\simeq\Map(\bL_{Y/k},f_*E)$ is the image of $d$ via $\Map(\bL_{X/k},E)\xrightarrow{f_*} \Map(f_*\bL_{X/k},f_*E)\xrightarrow[]{(-)\circ [f]^A} \Map(\bL_{Y/k},f_*E)$.  
    \end{enumerate}
\end{proposition}
\begin{proof}
\leavevmode
    \begin{enumerate}
        \item In the following commutative diagram, the big square and the right square are pullbacks of principal $\bB\bG_m$-bundles, so the left square is also a pullback.
         \[
     \xymatrix{
        X\ar[r]^-{i^X}\ar[d]^-{f}&X'\ar[d]^-{f'}\ar[r]^-{\pi^X}&X\ar[d]^-f\\
        Y\ar[r]^-{i^Y}&Y'\ar[r]^{\pi^Y}&Y
        }
    \] 
        \item If we ignore the pullback condition and principal $\bB\bG_m$-bundle structures  for the moment, the prescribed morphism $d_Y$ fits in a commutative diagram as in the proof of \cite[Theorem 2.12]{AHPS}. One way to describe this construction functorially in $E$ is as follows. Similarly to Lemma \ref{lemma-compute-fiber-Def-f-for-Gm-gerbes}(2), since $E$ is 0-twisted, we have a pullback square
        \[\tag{*}\label{*}
        \xymatrix{
\Map(\bL_{X/k},E)\ar[r]\ar[d]^-{\id}&\Map(\bL_{Y/k},f_*E)\ar[d]^-{f^*}\\
\Map(\bL_{X/k},E)\ar[r]^-{(-)\circ[f]}&\Map(f^*\bL_{Y/k},E)
        }
        \] where the top horizontal arrow is defined as in the statement of the proposition.
        As in the proof of Proposition \ref{proposition-omega-Def-f}, this pullback is the space of endomorphisms from $f'$ to itself which ``restrict back'' to $f$ upon composing with the morphism $f\to f'$ given by the left square of the diagram in the proof of part (1). 
        In particular, a point on the pullback gives us a diagram
        \[
        \xymatrix{
        X'\ar[r]^-{d}\ar[d]^-{f'}&X\ar[d]^-{f}\\
        Y'\ar[r]^-{d^Y}&Y
        }
        \] 
        by forgetting the restriction data. The diagram $(\ref{*})$ is functorial in $E$, so this justifies functoriality.

        Now, we check the resulting diagram of derived stacks is cartesian. 
        We have a commutative diagram
        \[\label{**}\tag{**}
     \xymatrix{
        X\ar[r]^-{i^X}\ar[d]^-{f}&X'\ar[d]^-{f'}\ar[r]^-{d}&X\ar[d]^-f\\
        Y\ar[r]^-{i^Y}&Y'\ar[r]^{d^Y}&Y
        }
    \] where the horizontal compositions are homotopic to the identity maps by the ``restriction data.'' 
    We let $P$ denote the pullback of the right square, and $\phi:X'\to P$ is a morphism defined by the universal property. We hope to prove $\phi$ is an equivalence.
    We consider base changing $\phi$ along $i^Y.$ Since the left square of the diagram is cartesian by part (1),  
    The resulting map can be written as $\phi\times_{Y'}Y:X\to  P\times_{Y'}Y$.
    By the pasting law of cartesian squares, 
    the square formed by 
    $\phi,\phi\times_{Y'}Y, i^X$, and $i^Y\times_{Y'}Y$ is cartesian. Since the big square is cartesian, $\phi\times_{Y'}Y$ is an equivalence. 
    Also observe that the maps $i^X, i^Y\times_{Y'}Y$ are square-zero thickenings. 
    By Lemma \ref{lemma-thickened-end-aut}, $\phi$ is an equivalence, so 
this shows the cartesian property if we forget the bundle structure. 

Lastly, we show that this cartesian square in fact respects the $\bB\bG_m$-action (i.e.\ it defines the base change square of principal $\bB\bG_m$-bundles).
Consider the following commutative diagram constructed using Example \ref{example-gerbe-banding-fiber-sequence}, where the vertical maps are induced by $i^Y$:
\[
\begin{tikzcd}[column sep={9em,between origins},row sep=4em]
\Map(Y',\bB^2\bG_m)
  \arrow[rr,""]
  \arrow[dr,"j'"]
  \arrow[dd,"i^Y_{\bB^2\bG_m}"]
&&
\Map(Y',\{*\})\simeq \{*\}
  \arrow[dd,equal]
  \arrow[dr]
\\
& \Map(Y', \bB\sAut(\bB\bG_m))
  \arrow[rr,"\tau'"]
  \arrow[dd,near start, "i^Y_{\bB\sAut(\bB\bG_m)}"]
&&
 \Map(Y',\bB\ZModTwo)\arrow[dd,"\simeq","i^Y_{\bB\ZModTwo}"']
\\
 \Map(Y,\bB^2\bG_m)
  \arrow[rr,near start, ""]
  \arrow[dr,,"j"]
&&
\Map(Y,\{*\})\simeq \{*\}
  \arrow[dr]
\\
& \Map(Y,\bB\sAut(\bB\bG_m))
  \arrow[rr,"\tau"]
&&
 \Map(Y,\bB\ZModTwo)
\end{tikzcd}
\]
In this diagram, $i^Y_{\bB\ZModTwo}$ is an equivalence by
Lemma \ref{lemma-invariance-etale-topos-thickening-bZ2Z}, and the two horizontal faces are cartesian. 
In the following argument, $d^Y_{\bB^2\bG_m}:\Map(Y,\bB^2\bG_m)\to \Map(Y',\bB^2\bG_m)$ is the map induced by $d^Y$, and we similarly define 
\[
d^Y_{\bB\sAut(\bB\bG_m)}:\Map(Y,\bB\sAut(\bB\bG_m))\to \Map(Y',\bB\sAut(\bB\bG_m))\]
and 
$d^Y_{\bB\ZModTwo}:\Map(Y,\bB\ZModTwo)\to \Map(Y',\bB\ZModTwo).$
View $f$ as defining a point $[f]\in\Map(Y,\bB^2\bG_m)$ and $f'$ as defining a point $[f']\in\Map(Y',\bB^2\bG_m)$. 
In view of the correspondence between $\bG_m$-gerbes and $\sAut(\bB\bG_m)$-bundles (Proposition \ref{proposition-gerbe-BAutBG}), the pullback square involving $d$ and $d^Y$ gives us a path $\ell:\Delta^1\to \Map(Y',\bB\sAut(\bB\bG_m))$ connecting the points $j'([f'])$ and $j'(d^Y_{\bB^2\bG_m}([f]))\simeq d^Y_{\bB\sAut(\bB\bG_m)}(j([f]))$. The goal is to lift $\ell$ to a path $\widetilde{\ell}:\Delta^1\to \Map(Y',\bB^2\bG_m)$ between $[f']$ and $d^Y_{\bB^2\bG_m}[f]$.  

The left square and the big square in (\ref{**}) are pullback diagrams of principal $\bB\bG_m$-bundles by part (1), so this gives us a path $\ell_0:\Delta^1\to \Map(Y,\bB^2\bG_m)$ between $ i^{Y}_{\bB^2\bG_m}([f'])$ and 
 $ i^{Y}_{\bB^2\bG_m}(d_Y^{\bB^2\bG_m}([f]))\simeq [f]$. Forgetting the $\bB\bG_m$-action in the diagram (\ref{**}) yields a path between $ j(i^{Y}_{\bB^2\bG_m}([f']))$ and 
 $ j(i^{Y}_{\bB^2\bG_m}(d^Y_{\bB^2\bG_m}([f])))\simeq j([f])$, which is $j\circ\ell_0$ up to homotopy. By construction of $\ell$, we have a homotopy between the morphisms $j\circ \ell_0,i_Y^{\bB\sAut(\bB\bG_m)}\circ \ell:\Delta^1\to \Map(Y,\bB\sAut(\bB\bG_m))$. Then, using the two vertical equivalences in the cubic diagram, we obtain a compatible morphism from $\Delta^1$ to all the vertices of the cube except $\Map(Y',\bB^2\bG_m)$, extending $\ell$. By the universal property of the top horizontal face, we obtain $\widetilde{\ell}:\Delta^1\to \Map(Y',\bB^2\bG_m)$. The two vertices of the path $\widetilde{\ell}$ (i.e.\ precompositions of $\widetilde{\ell}$ by the inclusion of the vertices into $\Delta^1$) satisfy the same universal properties as $[f']$ and $d^Y_{\bB^2\bG_m}([f])$, showing that the two vertices are $[f']$ and $d^Y_{\bB^2\bG_m}([f])$ (up to homotopy). 
   \end{enumerate}
\end{proof}

\subsection{Coarse space}
Fix a nice commutative group scheme $G$ over $\Lambda$.
Let $f:X\to Y$ be a $\bB G$-principal $\infty$-bundle of scalloped derived stacks.
\begin{definition}
    We define the \textbf{coarse space morphism $\coarse_f:\Def_X\to \Def_Y$ associated with $f$} as the composition
    \[
    \Def_X\simeq \Def_f\to \Def_Y,
    \] where the equivalence comes from Proposition \ref{proposition-equiv-FMP-f-X-Gmgerbe-BGmbundle}, and  $\Def_f\to \Def_Y$ is the forgetful map.
\end{definition}
\begin{remark}
    This terminology is not standard, but it is motivated by the fact that $Y$ is the coarse space of $X$ if $Y$ is an algebraic space. In our general setup, $Y$ can have nontrivial stabilizers, so that it's not always the coarse space in the traditional sense.   
\end{remark}
\begin{proposition}\label{proposition-Omega-coarse-space-map}
For any perfect $V\in\Mod_{k,\geq0}$, 
    the morphism $\coarse_f:\Omega\Def_X(\sqz_k(V))\to \Omega\Def_Y(\sqz_k(V))$ is identified with
    \[
    \Map(\bL_{X/k},\cO_X\otimes_kV)\xrightarrow{(-)\circ[f]} \Map(f^*\bL_{Y/k},\cO_X\otimes_kV) \simeq  \Map(\bL_{Y/k},\cO_Y\otimes_kV)
    \] 
    up to the identification of Lemma \ref{lemma-loopspace-Def-X-cotangent-complex}, where the equivalence comes from Lemma \ref{lemma-compute-fiber-Def-f-for-Gm-gerbes}.
\end{proposition}
\begin{proof}
   Under the identifications of Proposition \ref{proposition-omega-Def-f} and Lemma \ref{lemma-loopspace-Def-X-cotangent-complex}, the morphisms $\Def_f\to\Def_X$ and $\Def_f\to \Def_Y$ correspond to the two projections of the fiber product
   \[
\Map(\bL_{X/k},\cO_X\otimes_kV)\times_{\Map(f^*\bL_{Y/k},\cO_X\otimes_kV)}\Map(\bL_{Y/k},\cO_Y\otimes_kV).
   \] The first projection is an equivalence as in Lemma \ref{lemma-compute-fiber-Def-f-for-Gm-gerbes} or Proposition \ref{proposition-equiv-FMP-f-X-Gmgerbe-BGmbundle}.
   The result follows immediately. 
\end{proof}
\begin{example}\label{example-etale-group-coarse-space}
    If $G_k$ is finite \'etale over $k$ (for example, $G=\mu_n$ such that $n$ is a unit in $k$), $\coarse_f:\Def_X\simeq \Def_Y$ is an equivalence, so that $\Def_X,\Def_Y,\Def_f,\Def_{f}^{\Bund(\bB G)},\Def_f^{\Gerbe(G)}$ are all equivalent by Proposition \ref{proposition-equiv-FMP-f-X-Gmgerbe-BGmbundle}.
To see this, we observe that $\bL_{X/Y}\simeq0$ which implies the map in Proposition \ref{proposition-Omega-coarse-space-map} is an equivalence. Passing to tangent fibres, we see $\coarse_f$ is an equivalence.  
\end{example}
\subsection{The ``change of group'' morphism}
The remaining two subsections are not needed for our main results later, but we include them to develop the basic derived deformation theory of gerbes further.

Let $\phi:G\to H$ be a morphism of commutative $\infty$-groups of derived stacks over some base $S$.
We will recall how to associate a $H$-gerbe to a $G$-gerbe using $\phi$, similarly to \cite[Construction 3.8]{bergh-BS}. 
Given a principal $\bB G$-bundle $f:X\to Y$ over $S$, we can associate a principal $\bB H$-bundle $f\wedge_GH:X\wedge_G H\to Y$. This is the quotient of the diagonal $\bB G$-action on $X\times_S \bB H$, where $\bB G$ acts via the given action on $X$ and acts through the inverse of $\bB \phi$ on $\bB H$.
This construction is compatible with base change (as colimits are universal in $\infty$-topoi). 
By construction, we have a morphism $\phi_X:X\to X\times_S\bB H\to X\wedge_GH$ over $Y$, where the first morphism uses the base point of $\bB H$. 

When $G,H$ are nice commutative group schemes over $\Lambda$, $Y$ is a scalloped derived stack over $k$ and $f:X\to Y$ is a principal $\bB G$-bundle, this argument  defines a ``change of group'' morphism 
\[(-)\wedge_GH:\Def^{\Bund(\bB G)}_f\to \Def^{\Bund(\bB H)}_{f\wedge_GH}.\] 
Alternatively, we can view this map $(-)\wedge_GH$ as the composition of the classifying morphism to $\bB^2G$ with the map $\bB^2\phi:\bB^2G\to\bB^2 H$ (which can be constructed by applying the above recipe to the universal principal $\bB G$-bundle).
\begin{proposition}\label{proposition-change-of-group-CS}
    For $G,H,f$ as above, the following diagram commutes:
    \[
    \xymatrix{
\Def_{f}^{\Bund(\bB G)}\ar[d]^-{(-)\wedge_GH}\ar[r]^-{\simeq}&\Def_X\ar[rd]^-{\coarse_{f}}& \\
 \Def_{f\wedge_GH}^{\Bund(\bB H)}\ar[r]^-{\simeq}&\Def_{X\wedge_GH} \ar[r]^{\coarse_{f\wedge_GH}}  &\Def_Y
    }
    \]
    where the equivalences come from Proposition \ref{proposition-equiv-FMP-f-X-Gmgerbe-BGmbundle}.
\end{proposition}
\begin{proof}
    This is clear from the construction since the map $(-)\wedge_GH$ preserves the targets of the principal bundles. 
\end{proof}
\begin{proposition}
For $G,H,f$ as above, via the identifications of Proposition \ref{proposition-Omega-Def-f-as-BGm-bundle}, consider 
the morphism of tangent fibres $\RHom(\bL_{Y/(\bB^2 G)_k},\cO_Y)[1]\to \RHom(\bL_{Y/(\bB^2 H)_k},\cO_Y)[1]$ induced by the 
``change of group'' morphism
    $(-)\wedge_GH:\Def^{\Bund(\bB G)}_f\to \Def^{\Bund(\bB H)}_{f\wedge_GH}$. This morphism is induced by the morphism $\bL_{Y/(\bB^2 H)_k}\to \bL_{Y/(\bB^2 G)_k}$ coming from the commutative diagram 
\[
\xymatrix{
Y\ar[r]^-{[f]}\ar[dr]_-{[f\wedge_GH]}&(\bB^2G)_k\ar[d]^-{(\bB^2\phi)_k}\\
&(\bB^2H)_k
}
\]
\end{proposition}
\begin{proof}
If we unravel the identification used in 
    Proposition \ref{proposition-Omega-Def-f-as-BGm-bundle}, this is immediate from the fact that $(-)\wedge_GH$ arises from the composition of  $\bB^2\phi$.
\end{proof}
\begin{example}
    Let $f:X\to Y$ be a $\mu_n$-gerbe of scalloped derived stacks over a field $k$. Then, we have a morphism
    \[(-)\wedge_{\mu_n}\bG_m:\Def^{\Bund(\bB\mu_n)}_f\to \Def^{\Bund(\bB \bG_m)}_{f\wedge_{\mu_n}\bG_m}.\] 
    By Example \ref{example-etale-group-coarse-space} and Proposition \ref{proposition-change-of-group-CS}, if $n$ is invertible in $k$, $\coarse_{f\wedge_{\mu_n}\bG_m}:\Def_{X\wedge_{\mu_n}\bG_m}\to \Def_Y$ admits a section. 
\end{example}
\begin{remark}
    Via the ``change of group'' construction along $\mu_n\to\bG_m$, one can define twisted derived categories for principal $\bB \mu_n$-bundles.
\end{remark}
\subsection{Deformations over augmented derived $k$-algebras}
In this subsection, we work with derived deformations over augmented derived $k$-algebras where $k$ is a field (i.e.\ we take $\Lambda=k$). 
In this situation, any formal moduli problem $F$ 
is pointed: 
given any test derived ring $k\to A\to k$, $F(A)$ has a distinguished point (up to homotopy) defined by the augmentation morphism $\{*\}\simeq F(k)\to F(A)$. This way, any formal moduli problem in this context has a morphism $\triv_F:\{*\}\to F$ from the trivial formal moduli problem.
\begin{lemma}\label{lemma-T-of-fiber}
    For any morphism $f:F\to G$ of formal moduli problems in $\CAlg_{k//k}^{\an,\art}$, we have a fiber sequence of tangent fibers: 
    \[\rT_{\fib_{\triv_G}(f)}\to \rT_F\to\rT_G,\]
    where $\fib_{\triv_G}(f)$ is the fiber of $f$ over the base point $\triv_G$.
\end{lemma}
\begin{proof}
By construction, for any perfect $V\in{\Dqccn}(k)$, we have a pullback diagram
\[
\xymatrix{
\fib_{\triv_G}(f)(\sqz_k(V))\ar[r]\ar[d]&\{*\}\ar[d]\\
F(\sqz_k(V))\ar[r]&G(\sqz_k(V))
}
\]
This immediately implies the result.
\end{proof}
Let $G$ be a nice commutative group scheme over $k$, and let $f:X\to Y$ be a principal $\bB G$-bundle of scalloped derived stacks over $k$.
\begin{definition}
    In the above situation, we define the formal moduli problem $\Def_X^{\gerby}:=\fib_{\triv_{\Def_Y}}(\coarse_f)$,  which parametrizes purely gerby derived deformations of $X$.
\end{definition}
\begin{proposition}
In the above situation, we have $\rT_{\Def_X^{\gerby}}\simeq \RHom(\bL_{X/Y},\cO_X)[1]$, and applying the tangent fibre functor to the fiber sequence $\Def_X^{\gerby}\to\Def_X\xrightarrow{\coarse_f}\Def_Y$ yields the following fiber sequence:
\[
\RHom_X(\bL_{X/Y},\cO_X)[1]\to
\RHom_X(\bL_{X/k},\cO_X)[1]\to 
\RHom_Y(\bL_{Y/k},\cO_Y)[1]
\]
induced by the transitivity fiber sequence of cotangent complexes:
\[
f^*\bL_{Y/k}\to\bL_{X/k}\to\bL_{X/Y}
\]
\end{proposition}
\begin{proof}
    We apply Proposition \ref{proposition-Omega-coarse-space-map} and Lemma \ref{lemma-T-of-fiber} and then invoke the uniqueness of fibers. 
\end{proof}
\begin{example}\leavevmode
Let $G$ and $f:X\to Y$ be as above. Let $n$ be a positive integer.
\begin{enumerate}
    \item When $G=\bG_m$, we obtain
  $\rT_{\Def_X^{\gerby}}\simeq \rR\Gamma(X, \cO_X)[2]$.
  \item
  Assume $\characteristic(k)=0$ or that $n$ is not divisible by $\characteristic(k)>0$. When $G=\mu_n$, we obtain
  $\rT_{\Def_X^{\gerby}}\simeq 0$.
  \item 
  Assume $n$ is divisible by $\characteristic(k)>0$. 
  When $G=\mu_{n}$, we obtain
  \[\rT_{\Def_X^{\gerby}}\simeq \rR\Gamma(X,\cO_X)[1]\oplus \rR\Gamma(X,\cO_X)[2].\]
\end{enumerate}
\end{example}

\section{Twisted-to-noncommutative map}
\subsection{Generalities on deforming categories}
Let $\Lambda$ be a complete local noetherian ring with residue field $k$.
Let $\cC$ be a $k$-linear presentable stable $\infty$-category, 
and we will recall the general framework for formal deformations of $\cC$. The theory is well-documented in \cite[\S16.6]{SAG} and \cite[\S2.3]{Moore2}, so we will keep the review brief. 
As in \cite[Variant 2.16, Remark 2.17]{Moore2} we can define a functor 
$\CatDef_\cC^+:\CAlg_{\Lambda//k}^\an\to\widehat{\cS}$
parametrizing deformations of $\cC$ over derived rings in $\CAlg_{\Lambda//k}^\an$. We also have a variant
$\CatDef_\cC:\CAlg_{\Lambda//k}^{\an,\art}\to\widehat{\cS}$ which is the restriction to $\CAlg_{\Lambda//k}^{\an,\art}\subset \CAlg_{\Lambda//k}^{\an}$.
It's not known if the functor $\CatDef_{\cC}$ is a formal moduli problem in general, but it is known to be  a $2$-proximate formal moduli problem. By a theorem of Lurie, we can approximate it by a formal moduli problem, which we write by $\CatDef^{\wedge}_{\cC}:\CAlg_{\Lambda//k}^{\an,\art}\to\cS.$ It comes with a natural map $\CatDef_{\cC}\to \CatDef^{\wedge}_{\cC}$, which is $0$-truncated (i.e\ evaluating the morphism at any derived ring in $\CAlganArtLambdaK$ yields a morphism of spaces with 0-truncated homotopy fibers).
In particular, applying $\Omega^2$ gives us an equivalence $\Omega^2\CatDef_{\cC}\simeq \Omega^2\CatDef^{\wedge}_{\cC}$.  
\subsection{Construction of the map}\label{subsection-construction-twisted-to-noncommutative-map}
We let $f:X\to Y$ be a principal $\bB\bG_m$-bundle over a scalloped derived stack $Y$ which is over $k$.
In this subsection, we will construct 
\textbf{the twisted-to-noncommutative map} of formal moduli problems
\[
\Theta^\wedge: \Def_X\to \CatDef_{\Dqc^1(X/Y)}^{\wedge}.
\]
In fact, we will construct a morphism of functors (which we also abusively call the \textbf{twisted-to-noncommutative map}) 
\[
\Theta: \Def_X\to \CatDef_{\Dqc^1(X/Y)}
\] which then induces the above map by composition with the obvious map $\CatDef_{\Dqc^1(X/Y)}\to \CatDef_{\Dqc^1(X/Y)}^{\wedge}.$

By Remark \ref{remark-relative-base-change-tw-Dqc}, the construction of the twisted derived category is compatible with base change, so we obtain a morphism
\[\Dqc^1(-):\Def_f^{\Bund(\bB\bG_m)}\to\CatDef_{\Dqc^1(X/Y)}
\] similarly to \cite[Definition 4.2]{Moore2}.  
\begin{definition}[twisted-to-noncommutative map] 
We define the \newline \textbf{twisted-to-noncommutative map} $\Theta^\wedge$ of $f:X\to Y$ as the composition
    \[
    \Theta^\wedge:\Def_X\simeq \Def_f^{\Bund(\bB\bG_m)}\xrightarrow[]{\Dqc^1}\CatDef_{\Dqc^1(X/Y)}\to \CatDef_{\Dqc^1(X/Y)}^{\wedge},
    \] where the first equivalence comes from Proposition \ref{proposition-equiv-FMP-f-X-Gmgerbe-BGmbundle}.
\end{definition}

\subsection{Computation of the induced map on the tangent fibres}
For any $V\in\Mod_{k,\geq0}$, we will write $X^V:=X[\cO_X\otimes_kV]$, and we similarly use the notation $Y^V$ as well.
First, we recall that for any $E\in\Dqccn(X),$ $\Map_X(\bL_{X/k},E)\simeq\Map'_k(X[E],X)$ functorially in $E$. As in Lemma \ref{lemma-loopspace-Def-X-cotangent-complex}, $\Omega\Def_X(\sqz_k(V))\simeq\Map(\bL_{X/k},\cO_X\otimes_kV)$, and the underlying $k$-module spectrum of $T_{\Def_X}$ is given by $\RHom(\bL_{X/k},\cO_X)[1]$.
In the sequel, we will use these facts freely. 

\begin{remark}
    We observe from the twisted Fourier-Mukai theory (Example \ref{example-twisted-FMT-via-dualizability}) that the Hochschild cohomology (complex) $\HH^*(\Dqc^1(X/Y)/k)\in\Mod_{k}$ of $\Dqc^1(X/Y)$ is computed as \[\RHom_{\Dqc^{-1,1}(X\times_kX/Y\times_kY)}(\iota^X_*\cP_X,\iota^X_*\cP_X).\]
\end{remark}
\begin{lemma}\label{lemma-loopspace-CatDef}
For any perfect   $V\in\Mod_{k,\geq0}$, we have an equivalence
    \[
    \Omega\CatDef_{\Dqc^1(X/Y)}(\sqz_k(V))\simeq \Map'_k(\Dqc^1(X/Y),\Dqc^1(X^V/Y^V)),
    \] which is natural in $V$, where $\Map'$ means that we trivialize the composition with the map $i^*:\Dqc^1(X^V/Y^V)\to \Dqc^1(X/Y)$ by $\id$. 
\end{lemma}
\begin{proof}
This is proved in the same way as \cite[Proposition 2.25]{Moore2} using the base change property of the twisted derived category (Remark \ref{remark-relative-base-change-tw-Dqc}).
\end{proof}
\begin{proposition}\label{proposition-CatDef-hat-tangent-fibre}
    For each perfect $V\in\Mod_{k,\geq0}$, $\CatDef_{\Dqc^1(X/Y)}^{\wedge}(\sqz_k(V))$ is equivalent to
\[ \hspace*{-1em}
    \Omega\Map_{\Fun_k^{\rL}(\Dqc^1(X/Y),\Dqc^1(X/Y))}(\id,(-)\otimes_kV[3]) 
   \simeq 
    \Omega\Map_{\Dqc^{-1,1}(X\times_k X/Y\times_k Y)}(\iota^X_*\cP_X,\iota^X_*\cP_X\otimes_kV[3])
    \]
     naturally in $V$. In particular, the underlying $k$-module spectrum of $T_{\CatDef^{\wedge}_{\Dqc^1(X/Y)}}$ is given by $\HH^*(\Dqc^1(X/Y)/k)[2]\in\Mod_{k}.$
\end{proposition}
\begin{proof}
By the axiom of formal moduli problem, we have 
\[
\CatDef_{\Dqc^1(X/Y)}^{\wedge}(\sqz_k(V))\simeq \Omega^2\CatDef_{\Dqc^1(X/Y)}^{\wedge}(\sqz_k(V[2])).
\]
The latter space is canonically equivalent to $\Omega^2\CatDef_{\Dqc^1(X/Y)}(\sqz_k(V[2]))$ since the map $\CatDef_{\Dqc^1(X/Y)}\to \CatDef_{\Dqc^1(X/Y)}^{\wedge}$ is $0$-truncated. 
By Lemma \ref{lemma-loopspace-CatDef} and the noncommutative twisted Atiyah map, we have 
\[
\Omega^2\CatDef_{\Dqc^1(X/Y)}(\sqz_k(V[2]))
\simeq \Omega\Map'_k(\Dqc^1(X/Y),\Dqc^1(X^{V[2]}/Y^{V[2]}))
    \]
    \[
    \simeq
    \Omega\Map_{\Fun_k^{\rL}(\Dqc^1(X/Y),\Dqc^1(X/Y))}(\id,(-)\otimes_kV[3]).
    \]
     Since $Y$ is scalloped, by Examples \ref{example-twisted-FMT-via-dualizability} and \ref{example-twisted-FMT-tensor}, we can write it using twisted Fourier-Mukai kernels as $\Omega\Map_{\Dqc^{-1,1}(X\times_k X/Y\times_k Y)}(\iota^X_*\cP_X,\iota^X_*\cP_X\otimes_kV[3]).$
\end{proof}

\begin{proposition}\label{proposition-map-on-tangent-fibers-twisted-to-noncommutative}
    Consider the twisted-to-noncommutative map $\Def_X\to \CatDef_{\Dqc^1(X/Y)}^{\wedge}$. The induced map on the tangent fibres makes the following diagram commute:
    \[
    \xymatrix{
T_{\Def_X}\ar[r]\ar[d]\ar[r]^-{\simeq}&\RHom(\bL_{X/k},\cO_X)[1]\ar[d]^-{\widetilde{\bAt}^{\FM}_{\cO_X[1]}}\\
    T_{\CatDef_{\Dqc^1(X/Y)}^{\wedge}}\ar[r]^-{\simeq}&\HH^*(\Dqc^1(X/Y)/k)[2]
    }
    \]
    Here, the top horizontal map was explained at the beginning of the subsection, and the bottom one is Proposition \ref{proposition-CatDef-hat-tangent-fibre}. The right vertical map is the spectrally enriched twisted Atiyah map (Definition \ref{definition-spectrally-enriched-twisted-atiyah-map}).
    \end{proposition}
\begin{proof}
For each perfect $V\in\Mod_{k,\geq0},$ we consider the following commutative diagram:
\[
\makebox[\textwidth][c]{
\xymatrix{
\Def_X(\sqz_k(V))\ar[r]^-{\simeq}\ar[d]&
\Omega^2 \Def_X(\sqz_k(V[2]))\ar[r]^-{\simeq}\ar[d]&
\Omega\Map(\bL_{X/k},\cO_X\otimes V[2])\ar[d]\\
\CatDef_{\Dqc^1(X/Y)}(\sqz_k(V))\ar[r]\ar[d]&
\Omega^2\CatDef_{\Dqc^1(X/Y)}(\sqz_k(V[2]))\ar[r]^-{\simeq}\ar[d]^-{\simeq}&
\Omega\Map_{\Dqc^{-1,1}(X\times_k X/Y\times_k Y)}(\iota^X_*\cP_X,\iota^X_*\cP_X\otimes_kV[3])\ar@{=}[d]\\
\CatDef_{\Dqc^1(X/Y)}^{\wedge}(\sqz_k(V))\ar[r]^-{\simeq}&
\Omega^2\CatDef_{\Dqc^1(X/Y)}^{\wedge}(\sqz_k(V[2]))\ar[r]^-{\simeq}&
\Omega\Map_{\Dqc^{-1,1}(X\times_k X/Y\times_k Y)}(\iota^X_*\cP_X,\iota^X_*\cP_X\otimes_kV[3])
}
}
\]
The horizontal maps from the first column to the second column arise from the functoriality of $\Def_X,\CatDef_{\Dqc^1(X/Y)},\CatDef_{\Dqc^1(X/Y)}^{\wedge}$, and the top one and the bottom one are equivalences by the axiom of formal moduli problems. 
The equivalence $\Omega^2\CatDef_{\Dqc^1(X/Y)}(\sqz_k(V[2]))\simeq \Omega^2\CatDef_{\Dqc^1(X/Y)}^{\wedge}(\sqz_k(V[2]))$
is due to the $0$-truncatedness of the map $\CatDef_{\Dqc^1(X/Y)}\to \CatDef_{\Dqc^1(X/Y)}^{\wedge}$. 
The top right horizontal equivalence is as in the beginning of the subsection,
and the middle right horizontal equivalence comes from the combination of the noncommutative twisted Atiyah map and Lemma \ref{lemma-loopspace-CatDef}.
The top right vertical morphism is filled in to make the diagram commute.

We only need to check that this vertical map can be chosen to be $\Omega\widetilde{\bAt}_{\cO_X\otimes_k V[2]}$. If we unpack the constructions of the maps, we see that this is (the loopspace version of) the content of Proposition \ref{proposition-commutative-vs-noncommutative-atiyah-maps}.

\end{proof}
\begin{theorem}\label{theorem-tangent-fibre-trivial-gerbe-hkr}
    Let $Y$ be a smooth proper variety over a field $k$, and let $f:X=Y\times\bB\bG_m\to Y$ be the trivial principal $\bB\bG_m$-bundle. 
    Let $\phi:T_{\Def_X}\to T_{\CatDef_{\Dqc^1(X/Y)}^{\wedge}}\simeq T_{\CatDef_{\Dqc(Y)}^{\wedge}} $ be the morphism of tangent fibres induced by the twisted-to-noncommutative map. Then,
\begin{enumerate}
    \item the morphism $\phi$ can be identified with 
     \[I^{\HKR}_{\leq1}[2]:
\RHom_Y(\bL_{Y/k},\cO_Y[1])
\oplus \RHom_Y(\cO_Y,\cO_Y[2])
\to 
\RHom_{Y\times_k Y}(\Delta_*\cO_Y,\Delta_*\cO_Y[2]),\]
\item the morphism $\rH^{-1}(\phi)$ is an isomorphism,
  \item the morphism $\rH^0(\phi)$ is an isomorphism if $\rH^0(Y,\bigwedge^2 \rT_Y)=0$, and
  \item the morphism $\rH^1(\phi)$ is injective if $\characteristic(k)\neq 2$ or $\rH^0(Y,\bigwedge^2 \rT_Y)=0.$
\end{enumerate}
\end{theorem}
\begin{proof}
Part (1) follows from Proposition \ref{proposition-map-on-tangent-fibers-twisted-to-noncommutative} and
Theorem \ref{theorem-spectral-variant-twisted-atiyah-trivial-gerbe-formula-HKR}.
The remaining parts follow from 
    Proposition \ref{proposition-injectivity-truncated-HKR}.
\end{proof}
By the following lemma, we do not need to distinguish lifts of scalloped derived stacks as 1-Artin derived stacks and lifts as scalloped derived stacks.
\begin{lemma}\label{lemma-lift-scalloped-qcqs-property-preserved}
Fix $A\in\CAlganArtLambdaK$. Let $X$ be a qcqs scalloped 1-Artin derived stack over $k$. If $\cX$ is a deformation of $X$ over $A$ (i.e.\ it's a 1-Artin derived stack with  $\cX\times_Ak\simeq X$), then $\cX$ is qcqs and scalloped.
\end{lemma}
\begin{proof}
 We have $\cX_{\cl}\times_{\pi_0(A)}k\simeq X_{\cl}$ as classical locus commutes with fiber products. $\pi_0(A)$ is a nilpotent thickening of $k$ (note that $\pi_0(A)\to k$ is surjective as its composition with $\Lambda\to \pi_0(A)$ is the quotient map), so $\cX_{\cl}$ is also a nilpotent thickening of $X_{\cl}$ (although it's not a flat deformation in general). By Lemma \ref{lemma-qc-qs-check-on-classical-locus}, $X_{\cl}$ is qcqs, but being qcqs is invariant under infinitesimal thickenings (for classical stacks), so $\cX_{\cl}$ is also qcqs. Lemma \ref{lemma-qc-qs-check-on-classical-locus} then shows that $\cX$ is qcqs. 
 In view of Theorem \ref{theorem-khan-ravi-scalloped-iff-nicely-stabilized}, we will check $\cX$ is scalloped by studying the stabilizers of $\cX_{\cl}$. As above, $\cX_{\cl}$ is an infinitesimal thickening of $X_{\cl},$ so  $\cX_{\cl}$ and  $X_{\cl}$ have the same set of  geometric points. The stabilizers at these geometric points are identified along the closed immersion $X_{\cl}\to\cX_{\cl}$ as in the proof of Proposition \ref{proposition-existence-geometric-pushout}, so the stabilizers of $\cX_{\cl}$ are nice, as $X$ is scalloped.
\end{proof}
\begin{remark}
    Recall that qcqs 1-Artin derived stacks automatically have affine diagonals in our convention. Even if we adopt a different convention and allow potentially non-affine diagonals, when $X$ as in the lemma has an affine diagonal, the lift $\cX$ automatically has an affine diagonal as well (so the space of deformations is not made artificially smaller by our choice of convention).

\end{remark}
\begin{theorem}[{see \cite[Theorem B(1)]{Moore2}} for an untwisted analogue]\label{theorem-sqz-lift-twisted-Dqc}
    Let $Y$ be a smooth proper $k$-scheme. Assume $\characteristic(k)\neq2$ or $\rH^0(Y,\bigwedge^2\rT_Y)=0.$ Let $A'\to A$ be a {small} extension of classical Artinian local $\Lambda$-algebras with residue field $k.$
    Suppose $X:=Y\times\bB\bG_m$ lifts to $\cX$ over $A$. 
    Then, the principal $\bB\bG_m$-bundle $X\to Y$ lifts uniquely to $\cX\to \cY$ over $A$ up to equivalence. Moreover, 
    $\cX$ lifts to $A'$ if and only if $\Dqc^1(\cX/\cY)$ lifts to $A'$. 
\end{theorem}
\begin{proof}
Relevant lifts of derived stacks are automatically scalloped by Lemma \ref{lemma-lift-scalloped-qcqs-property-preserved}, so
the lifting of the principal bundle follows from Proposition \ref{proposition-equiv-FMP-f-X-Gmgerbe-BGmbundle}. 
The other statement is an immediate consequence of Theorem \ref{theorem-tangent-fibre-trivial-gerbe-hkr} (4) by standard derived deformation theory; this shows the obstruction class for lifting $\cX$ to $A'$ maps injectively to the obstruction class for lifting $\Dqc^1(\cX/\cY)$ to $A'$ (see \cite[Definition 3.12]{BranTael}).
\end{proof}

\begin{theorem}
    Consider the trivial principal $\bB\bG_m$-bundle $f:X=Y\times\bB\bG_m\to Y$ over a smooth proper variety $Y$ over $k$. Let $A\to k$ be a square-zero thickening as a classical $\Lambda$-algebra. 
    Then,  
 the morphism $\pi_0(\Theta(A)):\pi_0(\Def_X(A))\to \pi_0(\CatDef_{\Dqc^1(X/Y)}(A))\simeq \pi_0(\CatDef_{\Dqc(Y)}(A))$ is injective. In other words, if the twisted derived categories of two lifts of $X$ over $A$ admit an $A$-linear equivalence reducing to the identity functor of $\Dqc^1(X/Y)$, then these lifts of $X$ are equivalent. 
\end{theorem}
\begin{proof}
{We assume that $\pi_0(\Def_X(A))$ is nonempty, since the claim is trivial otherwise.
Let $I$ be the kernel of the thickening $A\to k$.
Then, we have a pullback square
\[
\xymatrix{
A\ar[r]\ar[d]&k\ar[d]\\
k\ar[r]&\sqz_k(I[1])
}
\]
of derived rings, and using the definition of formal moduli problems,
we have the following pullback square:
\[
\xymatrix{
\Def_X(A)\ar[r]\ar[d]&\Def_X(k)\ar[d]\\
\Def_X(k)\ar[r]&\Def_X(\sqz_k(I[1]))
}
\]
We have $\Def_X(k)\simeq\{*\}$, so the above diagram gives us  two points of $\Def_X(\sqz_k(I[1]))$ (where one corresponds to the trivial deformation of $X$ over $\sqz_k(I[1])$), and  
and if we choose a point of $\Def_X(A)$, this determines a path between them. Using this path, we can identify the pullback $\Def_X(A)$ with $\Omega\Def_X(\sqz
_k(I[1]))$, where the loopspace is based at the trivial deformation of $\Def_X(\sqz
_k(I[1]))$. 
Thus, after fixing a point of $\Def_X(A)$ (which accordingly determines a point of $\CatDef^{\wedge}_{\Dqc^1(X/Y)}(A)$), 
}
we have the following commutative diagram with horizontal equivalences 
\[
\xymatrix{
\Def_X(A)\ar[r]^-{\simeq}\ar[d]&\Omega\Def_X(\sqz_k(I[1]))\ar[d]\\
\CatDef_{\Dqc^1(X/Y)}^{\wedge}(A)\ar[r]^-\simeq& \Omega\CatDef_{\Dqc^1(X/Y)}^{\wedge}(\sqz_k(I[1]))
}
\]
The vertical maps come from the twisted-to-noncommutative map $\Theta^\wedge$. Taking $\pi_0$ of the right vertical map yields $\rH^0(I^{\HKR}_{\leq1}\otimes_k I[2])$ by Theorem \ref{theorem-tangent-fibre-trivial-gerbe-hkr}. This map is injective on $\pi_0$ by Lemma \ref{lemma-injectivity-truncated-HKR}, so the left vertical map is also injective on $\pi_0$. If we factor this map on $\pi_0$ as
\[
\pi_0(\Def_X(A))\to \pi_0(\CatDef_{\Dqc^1(X/Y)}(A))\to \pi_0( \CatDef_{\Dqc^1(X/Y)}^{\wedge}(A)),
\] the first map in the composition must also be injective.
\end{proof}
\begin{corollary}\label{corollary-unique-sqz-lift-FM-partner}
    Let $X,Y$ be smooth proper varieties over $k$, and assume we have a Fourier-Mukai equivalence $\Phi:\Dqc(X)\simeq\Dqc(Y)$. Let $A\to k$ be a square-zero thickening as a classical $\Lambda$-algebra.
    Let $\cX$ be a deformation of $X\times\bB\bG_m$ over $A$ (which comes with a principal $\bB\bG_m$-bundle structure $\cX\to \cX_{\coarse}$ lifting the projection $X\times\bB\bG_m\to X$). Then, there is at most one deformation $\cY$ of $Y\times\bB\bG_m$ over $A$ up to equivalence (which comes with a principal $\bB\bG_m$-bundle structure $\cY\to \cY_{\coarse}$ lifting the projection $Y\times\bB\bG_m\to Y$) which admits an $A$-linear Fourier-Mukai equivalence $\Dqc^1(\cX/\cX_{\coarse})\simeq \Dqc^1(\cY/\cY_{\coarse})$ reducing to 
    \[\Phi:\Dqc^1(X\times\bB\bG_m)\simeq\Dqc(X)\simeq \Dqc(Y)\simeq\Dqc^1(Y\times\bB\bG_m).\]
    In particular, for any (untwisted) lift $\cX$ of $X$ over $A$, there is at most one (untwisted) lift $\cY$ of $Y$ up to equivalence which admits an $A$-linear Fourier-Mukai equivalence $\Dqc(\cX)\simeq\Dqc(\cY)$ reducing to $\Phi:\Dqc(X)\simeq\Dqc(Y)$.
\end{corollary}
\begin{proof}
    For the first claim, suppose we have two such lifts $\cY^1\to \cY^1_{\coarse}$ and $\cY^2\to \cY^2_{\coarse}$ over $A$ which admit $A$-linear Fourier-Mukai equivalences
    \[
    \Phi^1:\Dqc^1(\cX/\cX_{\coarse})\simeq \Dqc^1(\cY^1/\cY^1_{\coarse}) ,\quad \Phi^2:\Dqc^1(\cX/\cX_{\coarse})\simeq \Dqc^1(\cY^2/\cY^2_{\coarse})
    \] reducing to $\Phi.$
    The reduction of $\Phi^2\circ (\Phi^1)^{-1}$ to $k$ is equivalent to the identity functor, so the theorem applies. 

    The untwisted statement follows from the twisted one since any untwisted lift $\cX$ of $X$ over $A$ yields the lift $\cX\times\bB\bG_m$ of $X\times\bB\bG_m$ over $A$, and $\Dqc^1(\cX\times\bB\bG_m/\cX)\simeq\Dqc(\cX)$.
\end{proof}
\begin{theorem}[{see \cite[Corollary 1.4]{Moore2}} for an untwisted analogue]\label{theorem-equivalence-classical-deformation-functors}
    \leavevmode
    \begin{enumerate}
        \item  Let $X$ be a smooth proper $k$-scheme with $\rH^0(X,\bigwedge^2\rT_X)=0.$ Then, the twisted-to-noncommutative map induces an isomorphism of classical deformation functors:
        \[
        \Def_{X\times\bB\bG_m}^{\cl}\simeq \CatDef_{\Dqc^1(X\times\bB\bG_m/X)}^{\wedge,\cl}\simeq \CatDef_{\Dqc(X)}^{\wedge,\cl}
        \] 
        \item When $X,Y$ are two such $k$-schemes, any (untwisted) derived equivalence between $X$ and $Y$ induces an isomorphism of classical deformation functors:
        \[
        \Def_{X\times\bB\bG_m}^{\cl}\simeq \Def_{Y\times\bB\bG_m}^{\cl}
        \] 
        \item Let $X,Y$ be smooth proper $k$-schemes that are derived equivalent, and assume that $\rH^0(Y,\bigwedge^2\rT_Y)=0.$ 
        Let $A$ be an artin local (classical) $\Lambda$-algebra with a map $A\to k$.
        If $X$ lifts to $A$ (or more generally if $X\times\bB\bG_m$ lifts to $A$), $Y$ also lifts to $A.$
    \end{enumerate}
\end{theorem}
\begin{proof}
The proofs of (1),(2) are similar to \cite[Corollary 1.4]{Moore2}, and the inputs needed to run the same argument are the surjectivity of $\rH^{-1}(\phi)$, the isomorphy of $\rH^0(\phi)$, and the injectivity of $\rH^1(\phi)$ from Theorem \ref{theorem-tangent-fibre-trivial-gerbe-hkr}. 

For (3), assume that $X$ lifts to some $\cX$ over $A$. Then, $\cX\times\bB\bG_m$ is a lift of $X\times\bB\bG_m$ to $A$ (and is qcqs scalloped by Lemma \ref{lemma-lift-scalloped-qcqs-property-preserved}), so
we obtain an element of $\CatDef_{\Dqc^1(X\times\bB\bG_m/X)}^{\wedge,\cl}(A)\simeq \CatDef_{\Dqc^1(Y\times\bB\bG_m/Y)}^{\wedge,\cl}(A)$. 
By part (1), we obtain a lift $\cY'$ of $Y\times\bB\bG_m$ to $A$. As in Proposition \ref{proposition-equiv-FMP-f-X-Gmgerbe-BGmbundle}, this 
lift comes with a morphism $\cY'\to\cY$ lifting the projection $Y\times\bB\bG_m\to Y$, so in particular, we have a lift $\cY$ of $Y$.
\end{proof}
\begin{remark}
In \S\ref{section-compactly-generated-variant}, we improve Theorem \ref{theorem-equivalence-classical-deformation-functors}(1) by considering compactly generated derived deformations of $\infty$-categories.
\end{remark}

\section{Compactly generated variant}\label{section-compactly-generated-variant}
In this section, we discuss a variant of $\CatDef_{\cC}$ where we impose compact generation, and we will use it to prove an additional result about deformations of twisted derived categories.

Let $\Lambda$ be a complete local noetherian ring with residue field $k$, and let $\cC$ be a compactly generated $k$-linear presentable stable $\infty$-category. Then, as in \cite[Variant 2.18]{Moore2}, we can define a functor
\[
\CatDef_{\cC}^\omega:\CAlganArtLambdaK\to \cS,
\] which is a subfunctor of $\CatDef_{\cC}$ parametrizing compactly generated deformations of $\cC.$
This variant is better-behaved than $\CatDef_{\cC}$ for the following reason:
\begin{proposition}[Lurie, {see \cite[Theorem 2.19]{Moore2}}]\label{proposition-Lurie-1-proximate-tamely-compactly-generated}
Assume $\cC$ is tamely compactly generated. Then, 
    $\CatDef_{\cC}^\omega$ is a 1-proximate formal moduli problem. 
\end{proposition}
In the sequel, assume that $\cC$ is tamely compactly generated. The formal moduli problem approximation of $\CatDef_{\cC}^\omega$ in fact coincides with the approximation of $\CatDef_{\cC}$ (one can see this from comparing tangent fibers for instance), 
so we can view the composition 
\[
\CatDef^\omega_{\cC}\to \CatDef_{\cC}\to \CatDef_{\cC}^{\wedge}
\] as the approximation morphism.
This composition morphism is $(-1)$-truncated (\cite[Theorem 16.4.2.1]{SAG}), and in particular evaluating it at test derived rings in $\CAlganArtLambdaK$ induces injections on $\pi_0$. 

When $f:X\to Y$ is a principal $\bB\bG_m$-bundle of scalloped derived stacks over $k$, 
we can use this construction to define a variant of the twisted to noncommutative map from \S\ref{subsection-construction-twisted-to-noncommutative-map}. 
Indeed, by Proposition \ref{proposition-tamely-compactly-generated}, the morphism $\Theta:\Def_X\to \CatDef_{\Dqc^1(X/Y)}$ factors through $\Theta^\omega:\Def_X\to \CatDef_{\Dqc^1(X/Y)}^\omega.$

\begin{theorem}[{Refinement of Theorem \ref{theorem-equivalence-classical-deformation-functors}(1)}]\label{theorem-compactly-generated-deformation-functor}
    Consider the trivial principal $\bB\bG_m$-bundle $f:X=Y\times\bB\bG_m\to Y$ over a smooth proper variety $Y$ over $k$. Assume $\rH^0(Y,\bigwedge^2 T_Y)=0.$
    
Let $\CatDef_{\Dqc(Y)}^{\omega,\cl}$ denote the composition
    \[
\CRing_{\Lambda//k}^{\art}\hookrightarrow\CAlganArtLambdaK\xrightarrow{\CatDef_{\Dqc(Y)}^{\omega}}{\cS}\xrightarrow{\pi_0}{\mathrm{Set}}
    \] where $\CRing_{\Lambda//k}^{\art}\subset \CAlganArtLambdaK$ is the full subcategory of classical rings.
Then, the morphism $\Def^{\cl}_{X}\to \CatDef_{\Dqc(Y)}^{\omega,\cl}$ induced by $\Theta^\omega$ is an isomorphism. In particular, $\CatDef_{\Dqc(Y)}^{\omega,\cl}$ is a deformation functor. 
\end{theorem}
\begin{proof}
Consider the commutative diagram of functors.
\[
\xymatrix{
\Def^{\cl}_{X}\ar[r]\ar[dr]&\CatDef_{\Dqc^1(X/Y)}^{\omega,\cl}\ar[d]\\
&\CatDef_{\Dqc^1(X/Y)}^{\wedge,\cl}
}
\]
The diagonal morphism is an isomorphism by
Theorem \ref{theorem-equivalence-classical-deformation-functors}(1). This implies in particular that the vertical map is surjective (if we evaluate at any test ring). 
The vertical morphism is also injective since the morphism $\CatDef_{\Dqc^1(X/Y)}^{\omega}\to \CatDef_{\Dqc^1(X/Y)}^{\wedge}$ is $(-1)$-truncated. Thus, all maps are isomorphisms. 
\end{proof}
\begin{remark}\label{remark-deformation-functor-CatDef-question}
     This theorem shows that the ``classicalization'' of the functor
     $\CatDef_{\Dqc(Y)}^{\omega}$ is a deformation functor, and one can ask for the derived analogue of this statement; namely, is
     $\CatDef_{\Dqc(Y)}^{\omega}$ a formal moduli problem? To the author's knowledge, it is only known to be a 1-proximate formal moduli problem (Proposition \ref{proposition-Lurie-1-proximate-tamely-compactly-generated}), but there appears to be no known counter-example to this derived statement.
\end{remark}

\addtocontents{toc}{\vspace{\normalbaselineskip}}
\bibliographystyle{amsalpha}
\bibliography{bib}

\end{document}